\mag=1200
\input amstex
\documentstyle{amsppt}
\NoBlackBoxes
\tolerance=1000



\input xypic

\def\wtilde{\mathaccent"0365 }

\define\emph#1{{\it #1\/}}

\def\lw#1#2{\lower#1mm\hbox{#2}}
\define\wwt{{w}}
\def\bopl{\bigoplus}
\def\vp{\vphantom{$a^b$}}
\def\kkk{\kern.2mm}

\define\tabboxone#1{\vtop{\hsize3.2cm
\tolerance10000
\parindent0pt
\baselineskip15mm
\lineskip5pt
\par
#1\par\null}}

\define\tabboxtwo#1{\vtop{\hsize2.2cm
\tolerance10000
\parindent0pt
\baselineskip15mm
\lineskip5pt
\par
#1\par\null}}

\define\tabboxthree#1{\vtop{\hsize1.1cm
\tolerance10000
\parindent0pt
\baselineskip4.5mm
\lineskip5pt
\par
#1\par\null}}

\define\PP{{\Bbb P}}
\define\FF{{\Bbb F}}
\define\CC{{\Bbb C}}
\define\wt{\operatorname{wt}}
\define\Proj{\operatorname{Proj}}

\topmatter

\title
Hyperelliptic and trigonal Fano threefolds
\endtitle

\author
V. V. Przyjalkowski, I. A. Cheltsov, and K. A. Shramov
\endauthor

\address
V.~A.~Steklov Math. Institute RAS
\endaddress

\footnote""{This work was partially supported by the Russian Foundation for
Basic Research (grant no.~04-01-00613).}

\email
przhijal\@mccme.ru,
\hfil\break\indent
\phantom{\it E-mail address: }cheltsov\@yahoo.com,
\hfil\break\indent
\phantom{\it E-mail address: }shramov\@mccme.ru
\endemail


\abstract
We classify Fano 3-folds with canonical Gorenstein singularities
whose anticanonical linear system has no base points but does not give an
embedding, and we classify anticanonically embedded Fano 3-folds with
canonical Gorenstein singularities which are not intersections of
quadrics. We also study the rationality questions for most of these
varieties.
\endabstract

\date
01/JUL/04
\enddate

\endtopmatter

\rightheadtext{Hyperelliptic and trigonal Fano threefolds}
\leftheadtext{V. V. Przyjalkowski, I. A. Cheltsov, and K. A. Shramov}

\document

\advance\baselineskip by 0.4pt

\head
\S\,1. Introduction
\endhead

Consider a Fano threefold $X$ with canonical Gorenstein
singularities\,\footnote{Canonical Gorenstein singularities are exactly the
rational Gorenstein singularities (see~\cite{70}).} (see~\cite{45},
\cite{112}, \cite{85}, \cite{39},~\cite{40}). Suppose that the
anticanonical linear system $|-K_X|$ is base point free. It is well
known that such varieties are divided into three classes.

1) Hyperelliptic varieties (that is, the morphism $\varphi_{|-K_X|}$ is
not an embedding). Then the intersection of two general divisors in
$|-K_X|$ is a hyperelliptic curve.

2) Trigonal varieties (that is, the morphism $\varphi_{|-K_X|}$ is an
embedding but its image is not an intersection of quadrics). Then the
intersection of two general divisors in~$|-K_X|$ is a trigonal curve or
the canonical image of a smooth plane quintic.

3) Varieties whose image under the embedding $\varphi_{|-K_X|}$ is an
intersection of quadrics.

We study varieties of the first two types. Theorems~1.5 and~1.6 give
complete classifications of hyperelliptic and trigonal varieties
respectively. Proposition~1.10 establishes rationality or non-rationality
for most of these varieties.

In the introduction we survey the modern state of the classification
problem of Fano threefolds with canonical Gorenstein singularities,
including the main results of this paper (Theorems~1.5 and~1.6 and
Proposition~1.10). The second section contains various known results
that are used in the proofs. In \S\S\,3 and~4 we prove Theorems~1.5 
and~1.6 respectively. In~\S\,5 we study the rationality
questions for elliptic and trigonal Fano threefolds.

The biregular classification of 3-folds whose curve sections are canonical
curves was considered by Fano~\cite{72}--\cite{75}. In the smooth case,
hyperplane sections of such 3-folds must be K3 surfaces by the adjunction
formula. Hence a natural generalization of the problem studied by Fano
is the biregular classification of 3-folds containing an ample effective
Cartier divisor which is a K3 surface with at most Du~Val singularities.
It turns out that, except for generalized cones (usual cones in the very
ample case) over K3 surfaces, all 3-folds of this class are Fano 3-folds
with canonical Gorenstein singularities (see \cite{45}, \cite{112},
\cite{85}, \cite{39},~\cite{40}).

A complete classification of smooth Fano 3-folds was obtained in
\cite{46}, \cite{12}, \cite{13},
\cite{104}, \cite{108}, \cite{107}, \cite{105}, where $105$ families of
smooth Fano 3-folds were found (see~\cite{88}). Moreover, every Fano
3-fold with terminal Gorenstein singularities is a deformation of a smooth
one (see~\cite{109}). However, there are Fano 3-folds with canonical
Gorenstein singularities that cannot be deformed into smooth Fano 3-folds.
For example, the weighted projective spaces
$\Bbb{P}(1^{3},3)$ and $\Bbb{P}(1^{2},4,6)$ are Fano 3-folds with
canonical Gorenstein singularities (see~\cite{68}, \cite{76})
but cannot be globally deformed into smooth varieties.

The classification of Fano 3-folds with canonical Gorenstein singularities
is nowadays far from being complete (see~\cite{107}). However, 4 important
steps are already done. The first step is the following result proved in
\cite{45} and~\cite{112} (see also \cite{48} and~\cite{49}).

\proclaim{Theorem 1.1}
Let $X$ be a Fano 3-fold with canonical Gorenstein singularities.
In particular, the anticanonical divisor $-K_{X}$ is an ample Cartier
divisor. Let $S$ be a sufficiently general surface in the complete linear
system $|-K_{X}|$. Then $S$ has at most Du~Val singularities.
\endproclaim

The second step is the following result of~\cite{107}. It is a natural
generalization of the classification of smooth Fano 3-folds with Picard
group $\Bbb{Z}$, which was obtained in~\cite{12} and~\cite{13}.

\proclaim{Theorem 1.2}
Let $X$ be a Fano $3$-fold with canonical Gorenstein singularities.
Suppose that the linear system $|-K_{X}|$ has no movable decomposition,
that is, the anticanonical divisor $-K_{X}$ is not rationally equivalent
to $A+B$ where $A$,~$B$ are Weil divisors whose complete linear systems
$|A|$,~$|B|$ have positive dimension. Then $X$ is one of the following
$3$-folds.

\rom{1)} A hypersurface of degree $6$ in~$\Bbb{P}(1^{4},3)$,
\ $-K_{X}^{3}=2$.

\rom{2)} A complete intersection of a quadric cone and a quartic
hypersurface in~$\Bbb{P}(1^{5},2)$,
$-K_{X}^{3}=\nomathbreak 4$.

\rom{3)} A quartic hypersurface in~$\Bbb{P}^{4}$, \ $-K_{X}^{3}=4$.

\rom{4)} A complete intersection of a quadric and a cubic in
$\Bbb{P}^{5}$, \ $-K_{X}^{3}=6$.

\rom{5)} A complete intersection of three quadrics in
$\Bbb{P}^{6}$, \ $-K_{X}^{3}=8$.

\rom{6)} An intersection of the Grassmannian $G(1,4)\subset\Bbb{P}^{9}$
with a linear subspace of codimension~$2$ and a quadric, $-K_{X}^{3}=10$.

\rom{7)} An intersection of the orthogonal Grassmannian
$OG(5,10)\subset\Bbb{P}^{15}$ with a linear subspace of
codimension~$7$, \ $-K_{X}^{3}=12$.

\rom{8)} An intersection of the Grassmannian $G(2,6)\subset\Bbb{P}^{14}$
with a linear subspace of codimension~$5$, \ $-K_{X}^{3}=14$.

\rom{9)} An intersection of the symplectic Grassmannian
$LG(3,6)\subset\Bbb{P}^{13}$ with a linear subspace of codimension
$3$, \ $-K_{X}^{3}=16$.

\rom{10)} An intersection of the $G_2$-homogeneous space
$\Sigma\subset\Bbb{P}^{13}$ with a linear subspace of codimension~$2$
\rom{(}see~\cite{88}, Example~$5.2.2$\rom{)}, \ $-K_{X}^{3}=18$.

\rom{11)} A Mukai--Umemura $3$-fold $V_{22}\subset\Bbb{P}^{13}$
\rom{(}see, for example,~\cite{108}, \cite{80},~\cite{88}\rom{)},
$-K_{X}^{3}=22$.
\endproclaim

The third step is the following boundedness result of~\cite{110}.

\proclaim{Theorem 1.3}
Let $X$ be a Fano $3$-fold with canonical Gorenstein singularities.
In particular, the anticanonical divisor $-K_{X}$ is an ample Cartier
divisor. Then $-K_{X}^{3}\leqslant 72$, and the equality implies that
either $X\cong\Bbb{P}(1^{3},3)$ or $X\cong\Bbb{P}(1^{2},4,6)$.
\endproclaim

The fourth step is the following result proved in~\cite{84}.

\proclaim{Theorem 1.4}
Let $X$ be a Fano $3$-fold with canonical Gorenstein singularities.
Suppose that the base locus of~$|-K_{X}|$ is non-empty. Then $X$ is one
of the following $3$-folds.

\rom{($B_{1}$)} A complete intersection of a quadric cone and a sextic in
$\Bbb{P}(1^{4},2,3)$,\linebreak  
$-K_{X}^{3}=2$.

\rom{($B_{2}$)} The blow up of a sextic in~$\Bbb{P}(1^{3},2,3)$ along a
curve of arithmetic genus~$1$, $-K_{X}^{3}=4$.

\rom{($B_{3}$)} $S_{1}\times\Bbb{P}^{1}$, where $S_{1}$ is a del~Pezzo
surface of degree~$1$ with Du~Val singularities, $-K_{X}^{3}=6$.

\rom{($B_{4}^m$)} The anticanonical model of the blow up of
$U_m$ along a curve~$\Gamma_0$, where $U_m$ is a double covering
$\pi\:U_m\!\to\Proj\bigl(\!\Cal{O}_{\Bbb{P}^1}(m)\oplus
\Cal{O}_{\Bbb{P}^1}(m-4)\oplus\nomathbreak \Cal{O}_{\Bbb{P}^1}\bigr)\!
=\Bbb{F}(m, m-4, 0)$ such that $-K_{U_m}=\pi^*M$ and $U_m$ has at worst
canonical singularities, and $\Gamma_0$ is a smooth rational complete
intersection contained in the smooth part of~$U_m$ such that
$\pi(\Gamma_0)$ is a complete intersection of a general divisor in~$|M|$
and the (unique) divisor in the linear system $|M-mF|$. Here $M$ is the
class of the tautological sheaf on~$\Bbb{F}(m, m-4, 0)$, and $F$ is the
class of a fibre of the natural projection to~$\Bbb{P}^1$. We also have
$3\leqslant m\leqslant 12$ and $-K_{X}^3=2m-2$.
\endproclaim

The purpose of this paper is to prove the following two results.

\proclaim{Theorem 1.5}
Let $X$ be a Fano $3$-fold with canonical Gorenstein singularities.
Suppose that the linear system $|-K_{X}|$ has no base points but the
induced morphism $\varphi_{|-K_{X}|}$ is not an embedding. Then
$X$ is one of the following $47$ Fano $3$-folds.

\rom{1) ($H_{1}$)} A hypersurface of degree $6$ in~$\Bbb{P}(1^{3},2,3)$, \
$-K_{X}^{3}=8$.

\rom{2) ($H_{2}$)} A hypersurface of degree $6$ in~$\Bbb{P}(1^{4},3)$, \
$-K_{X}^{3}=2$.

\rom{3) ($H_{3}$)} A complete intersection of a quadric cone and a quartic
in $\Bbb{P}(1^{5},2)$, \ $-K_{X}^{3}=4$.

\rom{4)} An anticanonical model of a \rom{``}weak Fano $3$-fold\rom{''}~$V$
with canonical Gorenstein singularities {\rm(}that is, $-K_{V}$ is a
numerically effective and big Cartier divisor and
$\varphi_{|-rK_{V}|}(V)=X$ for $r\gg 0${\rm)} such that
$V$ is a double covering of the rational scroll
$\Bbb{F}(d_{1},d_{2},d_{3})=\operatorname{Proj}
\bigl(\,\bopl_{i=1}^{3}\Cal{O}_{\Bbb{P}^{1}}(d_{i})\bigr)$
branched over a divisor rationally equivalent to
$4M+2\bigl(2-\sum_{i=1}^{3}d_{i}\bigr)L$, where $M$ is the class
of the tautological sheaf on~$\Bbb{F}(d_{1},d_{2},d_{3})$
and $L$ is the class of a fibre of the natural projection of
$\Bbb{F}(d_{1},d_{2},d_{3})$ to~$\Bbb{P}^{1}$. Here the
following cases are possible\rm:

($H_{4}$) \ $d_{1}=1$, \ $d_{2}=1$, \ $d_{3}=1$, \ $-K_{X}^{3}=6$;

($H_{5}$) \ $d_{1}=2$, \ $d_{2}=1$, \ $d_{3}=0$, \ $-K_{X}^{3}=6$;

($H_{6}$) \ $d_{1}=2$, \ $d_{2}=1$, \ $d_{3}=1$, \ $-K_{X}^{3}=8$;

($H_{7}$) \ $d_{1}=2$, \ $d_{2}=2$, \ $d_{3}=0$, \ $-K_{X}^{3}=8$;

($H_{8}$) \ $d_{1}=2$, \ $d_{2}=2$, \ $d_{3}=1$, \ $-K_{X}^{3}=10$;

($H_{9}$) \ $d_{1}=2$, \ $d_{2}=2$, \ $d_{3}=2$, \ $-K_{X}^{3}=12$;

($H_{10}$) \ $d_{1}=3$, \ $d_{2}=0$, \ $d_{3}=0$, \ $-K_{X}^{3}=6$;

($H_{11}$) \ $d_{1}=3$, \ $d_{2}=1$, \ $d_{3}=0$, \ $-K_{X}^{3}=8$;

($H_{12}$) \ $d_{1}=3$, \ $d_{2}=1$, \ $d_{3}=1$, \ $-K_{X}^{3}=10$;

($H_{13}$) \ $d_{1}=3$, \ $d_{2}=2$, \ $d_{3}=0$, \ $-K_{X}^{3}=10$;

($H_{14}$) \ $d_{1}=3$, \ $d_{2}=2$, \ $d_{3}=1$, \ $-K_{X}^{3}=12$;

($H_{15}$) \ $d_{1}=3$, \ $d_{2}=3$, \ $d_{3}=0$, \ $-K_{X}^{3}=12$;

($H_{16}$) \ $d_{1}=3$, \ $d_{2}=3$, \ $d_{3}=1$, \ $-K_{X}^{3}=14$;

($H_{17}$) \ $d_{1}=4$, \ $d_{2}=0$, \ $d_{3}=0$, \ $-K_{X}^{3}=8$;

($H_{18}$) \ $d_{1}=4$, \ $d_{2}=1$, \ $d_{3}=0$, \ $-K_{X}^{3}=10$;

($H_{19}$) \ $d_{1}=4$, \ $d_{2}=2$, \ $d_{3}=0$, \ $-K_{X}^{3}=12$;

($H_{20}$) \ $d_{1}=4$, \ $d_{2}=2$, \ $d_{3}=1$, \ $-K_{X}^{3}=14$;

($H_{21}$) \ $d_{1}=4$, \ $d_{2}=3$, \ $d_{3}=0$, \ $-K_{X}^{3}=14$;

($H_{22}$) \ $d_{1}=4$, \ $d_{2}=3$, \ $d_{3}=1$, \ $-K_{X}^{3}=16$;

($H_{23}$) \ $d_{1}=4$, \ $d_{2}=4$, \ $d_{3}=0$, \ $-K_{X}^{3}=16$;

($H_{24}$) \ $d_{1}=5$, \ $d_{2}=1$, \ $d_{3}=0$, \ $-K_{X}^{3}=12$;

($H_{25}$) \ $d_{1}=5$, \ $d_{2}=2$, \ $d_{3}=0$, \ $-K_{X}^{3}=14$;

($H_{26}$) \ $d_{1}=5$, \ $d_{2}=3$, \ $d_{3}=0$, \ $-K_{X}^{3}=16$;

($H_{27}$) \ $d_{1}=5$, \ $d_{2}=3$, \ $d_{3}=1$, \ $-K_{X}^{3}=18$;

($H_{28}$) \ $d_{1}=5$, \ $d_{2}=4$, \ $d_{3}=0$, \ $-K_{X}^{3}=18$;

($H_{29}$) \ $d_{1}=5$, \ $d_{2}=4$, \ $d_{3}=1$, \ $-K_{X}^{3}=20$;

($H_{30}$) \ $d_{1}=6$, \ $d_{2}=2$, \ $d_{3}=0$, \ $-K_{X}^{3}=16$;

($H_{31}$) \ $d_{1}=6$, \ $d_{2}=3$, \ $d_{3}=0$, \ $-K_{X}^{3}=18$;

($H_{32}$) \ $d_{1}=6$, \ $d_{2}=4$, \ $d_{3}=0$, \ $-K_{X}^{3}=20$;

($H_{33}$) \ $d_{1}=6$, \ $d_{2}=4$, \ $d_{3}=1$, \ $-K_{X}^{3}=22$;

($H_{34}$) \ $d_{1}=6$, \ $d_{2}=5$, \ $d_{3}=0$, \ $-K_{X}^{3}=22$;

($H_{35}$) \ $d_{1}=7$, \ $d_{2}=3$, \ $d_{3}=0$, \ $-K_{X}^{3}=20$;

($H_{36}$) \ $d_{1}=7$, \ $d_{2}=4$, \ $d_{3}=0$, \ $-K_{X}^{3}=22$;

($H_{37}$) \ $d_{1}=7$, \ $d_{2}=5$, \ $d_{3}=0$, \ $-K_{X}^{3}=24$;

($H_{38}$) \ $d_{1}=7$, \ $d_{2}=5$, \ $d_{3}=1$, \ $-K_{X}^{3}=26$;

($H_{39}$) \ $d_{1}=8$, \ $d_{2}=4$, \ $d_{3}=0$, \ $-K_{X}^{3}=24$;

($H_{40}$) \ $d_{1}=8$, \ $d_{2}=5$, \ $d_{3}=0$, \ $-K_{X}^{3}=26$;

($H_{41}$) \ $d_{1}=8$, \ $d_{2}=6$, \ $d_{3}=0$, \ $-K_{X}^{3}=28$;

($H_{42}$) \ $d_{1}=9$, \ $d_{2}=5$, \ $d_{3}=0$, \ $-K_{X}^{3}=28$;

($H_{43}$) \ $d_{1}=9$, \ $d_{2}=6$, \ $d_{3}=0$, \ $-K_{X}^{3}=30$;

($H_{44}$) \ $d_{1}=10$, \ $d_{2}=6$, \ $d_{3}=0$, \ $-K_{X}^{3}=32$;

($H_{45}$) \ $d_{1}=10$, \ $d_{2}=7$, \ $d_{3}=0$, \ $-K_{X}^{3}=34$;

($H_{46}$) \ $d_{1}=11$, \ $d_{2}=7$, \ $d_{3}=0$, \ $-K_{X}^{3}=36$;

($H_{47}$) \ $d_{1}=12$, \ $d_{2}=8$, \ $d_{3}=0$, \ $-K_{X}^{3}=40$.
\endproclaim

\goodbreak

\proclaim{Theorem 1.6}
Let $X$ be a Fano $3$-fold with canonical Gorenstein singularities
such that the linear system $|-K_{X}|$ has no base points and the induced
morphism $\varphi_{|-K_{X}|}\:X\to\Bbb{P}^{n}$ is an embedding, where
$n=-{\frac{K_{X}^{3}}{2}}+2$. Suppose that the anticanonical image
$\varphi_{|-K_{X}|}(X)\subset\Bbb{P}^{n}$ is not an intersection of
quadrics. Then $X$ is one of the following $69$ Fano $3$-folds.

\rom{1) ($T_{1}$)} A hypersurface of degree $4$ in~$\Bbb{P}^{4}$, \
$-K_{X}^{3}=4$.

\rom{2) ($T_{2}$)} A complete intersection of a quadric and a cubic
in~$\Bbb{P}^{5}$, \  $-K_{X}^{3}=6$.

\rom{3) ($T_{3}$)} An anticanonical image of a \rom{``}weak Fano
$3$-fold\rom{''}~$Y$ with canonical Gorenstein singularities
{\rm(}that is, $-K_{Y}$ is a numerically effective and big Cartier divisor
and $\varphi_{|-K_{Y}|}(Y)=\nomathbreak X${\rm)},
where $Y$ is a divisor in the rational scroll
$\operatorname{Proj}(\Cal{O}_{\Bbb{P}^{2}}(2)\oplus
\Cal{O}_{\Bbb{P}^{2}}\oplus\Cal{O}_{\Bbb{P}^{2}})$
and $Y$ is rationally equivalent to the divisor
$2T+F$. Here $T$ is the class of the tautological sheaf on
$\operatorname{Proj}(\Cal{O}_{\Bbb{P}^{2}}(2)\oplus
\Cal{O}_{\Bbb{P}^{2}}\oplus\Cal{O}_{\Bbb{P}^{2}})$ and
$F$ is the pull back of~$\Cal{O}_{\Bbb{P}^{2}}(1)$
under the natural projection onto~$\Bbb{P}^{2}$, \ $-K_{X}^{3}=10$.

\rom{4)} An anticanonical image of a \rom{``}weak Fano $3$-fold\rom{''}~$V$
with canonical Gorenstein singularities {\rm(}that is, $-K_{V}$ is a
numerically effective and big Cartier divisor and
$\varphi_{|-K_{V}|}(V)=\nomathbreak X${\rm)} such that $V$ is a divisor in
$\Bbb{F}(d_{1},d_{2},d_{3},d_{4})=\operatorname{Proj}\bigl(\,
\bopl_{i=1}^{4}\Cal{O}_{\Bbb{P}^{1}}(d_{i})\bigr)$
rationally equivalent to the divisor $3M+(2-\sum_{i=1}^{4}d_{i})L$,
where $M$ is the class of the tautological sheaf on
$\Bbb{F}(d_{1},d_{2},d_{3},d_{4})$ and $L$ is the class of a fibre of the
projection of~$\Bbb{F}(d_{1},d_{2},d_{3},d_{4})$ to~$\Bbb{P}^{1}$.
Here the following cases are possible\rm:

\rom{($T_{4}$)} \ $d_{1}=1$, \ $d_{2}=1$, \ $d_{3}=1$, \ $d_{4}=0$, \
$-K_{X}^{3}=8$;

\rom{($T_{5}$)} \ $d_{1}=1$, \ $d_{2}=1$, \ $d_{3}=1$, \ $d_{4}=1$, \
$-K_{X}^{3}=10$;

\rom{($T_{6}$)} \ $d_{1}=2$, \ $d_{2}=1$, \ $d_{3}=0$, \ $d_{4}=0$, \
$-K_{X}^{3}=8$;

\rom{($T_{7}$)} \ $d_{1}=2$, \ $d_{2}=1$, \ $d_{3}=1$, \ $d_{4}=0$, \
$-K_{X}^{3}=10$;

\rom{($T_{8}$)} \ $d_{1}=2$, \ $d_{2}=1$, \ $d_{3}=1$, \ $d_{4}=1$, \
$-K_{X}^{3}=12$;

\rom{($T_{9}$)} \ $d_{1}=2$, \ $d_{2}=2$, \ $d_{3}=0$, \ $d_{4}=0$, \
$-K_{X}^{3}=10$;

\rom{($T_{10}$)} \ $d_{1}=2$, \ $d_{2}=2$, \ $d_{3}=1$, \ $d_{4}=0$, \
$-K_{X}^{3}=12$;

\rom{($T_{11}$)} \ $d_{1}=2$, \ $d_{2}=2$, \ $d_{3}=1$, \ $d_{4}=1$, \
$-K_{X}^{3}=14$;

\rom{($T_{12}$)} \ $d_{1}=2$, \ $d_{2}=2$, \ $d_{3}=2$, \ $d_{4}=0$, \
$-K_{X}^{3}=14$;

\rom{($T_{13}$)} \ $d_{1}=2$, \ $d_{2}=2$, \ $d_{3}=2$, \ $d_{4}=1$, \
$-K_{X}^{3}=16$;

\rom{($T_{14}$)} \ $d_{1}=2$, \ $d_{2}=2$, \ $d_{3}=2$, \ $d_{4}=2$, \
$-K_{X}^{3}=18$;

\rom{($T_{15}$)} \ $d_{1}=3$, \ $d_{2}=1$, \ $d_{3}=0$, \ $d_{4}=0$, \
$-K_{X}^{3}=10$;

\rom{($T_{16}$)} \ $d_{1}=3$, \ $d_{2}=1$, \ $d_{3}=1$, \ $d_{4}=0$, \
$-K_{X}^{3}=12$;

\rom{($T_{17}$)} \ $d_{1}=3$, \ $d_{2}=2$, \ $d_{3}=0$, \ $d_{4}=0$, \
$-K_{X}^{3}=12$;

\rom{($T_{18}$)} \ $d_{1}=3$, \ $d_{2}=2$, \ $d_{3}=1$, \ $d_{4}=0$, \
$-K_{X}^{3}=14$;

\rom{($T_{19}$)} \ $d_{1}=3$, \ $d_{2}=2$, \ $d_{3}=1$, \ $d_{4}=1$, \
$-K_{X}^{3}=16$;

\rom{($T_{20}$)} \ $d_{1}=3$, \ $d_{2}=2$, \ $d_{3}=2$, \ $d_{4}=0$, \
$-K_{X}^{3}=16$;

\rom{($T_{21}$)} \ $d_{1}=3$, \ $d_{2}=2$, \ $d_{3}=2$, \ $d_{4}=1$, \
$-K_{X}^{3}=18$;

\rom{($T_{22}$)} \ $d_{1}=3$, \ $d_{2}=3$, \ $d_{3}=1$, \ $d_{4}=0$, \
$-K_{X}^{3}=16$;

\rom{($T_{23}$)} \ $d_{1}=3$, \ $d_{2}=3$, \ $d_{3}=2$, \ $d_{4}=0$, \
$-K_{X}^{3}=18$;

\rom{($T_{24}$)} \ $d_{1}=3$, \ $d_{2}=3$, \ $d_{3}=2$, \ $d_{4}=1$, \
$-K_{X}^{3}=20$;

\rom{($T_{25}$)} \ $d_{1}=4$, \ $d_{2}=1$, \ $d_{3}=0$, \ $d_{4}=0$, \
$-K_{X}^{3}=12$;

\rom{($T_{26}$)} \ $d_{1}=4$, \ $d_{2}=2$, \ $d_{3}=0$, \ $d_{4}=0$, \
$-K_{X}^{3}=14$;

\rom{($T_{27}$)} \ $d_{1}=4$, \ $d_{2}=2$, \ $d_{3}=1$, \ $d_{4}=0$, \
$-K_{X}^{3}=16$;

\rom{($T_{28}$)} \ $d_{1}=4$, \ $d_{2}=2$, \ $d_{3}=1$, \ $d_{4}=1$, \
$-K_{X}^{3}=18$;

\rom{($T_{29}$)} \ $d_{1}=4$, \ $d_{2}=2$, \ $d_{3}=2$, \ $d_{4}=0$, \
$-K_{X}^{3}=18$;

\rom{($T_{30}$)} \ $d_{1}=4$, \ $d_{2}=3$, \ $d_{3}=1$, \ $d_{4}=0$, \
$-K_{X}^{3}=18$;

\rom{($T_{31}$)} \ $d_{1}=4$, \ $d_{2}=3$, \ $d_{3}=2$, \ $d_{4}=0$, \
$-K_{X}^{3}=20$;

\rom{($T_{32}$)} \ $d_{1}=4$, \ $d_{2}=3$, \ $d_{3}=2$, \ $d_{4}=1$, \
$-K_{X}^{3}=22$;

\rom{($T_{33}$)} \ $d_{1}=4$, \ $d_{2}=3$, \ $d_{3}=3$, \ $d_{4}=0$, \
$-K_{X}^{3}=22$;

\rom{($T_{34}$)} \ $d_{1}=4$, \ $d_{2}=3$, \ $d_{3}=3$, \ $d_{4}=1$, \
$-K_{X}^{3}=24$;

\rom{($T_{35}$)} \ $d_{1}=4$, \ $d_{2}=4$, \ $d_{3}=2$, \ $d_{4}=0$, \
$-K_{X}^{3}=22$;

\rom{($T_{36}$)} \ $d_{1}=5$, \ $d_{2}=2$, \ $d_{3}=0$, \ $d_{4}=0$, \
$-K_{X}^{3}=16$;

\rom{($T_{37}$)} \ $d_{1}=5$, \ $d_{2}=2$, \ $d_{3}=1$, \ $d_{4}=0$, \
$-K_{X}^{3}=18$;

\rom{($T_{38}$)} \ $d_{1}=5$, \ $d_{2}=3$, \ $d_{3}=1$, \ $d_{4}=0$, \
$-K_{X}^{3}=20$;

\rom{($T_{39}$)} \ $d_{1}=5$, \ $d_{2}=3$, \ $d_{3}=2$, \ $d_{4}=0$, \
$-K_{X}^{3}=22$;

\rom{($T_{40}$)} \ $d_{1}=5$, \ $d_{2}=3$, \ $d_{3}=2$, \ $d_{4}=1$, \
$-K_{X}^{3}=24$;

\rom{($T_{41}$)} \ $d_{1}=5$, \ $d_{2}=3$, \ $d_{3}=3$, \ $d_{4}=0$, \
$-K_{X}^{3}=24$;

\rom{($T_{42}$)} \ $d_{1}=5$, \ $d_{2}=4$, \ $d_{3}=2$, \ $d_{4}=0$, \
$-K_{X}^{3}=24$;

\rom{($T_{43}$)} \ $d_{1}=5$, \ $d_{2}=4$, \ $d_{3}=3$, \ $d_{4}=0$, \
$-K_{X}^{3}=26$;

\rom{($T_{44}$)} \ $d_{1}=5$, \ $d_{2}=4$, \ $d_{3}=3$, \ $d_{4}=1$, \
$-K_{X}^{3}=28$;

\rom{($T_{45}$)} \ $d_{1}=6$, \ $d_{2}=2$, \ $d_{3}=0$, \ $d_{4}=0$, \
$-K_{X}^{3}=18$;

\rom{($T_{46}$)} \ $d_{1}=6$, \ $d_{2}=3$, \ $d_{3}=1$, \ $d_{4}=0$, \
$-K_{X}^{3}=22$;

\rom{($T_{47}$)} \ $d_{1}=6$, \ $d_{2}=3$, \ $d_{3}=2$, \ $d_{4}=0$, \
$-K_{X}^{3}=24$;

\rom{($T_{48}$)} \ $d_{1}=6$, \ $d_{2}=4$, \ $d_{3}=2$, \ $d_{4}=0$, \
$-K_{X}^{3}=26$;

\rom{($T_{49}$)} \ $d_{1}=6$, \ $d_{2}=4$, \ $d_{3}=3$, \ $d_{4}=0$, \
$-K_{X}^{3}=28$;

\rom{($T_{50}$)} \ $d_{1}=6$, \ $d_{2}=4$, \ $d_{3}=3$, \ $d_{4}=1$, \
$-K_{X}^{3}=30$;

\rom{($T_{51}$)} \ $d_{1}=6$, \ $d_{2}=4$, \ $d_{3}=4$, \ $d_{4}=0$, \
$-K_{X}^{3}=30$;

\rom{($T_{52}$)} \ $d_{1}=6$, \ $d_{2}=5$, \ $d_{3}=3$, \ $d_{4}=0$, \
$-K_{X}^{3}=30$;

\rom{($T_{53}$)} \ $d_{1}=7$, \ $d_{2}=3$, \ $d_{3}=1$, \ $d_{4}=0$, \
$-K_{X}^{3}=24$;

\rom{($T_{54}$)} \ $d_{1}=7$, \ $d_{2}=4$, \ $d_{3}=2$, \ $d_{4}=0$, \
$-K_{X}^{3}=28$;

\rom{($T_{55}$)} \ $d_{1}=7$, \ $d_{2}=4$, \ $d_{3}=3$, \ $d_{4}=0$, \
$-K_{X}^{3}=30$;

\rom{($T_{56}$)} \ $d_{1}=7$, \ $d_{2}=5$, \ $d_{3}=3$, \ $d_{4}=0$, \
$-K_{X}^{3}=32$;

\rom{($T_{57}$)} \ $d_{1}=7$, \ $d_{2}=5$, \ $d_{3}=4$, \ $d_{4}=0$, \
$-K_{X}^{3}=34$;

\rom{($T_{58}$)} \ $d_{1}=7$, \ $d_{2}=5$, \ $d_{3}=4$, \ $d_{4}=1$, \
$-K_{X}^{3}=36$;

\rom{($T_{59}$)} \ $d_{1}=8$, \ $d_{2}=4$, \ $d_{3}=2$, \ $d_{4}=0$, \
$-K_{X}^{3}=30$;

\rom{($T_{60}$)} \ $d_{1}=8$, \ $d_{2}=5$, \ $d_{3}=3$, \ $d_{4}=0$, \
$-K_{X}^{3}=34$;

\rom{($T_{61}$)} \ $d_{1}=8$, \ $d_{2}=5$, \ $d_{3}=4$, \ $d_{4}=0$, \
$-K_{X}^{3}=36$;

\rom{($T_{62}$)} \ $d_{1}=8$, \ $d_{2}=6$, \ $d_{3}=4$, \ $d_{4}=0$, \
$-K_{X}^{3}=38$;

\rom{($T_{63}$)} \ $d_{1}=9$, \ $d_{2}=5$, \ $d_{3}=3$, \ $d_{4}=0$, \
$-K_{X}^{3}=36$;

\rom{($T_{64}$)} \ $d_{1}=9$, \ $d_{2}=6$, \ $d_{3}=4$, \ $d_{4}=0$, \
$-K_{X}^{3}=40$;

\rom{($T_{65}$)} \ $d_{1}=9$, \ $d_{2}=6$, \ $d_{3}=5$, \ $d_{4}=0$, \
$-K_{X}^{3}=42$;

\rom{($T_{66}$)} \ $d_{1}=10$, \ $d_{2}=6$, \ $d_{3}=4$, \ $d_{4}=0$, \
$-K_{X}^{3}=42$;

\rom{($T_{67}$)} \ $d_{1}=10$, \ $d_{2}=7$, \ $d_{3}=5$, \ $d_{4}=0$, \
$-K_{X}^{3}=46$;

\rom{($T_{68}$)} \ $d_{1}=11$, \ $d_{2}=7$, \ $d_{3}=5$, \ $d_{4}=0$, \
$-K_{X}^{3}=48$;

\rom{($T_{69}$)} \ $d_{1}=12$, \ $d_{2}=8$, \ $d_{3}=6$, \ $d_{4}=0$, \
$-K_{X}^{3}=54$.
\endproclaim

In the smooth case Theorems~1.4--1.6 were proved by Iskovskikh
(see~\cite{15}, \cite{17}).

\goodbreak

\remark{Remark 1.7}
For any Fano 3-fold $X$ with canonical Gorenstein singularities, there is
a birational morphism $f\:V\to X$ (called the {\it terminal modification\/}
of~$X$) such that $K_{V}\sim f^{*}(K_{V})$ and $V$ has terminal Gorenstein
singularities. The existence of~$f$ follows from the Minimal Model Program
and the contraction theorem (see~\cite{93}). On the other hand, if~$V$ is
any ``weak Fano 3-fold'' (that is, a variety whose anticanonical class
$-K_{V}$ is numerically effective and big) with canonical Gorenstein
singularities, then the contraction theorem implies that there is a
birational morphism $f\:V\to X$ such that $X$ is a Fano 3-fold with
canonical Gorenstein singularities.
\endremark

\medskip

In what follows we use the symbols $B_{k}$, $B_{4}^m$, $H_{i}$, and $T_{j}$
to denote the corresponding Fano 3-folds described in Theorems~1.4, 1.5,
and~1.6.

\remark{Remark 1.8}
The 3-fold $H_{1}$ is a double covering of a cone over the Veronese surface,
and $H_{2}$ is a double covering of~$\Bbb{P}^{3}$ ramified in a sextic
surface (which may be singular). The 3-fold $H_{3}$ is a double covering
of a quadric 3-fold (possibly singular), and $H_{4}$ is a double covering
of $\Bbb{P}^{1}\times\Bbb{P}^{2}$ branched over a divisor of bidegree
$(2,4)$. The 3-fold $H_{6}$ is the blow up of a hypersurface of degree~4
in $\Bbb{P}(1^{4},2)$ (or, equivalently, the blow up of a double covering
of $\Bbb{P}^{3}$ branched over a quartic surface) along the intersection
of two different divisors in the half-anticanonical linear system.
The 3-fold $H_{9}$ is isomorphic to the product $\Bbb{P}^{1}\times S_{2}$,
where $S_{2}$ is a del~Pezzo surface of degree~$2$ with at worst Du~Val
singularities. The 3-fold $H_{10}$ is a hypersurface of degree $10$ in
$\Bbb{P}(1^{2},3^{2},5)$, \ $H_{17}$ is a hypersurface of degree $12$ in
$\Bbb{P}(1^{2},4^{2},6)$, and $T_{3}$ is a hypersurface of degree $5$ in
$\Bbb{P}(1^{3},2^{2})$. The 3-fold $T_{5}$ is a divisor of bidegree
$(1,3)$ in~$\Bbb{P}^{1}\times\Bbb{P}^{3}$. The 3-fold $T_{8}$ is obtained
by blowing up a plane cubic curve on a cubic 3-fold in~$\Bbb{P}^{4}$.
The 3-fold $T_{14}$ is the product $\Bbb{P}^{1}\times S_{3}$, where
$S_{3}$ is a cubic surface in~$\Bbb{P}^{3}$ with at worst Du~Val
singularities.
\endremark

\remark{Remark 1.9}
The 3-folds $H_{1}$, $H_{2}$, $H_{3}$, $H_{4}$, $H_{6}$,
$H_{9}$, $T_{1}$, $T_{2}$, $T_{5}$, $T_{8}$ and~$T_{14}$ are the only
3-folds among $H_{i}$ and $T_{j}$ that can be chosen smooth
\rom(see~\cite{15}, \cite{88}\rom).
\endremark

\medskip

The 3-folds $H_{i}$ and $T_{j}$ are rationally connected (see~\cite{95}).
Moreover, the majority of~$H_{i}$ and $T_{j}$ must be rational, although
some of them are definitely not. For example, it is well known that
sufficiently general 3-folds $H_{1}$, $H_{2}$, $H_{3}$, $H_{4}$, $H_{6}$,
$T_{1}$, $T_{2}$, $T_{8}$ are non-rational (see \cite{22}, \cite{62},
\cite{54}, \cite{37}, \cite{16}, \cite{52}, \cite{98}).
Their non-rationality can also be proved in the smooth and some singular
cases (see~\cite{31}--\cite{34},
\cite{28}, \cite{29}, \cite{111}, \cite{5},
\cite{82}, \cite{63}, \cite{64}, \cite{6},
\cite{102}, \cite{7}, \cite{59}). However, all of these cases include
examples of rational singular 3-folds $H_{i}$ and $T_{j}$ even when the
singularities are isolated ordinary double points (see~\cite{90}). The 3-folds
$H_{9}$, $T_{5}$ and~$T_{14}$ are always rational by Remark~1.8. In this
paper we prove the following result.

\proclaim{Proposition 1.10}
The $3$-folds $H_{i}$ and $T_{j}$ are rational for
$i\in\{8, 9, 22,\allowmathbreak
26, \allowmathbreak
27, 28, 29, 31, \dots, 47\}$ and
$j\in\{5, 10, 11, 12, 13, 14, 17, \dots, 69\}$. On the other hand,
sufficiently general $3$-folds $H_{i}$ and~$T_{j}$ are non-rational for
$i\leqslant 7$ and $j\in\{1, 2, 3, 4, 6, 7, 8, 9\}$.
\endproclaim

There are birational relations between some of the 3-folds
 $H_{i}$, $T_{j}$, $B_{k}$ and~$B_{4}^m$. The simplest example is a
projection from a $\operatorname{cDV}$-\allowbreak point: the anticanonical
model of the blow up of a $\operatorname{cDV}$-point on any of the 3-folds
$H_{i}$ or $T_{j}$ of anticanonical degree $d\geqslant 4$ must be one of
the 3-folds $H_{i}$, $T_{j}$, $B_{k}$ or $B_{4}^m$ of anticanonical degree
$d-2$. For example, $B_{4}^{4}$ is birationally isomorphic to~$H_{17}$, \
$H_{5}$ is birationally isomorphic to
$H_{1}$ with a~$\operatorname{cDV}$-point
(see~\cite{5}, Lemma~3.4), and the 3-folds $T_{1}$ and $T_{2}$ having a
$\operatorname{cDV}$-point are birationally isomorphic to the singular
3-folds $H_{2}$ and~$T_{1}$ respectively. Moreover, there are many
non-obvious birational transformations of the 3-folds $H_{i}$ and~$T_{j}$.

\example{Example 1.11}
In the notation of Theorem~1.6, let $X$ be a sufficiently general 3-fold
$T_{7}$, and let $V$ be the corresponding weak Fano 3-fold
$V\subset\Bbb{F}(2,1,1,0)$. Then $V$ is smooth (see  the proof of
Theorem~1.6) and $-K_{V}$ has trivial intersection with only one 
rational curve $Y_{4}\subset\Bbb{F}(2,1,1,0)$ (see Corollary~2.20).
It follows that the birational morphism $\varphi_{|-K_{V}|}\:V\to X$
contracts the curve $Y_{4}$ to an ordinary double point of~$X$.
Let $f\:V\dasharrow \widetilde{V}$ be a flop in the curve~$Y_{4}$.
Then one can find a birational morphism $g\:\widetilde{V}\to Y$
such that there is a double covering $\pi\:Y\to\Bbb{P}^{3}$
branched over a smooth hypersurface of degree~$4$, that is, $Y$ is a
double space of index~2 (see~\cite{60}, \cite{61}).
Moreover, the birational morphism~$g$ is a blow up of a smooth rational
curve $C\subset Y$ with $-K_{Y}\cdot C=2$. All these constructions of
birational maps are easily seen to be reversible. More precisely, let
$\pi\:Y\to\Bbb{P}^{3}$ be any double covering branched over a smooth quartic
surface, and let $C\subset Y$ be any non-singular rational curve with
$-K_{Y}\cdot C=2$. Then one always can construct the corresponding Fano
3-fold $T_{7}$ (see~\cite{56}, \S\,4.4.1).
\endexample

There are only two Fano 3-folds with canonical Gorenstein singularities
whose anticanonical divisor is divisible (in the Picard group) by an
integer greater than~2. These are $\Bbb{P}^{3}$ and a quartic 3-fold
$Q\subset\Bbb{P}^{4}$. Fano 3-folds with canonical Gorenstein singularities
whose anticanonical divisor is divisible by~2 are called
{\it del~Pezzo $3$-folds\/} (see~\cite{88}, Theorem~3.3.1).
It is easy to prove explicitly that $H_1$ is the only del~Pezzo 3-fold
among the 3-folds $H_{i}$ and~$T_{j}$. This is also confirmed by the
classification of del~Pezzo 3-folds
(see~\cite{78}, \cite{79}, \cite{57}, \cite{117}).

\remark{Remark 1.12}
The 3-folds $H_{i}$ and $T_{j}$ are naturally birationally isomorphic to
del~Pezzo fibrations of degree $2$ and~$3$ respectively, except for the
following cases: $H_{1}$, $H_{2}$, $H_{3}$, $T_{1}$, $T_{2}$ and~$T_{3}$.
On the other hand, sufficiently general 3-folds $H_{1}$,
$H_{2}$, $H_{3}$, $T_{1}$ and~$T_{2}$ are not birationally isomorphic to
any del~Pezzo fibration of degree $2$ or~$3$
(see~\cite{22}, \cite{16}, \cite{29}, \cite{63}, \cite{6}, \cite{7}).
\endremark

\remark{Remark 1.13}
It is well known that 3-folds with a pencil of del~Pezzo surfaces of degree
$2$ or~$3$ are unirational (see~\cite{23}--\cite{25}).
Therefore the 3-folds $H_{i}$ and~$T_{j}$ are unirational for
$i\geqslant 4$ and~$j\geqslant 4$. The 3-fold $H_{3}$ is also known to be
unirational (see~\cite{115}, \cite{16}). The proof of Proposition~5.5 below
implies that the 3-fold $T_{3}$ is unirational since it is birationally
equivalent to a conic bundle with a rational multisection. The 3-fold $T_{2}$
is also unirational (see~\cite{71},~\cite{115},
\cite{36},~\cite{16}). However, it is still unknown whether a general
quartic 3-fold $T_{1}$ is unirational or not, despite several examples
of smooth unirational quartic 3-folds
(see~\cite{118}, \cite{22}, \cite{16}, \cite{100}).
Unfortunately, nothing is known about the unirationality of general 3-folds
$H_{1}$ and~$H_{2}$. It is expected that general 3-folds $H_{2}$ are not
unirational (see~\cite{97}, Conjecture~4.1.6).
\endremark

\medskip

The authors are very grateful to M.~M.~Grinenko, V.~A.~Iskovskikh,
N.~F.~Zak, A.~Corti, S.~A.~Kudryavtsev, V.~S.~Kulikov, J.~Park,
Yu.~G.~Prokhorov, A.~V.~Pukh\-li\-kov, D.~A.~Stepanov and
V.~V.~Shokurov for fruitful conversations.

\head
\S\,2. Preliminaries
\endhead

In what follows all varieties are assumed to be projective, normal, and
defined over~$\Bbb{C}$.

\proclaim{Proposition 2.1 \rm(\cite{96}, Proposition~3.1.6)}
Suppose that $\rho\:V\to X$ is a finite morphism,
$D_{X}$ is an effective $\Bbb{Q}$\kkk-divisor on~$X$, and
$D_{V}=\rho^{*}(D_{X})-K_{V/X}$, that is,
$K_{V}+D_{V}=\rho^{*}(K_{X}+D_{X})$. The singularities of the log
pair $(V, D_{V})$ are Kawamata log terminal \rom{(}see~\cite{93},
\cite{96}\rom{)} if and only if the singularities  of the log pair
$(X,D_{X})$ are Kawamata log terminal.
\endproclaim

We note that Kawamata log terminal singularities are canonical if the
canonical divisor is a Cartier divisor. Hence Proposition~2.1 yields
the following result.

\proclaim{Corollary 2.2}
Let $X$ be a smooth variety, and
$\rho\:V\to X$ a double covering branched over a reduced effective divisor
$D\subset X$. The singularities of~$V$ are canonical if and only if the
singularities of the log pair $\bigl(X,{\frac{1}{2}}D\bigr)$ are Kawamata
log terminal.
\endproclaim

\proclaim{Theorem 2.3 \rm(\cite{96}, Theorem~4.5.1)}
Let $X$ be a smooth variety, $\Cal{H}$ a linear system on~$X$ whose base
locus has codimension at least~$2$, and $D$ a sufficiently general divisor
in~$\Cal{H}$. Suppose that for every point $x\in X$ there is a divisor
$H\in\Cal{H}$ such that the singularities of~$(X, H)$ are canonical in the
neighbourhood of~$x$. Then the singularities of the log pair $(X, D)$
are canonical.
\endproclaim

The next result is Theorem~7.9 of~\cite{96}, which was proved in~\cite{119}.

\proclaim{Theorem 2.4}
Let $X$ be a normal variety such that $\omega_{X}$ is locally free, and let
$S\subset X$ be an effective Cartier divisor on~$X$. Then~$S$ has canonical
singularities if and only if the singularities of the log pair $(X, S)$
are canonical.
\endproclaim

Theorems~2.3 and~2.4 yield the following result.

\proclaim{Corollary 2.5}
Let $X$ be a smooth variety, $\Cal{H}$ a linear system on~$X$ whose base
locus has codimension at least~$2$, and $D$ a sufficiently general divisor
in~$\Cal{H}$. Suppose that for every point $x\in X$ there is a divisor
$H\in\Cal{H}$ whose singularities are canonical in a neighbourhood of~$x$.
Then $D$ has canonical singularities.
\endproclaim

The following result is implied by Theorem~4.8 of~\cite{96}.

\proclaim{Theorem 2.6}
Let $X$ be a smooth variety, $\Cal{H}$ a linear system on~$X$, \ $D$
a sufficiently general divisor in~$\Cal{H}$, and
$\lambda\in\Bbb{Q}\cap [0,1)$. Suppose that for every point $x\in X$
there is $H\in\Cal{H}$ such that the log pair $(X, \lambda H)$ has
Kawamata log terminal singularities in the neighbourhood of~$x$.
Then the singularities of~$(X,\lambda D)$ are Kawamata log terminal.
\endproclaim

Corollary~2.2 and Theorem~2.6 imply the following result.

\proclaim{Corollary 2.7}
Let $X$ be a smooth variety, $\Cal{H}$ a linear system on~$X$, and
$D$ a sufficiently general divisor in~$\Cal{H}$. Suppose that for every
point $O\in X$ there is an effective reduced divisor
$H\in\Cal{H}$ such that there is a double covering
$\beta\:Y\to X$ branched over $H\subset X$ with $Y$ having canonical
singularities in the neighbourhood of~$\beta^{-1}(O)$. Let $\rho\:V\to X$
be the double covering branched over $D\subset X$. Then $V$ has canonical
singularities.
\endproclaim

\remark{Remark 2.8}
It is also easy to deduce Corollary~2.7 from Corollary~2.5.
Indeed, in the notation of Corollary~2.7, let $B\subset X$ be a
divisor with $D\sim 2B$. We put
$U=\operatorname{Proj}(\Cal{O}_{X}\oplus\Cal{O}_{X}(B))$. Let
$M$ be the tautological line bundle on~$U$, and let $f\:U\to X$ be the
natural projection. Then $Y$ may be regarded as a divisor on~$U$
in the linear system $|2M|$ such that $\beta=f|_{Y}$. We may assume
that $\Cal{H}=|H|$ without loss of generality. Hence we can identify
$V$ with a sufficiently general divisor in the linear system~$|2M|$.
The base locus of~$|2M|$ is contained in~$Y\cap f^{-1}(H)$ because
$2S+f^{-1}(H)\sim 2M$ and $S\cap Y=\varnothing$, where $S\sim M-f^{*}(B)$
is a \emph{negative section} of~$f\:U\to X$. Therefore the base locus of
$|2M|$ has codimension at least~$2$, and the singularities of~$V\in |2M|$
are canonical by Corollary~2.5.
\endremark

We recall the following classical result (see~\cite{53},
\cite{1},~\cite{114}).

\proclaim{Theorem 2.9}
Suppose that $X$ is a normal algebraic surface and $O\in\nomathbreak X$
is an isolated singular point such that the singularities of~$X$ are
canonical in the neighbourhood of~$O$, that is, $O$ is a Du~Val singular
point on~$X$. Then $O\in\nomathbreak  X$ is a hypersurface quasi-homogeneous
singularity and is locally isomorphic to the singularity
$(0,0,0)\in\Bbb{C}^{3}\cong\operatorname{Spec}(\Bbb{C}[x,y,z])$
of one of the following types\rm:

\rom{($\Bbb{A}_{n}$)} \ $x^{2}+y^{2}+z^{n+1}=0$, \
$\operatorname{wt}(x)=n+1$, \
$\operatorname{wt}(y)=n+1$, \
$\operatorname{wt}(z)=2$, \
$n\geqslant 1$;

\rom{($\Bbb{D}_{n}$)} \ $x^{2}+y^{2}z+z^{n-1}=0$, \
$\operatorname{wt}(x)=n-1$, \ $\operatorname{wt}(y)=n-2$, \
$\operatorname{wt}(z)=2$, \ $n\geqslant 4$;

\rom{($\Bbb{E}_{6}$)} \ $x^{2}+y^{3}+z^{4}=0$, \
$\operatorname{wt}(x)=6$, \
$\operatorname{wt}(y)=4$, \ $\operatorname{wt}(z)=3$;

\rom{($\Bbb{E}_{7}$)} \ $x^{2}+y^{3}+yz^{3}=0$, \
$\operatorname{wt}(x)=9$, \ $\operatorname{wt}(y)=6$, \
$\operatorname{wt}(z)=4$;

\rom{($\Bbb{E}_{8}$)} \ $x^{2}+y^{3}+z^{5}=0$, \
$\operatorname{wt}(x)=15$, \ $\operatorname{wt}(y)=10$, \
$\operatorname{wt}(z)=6$.
\endproclaim

The following result is proved in \S\,12.3, \S\,12.6
and~\S\,13.1 of~\cite{1}.

\proclaim{Theorem 2.10}
Let $X\subset\Bbb{C}^{3}\cong\operatorname{Spec}(\Bbb{C}[x,y,z])$ be a
hypersurface $f(x,y,z)=0$ such that the origin $O\in\Bbb{C}^{3}$ is an
isolated singular point of~$X$. Write
$$
f(x,y,z)=f_{d}(x,y,z)+f_{d+1}(x,y,z)+\ldots,
$$
where $f_{i}(x,y,z)$ is a quasi-homogeneous polynomial of quasi-homogeneous
degree $i\geqslant 2$ with respect to positive integer weights
$\operatorname{wt}(x)=a$, \
$\operatorname{wt}(y)=b$, \
$\operatorname{wt}(z)=c$.
Suppose that the origin $O\in\Bbb{C}^{3}$ is an isolated singular point
of the hypersurface $f_{d}(x,y,z)=0$, where
$2a\leqslant d$, \ $2b\leqslant d$, \
$2c\leqslant d$ and~$a+b+c>d$.

\rom{1)} If $(a,b,c)=(n+1,n+1,2)$, then $O\in X$ is a singularity of
type $\Bbb{A}_{n}$.

\rom{2)} If $(a,b,c)=(n-1,n-2,2)$, then $O\in X$ is a singularity of type
$\Bbb{D}_{n}$.

\rom{3)} If $(a,b,c)=(6,4,3)$, then $O\in X$ is a singularity of type
$\Bbb{E}_{6}$.

\rom{4)} If $(a,b,c)=(9,6,4)$, then $O\in X$ is a singularity of type
$\Bbb{E}_{7}$.

\rom{5)} If $(a,b,c)=(15,10,6)$, then $O\in X$ is a singularity of type
$\Bbb{E}_{8}$.
\endproclaim

The following result is due to Enriques (see~\cite{77},
\cite{15}, \cite{69}, \cite{17}).

\proclaim{Theorem 2.11}
Let $X\subset\Bbb{P}^{n}$ be a variety of degree
$n-\operatorname{dim}(X)+1$ such that $X$ is not contained in any
hyperplane. Then $X$ is one of the following varieties\rom:

\rom{1)} a projective space $\Bbb{P}^{n}$\rom;

\rom{2)} a quadric hypersurface in~$\Bbb{P}^{n}$\rom;

\rom{3)} the image of a rational scroll
$\Bbb{F}(d_{1}, \dots, d_{k})=
\allowmathbreak
\operatorname{Proj}\bigl(\,
\bopl_{i=1}^{k}\Cal{O}_{\Bbb{P}^{1}}(d_{i})\bigr)$
under the map given by the tautological line bundle, where
$0\ne d_{1}\geqslant\dots\geqslant d_{k}\geqslant 0$ and
$n+1=
\allowmathbreak
\sum_{i=1}^{k}(d_{i}+1)$\rom;

\rom{4)} a Veronese surface in~$\Bbb{P}^{5}$ when $n=5$\rom;

\rom{5)} a cone in~$\Bbb{P}^{n}$ over the Veronese surface
in~$\Bbb{P}^{5}$.
\endproclaim

\goodbreak
It is easy to see that the varieties in Theorem~2.11 have the smallest
possible degree among all varieties of the same dimension in~$\Bbb{P}^{n}$.

Using the Kawamata--Viehweg vanishing theorem (see~\cite{91},
\cite{120}) and elementary properties of linear systems on K3 surfaces
(see~\cite{116}), we get the following well-known result
(see~\cite{15}, \cite{17},~\cite{88}).

\proclaim{Theorem 2.12}
Let $X$ be a Fano $3$-fold with canonical Gorenstein singularities such that
the linear system $|-K_{X}|$ has no base points but the anticanonical divisor
$-K_{X}$ is not very ample. Then
$\varphi_{|-K_{X}|}\:X\to V\subset\Bbb{P}^{n}$ is a double covering
and $V\subset\Bbb{P}^{n}$ is a subvariety of minimal degree, that is,
$\operatorname{deg}(V)=n-2$, where $n=-{\frac{1}{2}}K_{X}^{3}+2$.
\endproclaim

The following result is a theorem of Noether--Enriques--Petri
(see \cite{44}, \cite{81}).

\proclaim{Theorem 2.13}
Let $C\subset\Bbb{P}^{g-1}$ be a canonically embedded smooth irreducible
curve whose genus $g(C)$ is at least~$3$. Then the following assertions
hold.

\rom{1)} The curve $C\subset\Bbb{P}^{g-1}$ is projectively normal.

\rom{2)} If $g(C)=3$, then $C$ is a plane quartic curve.

\rom{3)} If $g(C)\geqslant 4$, then the graded ideal $I_{C}$
of the curve $C\subset\Bbb{P}^{g-1}$ is generated by the components
of degree $2$ and~$3$, that is, the curve $C\subset\Bbb{P}^{g-1}$
is cut out by quadrics and cubics in~$\Bbb{P}^{g-1}$ in the
scheme-theoretic sense.

\rom{4)} If $g(C)\geqslant 4$, then the graded ideal $I_{C}$ of the
curve $C\subset\Bbb{P}^{g-1}$ is generated by the component of degree~$2$
except for the following two cases\rom:

-- the curve $C$ is trigonal, that is, there is a map
$\psi\:C\to\Bbb{P}^{1}$ of degree~$3$\rom;

-- the curve $C$ is isomorphic to a smooth plane quintic
\rom{(}in particular, $g(C)=\nomathbreak 6$\rom{)}.

\rom{5)} In the trigonal case, quadrics through $C$ in~$\Bbb{P}^{g-1}$
cut out either an irreducible \rom{(}possibly singular\rom{)} quadric
surface when $g(C)=4$,  or a smooth irreducible surface of degree
$g-2$ which is the image of
$\operatorname{Proj}(\Cal{O}_{\Bbb{P}^{1}}(d_{1})
\oplus\Cal{O}_{\Bbb{P}^{1}}(d_{2}))$
under the map given by the tautological line bundle, where
$d_{1}\geqslant d_{2}>0$ and~$g=d_{1}+d_{2}+2$.

\rom{6)} If $C$ is isomorphic to a smooth plane quintic, then quadrics
through $C$ in~$\Bbb{P}^{5}$ cut out a Veronese surface.
\endproclaim

The following result is a corollary of Theorem~2.13
(see~\cite{15}, \cite{17},~\cite{88}).

\proclaim{Theorem 2.14}
Let $X\subset\Bbb{P}^{n}$ be an anticanonically embedded Fano $3$-fold
with canonical singularities, that is,
$-K_{X}\sim\Cal{O}_{\Bbb{P}^{n}}(1)|_{X}$ and
$n=-{\frac{1}{2}}K_{X}^{3}+2$. Then the following assertions hold.

\rom{1)} The $3$-fold $X$ is projectively normal in~$\Bbb{P}^{n}$.

\rom{2)} If $-K_{X}^{3}=4$, then $X$ is a quartic $3$-fold in~$\Bbb{P}^{4}$.

\rom{3)} If $-K_{X}^{3}\geqslant 6$, then the graded ideal $I_{X}$ of the
$3$-fold $X\subset\Bbb{P}^{n}$ is generated by the components of degree
$2$ and~$3$.

\rom{4)} If $-K_{X}^{3}\geqslant 6$, then the graded ideal $I_{X}$ of the
$3$-fold $X\subset\Bbb{P}^{n}$ is generated by the component of degree~$2$
except for the case when, for a general linear subspace
$\Pi\subset\Bbb{P}^{n}$ of codimension~$2$, the curve $X\cap\Pi$
is either a canonically embedded smooth trigonal curve or a canonically
embedded smooth plane quintic curve and
$\operatorname{deg}(X\subset\Bbb{P}^{n})=10$.

\rom{5)} In the trigonal case, quadrics through $X$ in~$\Bbb{P}^{n}$
cut out either an irreducible \rom{(}possibly singular\rom{)} quadric
$4$-fold when $-K_{X}^{3}=6$, or a $4$-fold of degree $n-3$ which is
the image of a rational scroll
$\operatorname{Proj}\bigl(\,\bopl_{i=1}^{4}
\Cal{O}_{\Bbb{P}^{1}}(d_{i})\bigr)$
under the map given by the tautological line bundle, where
$0\ne d_{1}\geqslant\dots\geqslant d_{4}\geqslant 0$
and $n+1=\sum_{i=1}^{4}(d_{i}+1)$.

\rom{6)} If $X\cap\Pi$ is a canonically embedded plane quintic, then quadrics
through $X$ in~$\Bbb{P}^{7}$ cut out a $4$-dimensional cone over a Veronese
surface.
\endproclaim

\proclaim{Proposition 2.15 \rm(\cite{110}, Claim~6.9)}
Let $X$ be a $3$-fold with composite Du~Val {\rm(}$\operatorname{cDV}${\rm)}
singularities, $\Gamma\subset\operatorname{Sing}(X)$ a smooth curve regarded
as a reduced subscheme of~$X$, and $f\:V\to X$ the blow up of~$\Gamma$.
Then $V$ has at most $\operatorname{cDV}$-singularities and
$K_{V}\sim f^{*}(K_{X})$, that is, the map $f$ is crepant.
\endproclaim

The following result is proved in~\cite{47} and is a special case of a
conjectural rationality criterion for standard 3-dimensional conic bundles
(see~\cite{86}, \cite{87}, \cite{18},~\cite{19}).

\proclaim{Theorem 2.16}
Suppose that $Y$ is a smooth $3$-fold, $Z$ is either $\Bbb{P}^{2}$
or a minimal rational ruled surface $\Bbb{F}_{r}$, and
$\xi\:Y\to Z$ is a conic bundle with $\operatorname{Pic}(Y/Z)=\Bbb{Z}$ and
$|2K_{Z}+\Delta|\ne\varnothing$, where $\Delta\subset Z$ is the degeneration
divisor of~$\xi\:Y\to Z$. Then $Y$ is non-rational.
\endproclaim

\remark{Remark 2.17}
In the notation of Theorem~2.16, the hypothesis
$|2K_{Z}+\Delta|=\varnothing$ implies that the $3$-fold $Y$ is rational
except for the case when there is a commutative diagram
$$
\xymatrix{
X\ar@{->}[d]_{\chi}\ar@{-->}[rr]^{\alpha}&& Y\ar[d]^{\xi}\\
\Bbb{P}^2\ar@{-->}[rr]^{\beta}&&Z}
$$
where $\alpha$ and $\beta$ are birational maps, $X$ is a smooth 3-fold,
and $\chi\:X\to\Bbb{P}^2$ is a conic bundle with
$\operatorname{Pic}(X/\Bbb{P}^2)=\Bbb{Z}$ whose degeneration divisor
$D\subset\Bbb{P}^2$ is a quintic curve and the double covering
$\psi\:\widetilde{D}\to D$ induced by~$\chi$ corresponds to an even
$\theta$-characteristic (see~\cite{18}).
\endremark

The following result was proved in~\cite{101}. It is a particular case of
a more general result in~\cite{95} (see also \cite{94}, \cite{98}), which
generalizes the standard degeneration technique (see~\cite{21}).

\proclaim{Theorem 2.18}
Let $\xi\:Y\to Z$ be a flat proper morphism with irreducible and reduced
geometric fibres. Then there are countably many closed subsets
$Z_{i}\subset Z$ such that the fibre $\xi^{-1}(s)$ over a closed point
$s\in Z$ is ruled if and only if~$s\in\bigcup Z_{i}$.
\endproclaim

\proclaim{Proposition 2.19~\rm\cite{114}}
Let $V$ be a rational scroll
$\operatorname{Proj}\bigl(\,\sum_{i=1}^{k}\Cal{O}_{\Bbb{P}^1}(d_i)\bigr)$
and let $f\:V\to\Bbb{P}^1$ be the natural projection. Then
$\operatorname{Pic}(V)\cong\Bbb{Z} M\oplus\Bbb{Z}L$, where $M$ is the
tautological line bundle on~$V$ and $L$ is the class of a fibre of~$f$.
Let $(t_1:t_2)$ be homogeneous coordinates on the base~$\Bbb{P}^1$, and let
$(x_1:\dots:x_k)$ be the homogeneous coordinates \rom{(}corresponding to
the coordinates on~$\sum_{i=1}^k \Cal O_{\PP^1}(d_i)$\rom{)} on the fibre
of~$f$, which is isomorphic to~$\Bbb{P}^{k-1}$. Then $|aM+bL|$ is generated
by bihomogeneous coordinates
$$
c_{i_{1}, \dots, i_{k}}x_1^{i_{1}} x_2^{i_{2}}\dots x_k^{i_{k}},
$$
where $\sum_{j=1}^{k}i_{j}=a$, \ $i_{j}\geqslant 0$, and
$c_{i_{1},\dots, i_{k}}=c_{i_{1},\dots, i_{k}}(t_{0}:t_{2})$ is a
homogeneous polynomial of degree $b+\sum_{j=1}^{k}i_{j}d_{j}$.
\endproclaim

Proposition~2.19 implies the following result, which is known as a lemma
of~Reid.

\proclaim{Corollary 2.20}
Let $V$ be a $k$-dimensional rational scroll
$\operatorname{Proj}\bigl(\,\sum_{i=1}^{k}\Cal{O}_{\Bbb{P}^1}(d_i)\bigr)$
with $d_{1}\geqslant\dots\geqslant d_{k}\geqslant 0$, and let
$Y_{j}\subset V$ be the \rom{``}negative rational subscroll\rom{''}
$\operatorname{Proj}\bigl(\,
\bopl_{i=j}^{k}\Cal{O}_{\Bbb{P}^1}(d_i)\bigr)$, which corresponds to
the natural projection
$$
\bopl_{i=1}^{k}\Cal{O}_{\Bbb{P}^1}(d_i)\to
\bopl_{i=j}^{k}\Cal{O}_{\Bbb{P}^1}(d_i).
$$
Take an effective divisor $D\subset V$ that is rationally equivalent to
$aM+bL$, where $M$ is the tautological line bundle on~$V$, \ $L$ is
a fibre of the natural projection to~$\Bbb{P}^1$, and $a, b\in\Bbb{Z}$.
We have $\operatorname{mult}_{Y_{j}}(D)\geqslant q$ for $q\in\Bbb{N}$
if and only if~$ad_{j}+b+(d_{1}-d_{j})(q-1)<0$.
\endproclaim

The following result is implied by the Riemann--Roch theorem
(see~\cite{15}, \cite{17}, \cite{88}),
the Kawamata--Viehweg vanishing theorem (see~\cite{91}, \cite{120}),
the rationality of canonical singularities (see~\cite{93}, \cite{96}),
and the global-to-local spectral sequence.

\proclaim{Proposition 2.21}
Let $X$ be a Fano $3$-fold with canonical Gorenstein singularities. Then
$$
h^{0}\bigl(\Cal{O}_{X}(-mK_{X})\bigr)={\frac{m(m+1)(2m+1)}{12}}\,
(-K_{X})^{3}+2m+1.
$$
\endproclaim

The following two results are well known. Their proof can be found in
\cite{8}, \cite{23}, \cite{24}, \cite{9}, \cite{10}, \cite{11}, \cite{14}.
For their modern proof see \cite{103} and~\cite{20}.

\proclaim{Theorem 2.22}
Let $W$ be a smooth minimal geometrically irreducible and geometrically
rational surface defined over a perfect field~$\Bbb{F}$. This means that
no curve on~$W$ can be contracted to a smooth point over~$\Bbb{F}$ and
$W$ is irreducible and rational over~$\overline{\Bbb{F}}$. Then either
$\operatorname{Pic}(W)\cong\Bbb{Z}$ and $W$ is a smooth del~Pezzo surface
or $\operatorname{Pic}(W)\cong\Bbb{Z}\oplus\Bbb{Z}$ and $W$ is a conic
bundle $\pi\:W\to Z$.
\endproclaim

\proclaim{Theorem 2.23}
Let $W$ be a smooth minimal geometrically irreducible and geometrically
rational surface defined over a perfect field~$\Bbb{F}$. The surface $W$ is
rational over $\Bbb{F}$ if and only if~$W$ has an $\Bbb{F}$-point and
$K_{W}^{2}\geqslant 5$.
\endproclaim

\goodbreak
\proclaim{Theorem 2.24 \rm\cite{23}}
Let $W$ be a smooth geometrically irreducible and geometrically
rational surface defined over a $C_{1}$-field $\Bbb{F}$, say, over
$\Bbb{F}=\Bbb{C}(x)$. Then $W$ has an $\Bbb{F}$-point.
\endproclaim

\proclaim{Theorem 2.25 \rm\cite{95}}
Let $Y$ be a projective variety, and let $g\:Y\to R$ be a morphism with a
section onto a smooth curve~$R$. Suppose that there is a set
$\{r_{1},\dots,r_{k}\}\subset R$
of closed points such that each fibre $Y_{i}=g^{-1}(r_{i})$ is
smooth and separably rationally connected. Then for every set of
closed points $y_{i}\in Y_{i}$ there is a section $C\subset Y$
of the morphism $g$ passing through each point $y_{i}$.
\endproclaim

\head
\S\,3. Proof of Theorem~1.5
\endhead

Let $X$ be a Fano 3-fold with canonical Gorenstein singularities
such that the linear system $|-K_{X}|$ has no base points but the induced
morphism $\varphi_{|-K_{X}|}$ is not an embedding.
Then $\varphi_{|-K_{X}|}\:X\to
\allowmathbreak
Y\subset\Bbb{P}^{n}$ is a double covering and
$\operatorname{deg}(Y\subset\Bbb{P}^{n})=n-2$, where
$n=-{\frac{1}{2}}K_{X}^{3}+2$.

\remark{Remark 3.1}
If $-K_{X}^{3}=2$, then the 3-fold $Y$ is nothing but $\Bbb{P}^{3}$ and
$\varphi_{|-K_{X}|}$ is a double covering ramified in a sextic surface
(possibly singular). In this case, $X$ may be regarded as a hypersurface
of degree~$6$ in~$\Bbb{P}(1^{4},3)$. Birational geometry of such varieties
$X$ was studied in~\cite{54}, \cite{37}, \cite{16},
\cite{27}, \cite{111}, \cite{82}, \cite{42},~\cite{59}.
\endremark

\remark{Remark 3.2}
If $-K_{X}^{3}=4$, then $Y$ is a quadric (possibly singular) in
$\Bbb{P}^{4}$ and $\varphi_{|-K_{X}|}$ is a double covering branched over
a surface that is cut on~$Y$ by a quartic hypersurface in~$\Bbb{P}^{4}$.
In this case, $X$ may be regarded as a complete intersection of a quadric
cone and a quartic in~$\Bbb{P}(1^{5},2)$.
Birational geometry of such varieties $X$ was studied in
\cite{54}, \cite{37}, \cite{16}, \cite{27}, \cite{3},~\cite{4}.
\endremark

\medskip
Thus we may assume that $-K_{X}^{3}\geqslant 6$. Hence Theorem~2.11
implies that either $-K_{X}^{3}=8$ and $Y\subset\Bbb{P}^{6}$ is a cone over
a Veronese surface $F_{4}\subset\nomathbreak\Bbb{P}^{5}$ or $Y$ is the image
of a rational scroll $\Bbb{F}(d_{1},d_{2},d_{3})=\operatorname{Proj}\bigl(\,
\bopl_{i=1}^{3}\Cal{O}_{\Bbb{P}^{1}}(d_{i})\bigr)$
under the map given by the tautological line bundle, where
$0\ne d_{1}\geqslant\dots\geqslant d_{3}\geqslant 0$
and $-K_{X}^{3}=2(d_{1}+d_{2}+d_{3})$.

\proclaim{Lemma 3.3}
Suppose that $Y$ is a cone over a Veronese surface $F_{4}$ with vertex~$O$.
Then $X$ is a hypersurface of degree~$6$ in~$\Bbb{P}(1^{3},2,3)$.
\endproclaim

\demo{Proof}
We have $Y\cong\Bbb{P}(1^{3},2)$. The double covering
$\varphi_{|-K_{X}|}$ is branched over the vertex~$O$ because $O$ is not
a Gorenstein point of~$Y$. On the other hand, the equation
$-K_{X}^{3}=8$ implies that the double covering $\varphi_{|-K_{X}|}$
is branched over a divisor $D\subset Y$ such that $D\sim\Cal{O}_{Y}(6)$.
Similarly to the smooth case (see~\cite{15}), it follows (see~\cite{117})
that $X$ is a hypersurface of degree $6$ in~$\Bbb{P}(1^{3},2,3)$.
\enddemo

Hence we may assume that there is a birational morphism $f$: $U\to Y$
for some $U=\operatorname{Proj}\bigl(\,
\bopl_{i=1}^{3}\Cal{O}_{\Bbb{P}^{1}}(d_{i})\bigr)$,
and we have $f=\varphi_{|M|}$, where $M$ is the tautological line bundle
on~$U$, \  $0\ne d_{1}\geqslant\dots\geqslant d_{3}\geqslant 0$ and
$-K_{X}^{3}=2(d_{1}+d_{2}+d_{3})\geqslant 6$.

\proclaim{Lemma 3.4}
Suppose that $d_{2}=d_{3}=0$.
Then $X$ is either a hypersurface of degree $10$ in
$\Bbb{P}(1^{2},3^{2},5)$ or a hypersurface of degree $12$
in~$\Bbb{P}(1^{2},4^{2},6)$.
\endproclaim

\demo{Proof}
Take a sufficiently general divisor $H\in |-K_{X}|$.  Then $H$ is a K3
surface with Du~Val singularities and $f(H)$ is a cone in~$\Bbb{P}^{n-1}$
over a rational normal curve. Moreover, the restriction map
$$
H^{0}\bigl(\Cal{O}_{X}(-K_{X})\bigr)\to
H^{0}\bigl(\Cal{O}_{H}(-K_{X}|_{H})\bigr)
$$
is surjective since $H^{1}(\Cal{O}_{X})=0$. Hence the equations
$d_{1}=d_{2}=0$ imply that $-K_{X}^{3}\leqslant 8$ by~\cite{116}.

Thus there are two possible cases: $d_{1}=3$ and $d_{1}=4$. We have
$Y\cong\Bbb{P}(1^{2},3^{2})$ in the first case and
$Y\cong\Bbb{P}(1^{2},4^{2})$ in the second case. To get the
desired result, we now proceed as in the proof of Lemma~3.3.
\enddemo

\remark{Remark 3.5}
One can use basic properties of hypersurfaces in the weighted projective
spaces (see~\cite{76}) to prove the existence of a hypersurface of degree
$10$ in~$\Bbb{P}(1^{2},3^{2},5)$ and a hypersurface of degree $12$ in
$\Bbb{P}(1^{2},4^{2},6)$ having only canonical Gorenstein singularities.
However, we shall prove this in a different and more geometric way
together with the proof for other possible cases.
\endremark

\medskip
Let $V$ be the normalization of the fibred product
$X\times_{Y}U$, \ $\pi\:V\to U$ the double covering induced by
$\varphi_{|-K_{X}|}\:X\to Y$, and~$h\:V\to X$ the birational morphism induced
by $f\:U\to Y$.

\proclaim{Lemma 3.6}
The $3$-fold $V$ has canonical Gorenstein singularities, the anticanonical
divisor $-K_{V}$ is numerically effective and big, and
$K_{V}\sim h^{*}(K_{X})$, that is, the map $h$ is crepant.
\endproclaim

\demo{Proof}
If $d_{2}\ne 0$, then the 3-folds $X$ and~$V$ are isomorphic in
codimension~$2$, which easily yields the lemma (compare~\cite{92}).

Thus we may assume that $d_{2}=0$. Then $f\:U\to Y$
contracts a divisor $D\subset U$ to a curve $C\cong\Bbb{P}^{1}$, and
Lemma~3.4 implies that either $d_{1}=3$ or $d_{1}=4$. In both cases
$\varphi_{|-K_{X}|}$ must be ramified in the curve~$C$ since $Y$
is non-Gorenstein at a general point of~$C$.

Let $R\subset U$ be the ramification divisor of~$\pi\:V\to U$, \
$M$ the tautological line bundle on~$U$, and $L$ a fibre of the
natural projection of~$U$ to~$\Bbb{P}^{1}$. Then the equivalences
$$
-K_{X}\sim \varphi_{|-K_{X}|}^{*}\bigl(\Cal{O}_{{\Bbb P}^{n}}(1)|_{Y}\bigr),
\qquad
M\sim f^{*}\bigl(\Cal{O}_{{\Bbb P}^{n}}(1)|_{Y}\bigr)
$$
imply that $R\sim 4M-2(d_{1}-2)L+aD$ for some $a\in\Bbb{Z}$. Moreover, we
have $a\geqslant 0$ since the singularities of~$X$ are canonical.
Suppose that $a>0$. Then Corollary~2.20 shows that $R=2D\cup S$, where
$S$ is an effective divisor on~$U$ since $D\sim M-d_{1}L$. This contradicts
the normality of~$V$. Thus $a=0$, \  $R\sim4M-2(d_{1}-\nomathbreak 2)L$,
and $-K_{V}\sim f^{*}(M)\sim h^{*}(K_{X})$, which easily yields the
desired assertion.
\enddemo

\goodbreak
Let $D\subset U$ be the ramification divisor of~$\pi\:V\to U$, \
$M$ the tautological line bundle on~$U$, and $L$ a fibre of the natural
projection of~$U$ to~$\Bbb{P}^{1}$. Then $-K_{V}\sim
\pi^{*}(M)$ by construction. Hence,
$$
D\sim 4M-2(d_{1}+d_{2}+d_{3}-2)L.
$$

\goodbreak
Let $Y_{2}\subset V$ and $Y_{3}\subset V$ be the subscrolls corresponding to
the natural projections
$$
\bopl_{i=1}^{3}\Cal{O}_{\Bbb{P}^1}(d_i)\to
\Cal{O}_{\Bbb{P}^1}(d_2)\oplus\Cal{O}_{\Bbb{P}^1}(d_3),
\qquad
\bopl_{i=1}^{3}\Cal{O}_{\Bbb{P}^1}(d_i)\to\Cal{O}_{\Bbb{P}^1}(d_3).
$$
Then $Y_{2}\cong\operatorname{Proj}(\Cal{O}_{\Bbb{P}^1}(d_2)\oplus
\Cal{O}_{\Bbb{P}^1}(d_3))$ and~$Y_{3}\cong\Bbb{P}^{1}$.

\proclaim{Lemma 3.7}
We have $\operatorname{mult}_{Y_{2}}(D)\leqslant 1$ and
$\operatorname{mult}_{Y_{3}}(D)\leqslant 3$.
\endproclaim

\demo{Proof}
The first inequality follows from the normality of~$V$. Suppose that
$d=\operatorname{mult}_{Y_{3}}(D)\geqslant 2$. Then the local equation
of $V$ in the neighbourhood of a generic point of the curve
$C=\pi^{-1}(Y_{3})$ is
$$
\omega^{2}=f_{d}(x,y)+f_{d+1}(x,y)+\dots
\subset\operatorname{Spec}(\Bbb{C}[x,y,z,\omega]),
$$
where $x=y=0$ are local equations of the curve~$C$, and $f_{i}(x,y)$
is a homogeneous polynomial of degree~$i$. On the other hand, the
singularities of~$V$ at the general point of~$C$ must be locally isomorphic
to one of the following types of singularities:
$\Bbb{C}\times\Bbb{A}_{n}$, $\Bbb{C}\times\Bbb{D}_{n}$,
$\Bbb{C}\times\Bbb{E}_{6}$, $\Bbb{C}\times\Bbb{E}_{7}$,
$\Bbb{C}\times\Bbb{E}_{8}$. Then Theorem~2.9 implies that $d\leqslant 3$.
\enddemo

Therefore Corollary~2.20 implies that
$$
d_1-d_3-2d_2+4\geqslant 0, \qquad d_2-d_1-2d_3+4\geqslant0.
$$
We also have $0\ne d_{1}\geqslant\dots\geqslant d_{3}\geqslant 0$ and
$d_{1}+d_{2}+d_{3}\geqslant 3$ by assumption. The resulting inequalities
determine $44$ different rational scrolls
$$
\Bbb{F}(d_{1},d_{2},d_{3})=\operatorname{Proj}\biggl(\,
\bopl_{i=1}^{3}\Cal{O}_{\Bbb{P}^{1}}(d_{i})\biggr)
$$
with the ramification divisor
$$
D\sim 4M-2(d_{1}+d_{2}+d_{3}-2)L,
$$
where $M$ is the tautological line bundle on
$\Bbb{F}(d_{1},d_{2},d_{3})$ and $L$ is a fibre of the natural projection
to~$\Bbb{P}^{1}$.

\remark{Remark 3.8}
The 3-fold $X$ is an anticanonical model of the 3-fold $V$, that is,
$X\cong\varphi_{|-rK_{V}|}(V)$ for $r\gg 0$. Thus the contraction theorem
(see~\cite{93}) implies that $X$ is a Fano 3-fold with canonical Gorenstein
singularities if and only if~$V$ has canonical Gorenstein singularities
and $-K_{V}$ is numerically effective and big. On the other hand,
$V$ is uniquely determined by the rational scroll
$\Bbb{F}(d_{1},d_{2},d_{3})$ and the ramification divisor
$D\in|4M-2(d_{1}+d_{2}+d_{3}-2)L|$. The only trouble is that the linear
system $|4M-2(d_{1}+d_{2}+d_{3}-2)L|$ may contain no divisor~$D$ such that
the corresponding double covering~$V$ has canonical singularities.
\endremark

\goodbreak
\medskip
In the rest of this section we explicitly show that, in each of the
cases obtained, the linear system $|4M-2(d_{1}+d_{2}+d_{3}-2)L|$ contains
a divisor $D$ such that $V$ has canonical singularities. We shall use
Corollary~2.7 together with Proposition~2.19. This verification will
complete the proof of Theorem~1.5.

\remark{Remark 3.9}
The same ideas were actually used in the classification of smooth
hyperelliptic Fano 3-folds (see~\cite{15}, \cite{17}). In the smooth case,
the corresponding inequalities are much stronger and the calculations are
much shorter. This method was also used in~\cite{41} to find an effective
bound of the degree of hyperelliptic Fano 3-folds with canonical Gorenstein
singularities, but there is a gap in the proof of Lemma~3.2 in~\cite{41}.
Namely, the stronger inequality $\operatorname{mult}_{Y_{3}}(D)\leqslant 2$
was used there instead of the inequality
$\operatorname{mult}_{Y_{3}}(D)\leqslant 3$. This gave a wrong bound
$-K_{X}^{3}\leqslant 16$ instead of the right one
$-K_{X}^{3}\leqslant 40$, which is {\it a posteriori\/} seen to be sharp.
Nevertheless, one can use the estimate $-K_{X}^{3}\leqslant 40$ to prove
the main result of~\cite{41}. But this result became obsolete now because
of~\cite{110}.
\endremark

\medskip
Let us consider one of the possible cases in full detail.

\example{Example 3.10}
Let $\pi\:V\to \Bbb{F}(6,2,0)$ be a double covering branched
over a sufficiently general divisor $D\subset\Bbb{F}(6,2,0)$
such that $D\sim 4M-12L$, where $M$ is the class of the tautological
line bundle on~$\Bbb{F}(6,2,0)$ and $L$ is the class of a fibre of the
projection of~$\Bbb{F}(6,2,0)$ to~$\Bbb{P}^{1}$. We must show that the
3-fold $V$ has canonical singularities.

By Proposition~2.19, the divisor
$D$ is given by the zeros of the bihomogeneous polynomial
$$
\align
&\alpha_{12}(t_1,t_2)\,x_1^4+\alpha_8(t_1,t_2)\,x_1^3x_2
+\alpha_{6}(t_1,t_2)\,x_1^3x_3+\alpha_4(t_1,t_2)\,x_1^2x_2^2
\\
&\qquad\qquad
+\alpha_2(t_1,t_2)\,x_1^2x_2x_3+ \alpha_0^1(t_1,t_2)\,x_1^2x_3^2+
\alpha_0^2(t_1,t_2)\,x_1x_2^3,
\endalign
$$
where $\alpha_d(t_1,t_2)$ (or $\alpha_d^i(t_1,t_2)$) is an arbitrary form
of degree~$d$. We define a surface $E\subset\Bbb{F}(6,2,0)$
and  a curve $C\subset\Bbb{F}(6,2,0)$ by the equations $x_1=0$
and $x_1=x_2=0$ respectively. The base locus of the linear system
$|4M-12L|$ equals~$E$. (In particular, $D\setminus E$ and
$V\setminus\pi^{-1}(E)$ are smooth by the Bertini theorem.)
The automorphism group of~$E\cong\FF(4,0)$ acts transitively
on~$E\setminus C$, that is, all points of~$E\setminus C$ are
mapped to each other by changes of the coordinates
$t_1$, $t_2$, $x_4$, $x_5$. By Lemma~3.7, the divisor $D$ has
multiplicity~$1$ at a general point of~$E$. Hence for every point of
$E\setminus C$ there is a divisor $D'$ whose multiplicity at this point
is equal to~$1$: it suffices to make an appropriate change of coordinates
in the equation of~$D$. The singularities of a general divisor
$D$ on~$E\setminus C$ are canonical by Corollary~2.7, and it suffices
to prove that for every point $p$ of~$C$ there is a divisor~$D$
such that the corresponding variety $V$ has a canonical singularity in the
neighbourhood of the point $\pi^{-1}(p)\in\pi^{-1}(C)$.

Let $Y$ be a fibre of the projection of~$\Bbb{F}(6,2,0)$ to
$\Bbb{P}^{1}$ over a sufficiently general point $P\in\Bbb{P}^{1}$.
We put $Z=\pi^{-1}(Y)$. Then $Z$ is a del~Pezzo surface of degree~$2$.
Moreover, the only possible singular point of~$Z$ is $O=\pi^{-1}(C\cap Y)$.
Let us prove that $O$ is a Du~Val point on~$Z$. This already implies that
the singularities of~$Y$ are canonical.

\goodbreak
Suppose that the point $P\in\Bbb{P}^{1}$ has homogeneous coordinates
$(\gamma:\delta)$. Then $Z$ may be given as a hypersurface
$$
\omega^{2}=\alpha_{12}x_1^4+\alpha_8x_1^3x_2
+\alpha_{6}x_1^3x_3+\alpha_4x_1^2x_2^2+\alpha_{2}x_1^2x_2x_3
+\alpha_{0}^1x_1^2x_3^2+\alpha_{0}^2x_1x_2^3
$$
in
$\Bbb{P}(1,1,1,2)\cong\operatorname{Proj}(\Bbb{C}[x_{1},x_{2},x_{3},\omega])$,
where $\alpha_d^{i}=\alpha_d^{i}(\gamma,\delta)$ (respectively,
$\alpha_d=\alpha_d(\gamma,\delta)$). Since $Y$ is general, we have
$\alpha_d^{i}\ne 0$ (resp\. $\alpha_d\ne 0$) for all~$d$ and~$i$.
Therefore we may assume for convenience that $\alpha_d^{i}=1$
(resp\. $\alpha_d=1$) for all~$d$ and~$i$.

Let $\omega=x$,  $x_1=y$,  $x_2=z$ and~$x_3=1$. Then the local equation of
$Z$ in a neighbourhood of~$O$ is given by
$$
x^2+y^4+y^3z+y^3+y^2z^2+y^2z+y^2+yz^3=0.
$$
Let $\operatorname{wt}(x)=3$, \ $\operatorname{wt}(y)=3$
and $\operatorname{wt}(z)=1$.
Then $\operatorname{wt}(x^2+y^2+yz^3)=6$, \
$\operatorname{wt}(y^4)=12$, \
$\operatorname{wt}(y^3z)=10$, \ $\operatorname{wt}(y^3)=9$, \
$\operatorname{wt}(y^2z^2)=8$ and~$\operatorname{wt}(y^2z)=7$. Moreover,
the equation $x^2+y^2+yz^3=0$ determines an isolated point. Hence
Theorem~2.10 implies that the singularity of~$Z$ at~$O$ is locally
isomorphic to a Du~Val singularity of type~$\Bbb{A}_5$.
Hence the 3-fold $V$ has a singularity of type $\Bbb{A}_5\times\Bbb{C}$
at a general point of the 
\linebreak
curve $\pi^{-1}(C)$.

Using the generality in the choice of~$D$, we may actually assume
that the point $P\in\Bbb{P}^{1}$ is not just a general point but an
arbitrary point of~$\Bbb{P}^{1}$. In other words, given any point
$P\in\Bbb{P}^{1}$, one can find homogeneous polynomials $\alpha_d$
such that $\alpha_d(P)\ne 0$ and repeat all the previous arguments
in the neighbourhood of the corresponding point $O=\pi^{-1}(C\cap Y)$.
Hence the singularities of~$V$ are canonical 
\linebreak
by Corollary~2.7.
\endexample

In the rest of the section we consider the other possible cases
following the pattern of Example~3.10. The differences appear only
in the numerical characteristics of varieties, their equations,
types of singularities etc. They are surveyed 
\linebreak
in Table~1.

\looseness=1
Table~1 is organized as follows. The first column contains labels of
the varieties $V$ in the notation of Theorem~1.5. The second column
yields a triple $(d_1, d_2, d_3)$ such that there is a double covering
$\pi\:V\to\FF(d_1, d_2, d_3)$, which is branched over a divisor~$D$.
The third column displays the number~$b$ such that $|D|=|4M+bL|$.
The corresponding linear system appears to be base point free in the
cases $H_{4}$, $H_{6}$ and~$H_{9}$. Then the divisor~$D$ is non-singular
by the Bertini theorem and, therefore, $V$ is non-singular (and we do not
need the information of the other columns). In all other cases, the set
$\operatorname{Bs}|D|$ of base points is either the curve $C=Y_3$
given by~$x_1=x_2=0$ (then the Bertini theorem shows that $D$ and $V$ are
smooth outside $C$ and~$\pi^{-1}(C)$ respectively, so it suffices to study
the singularities of~$V$ over~$C$ only)
or the surface $E=Y_2$ given by~$x_1=0$
(then the divisor $D$ has multiplicity~1 at a general point of~$E$
and, since the automorphism group of~$E\cong\FF(d_2, d_3)$ is transitive
on $E\setminus C$ (compare Example~3.10), Corollary~2.7 shows that it
again suffices to study the singularities of~$V$ over~$C$ only).

\looseness=1
The fourth column contains equations of general divisors
$D$ in the linear system $|4M+bL|$, and the fifth column yields
an equation of the fibre $Z$ of the projection $V\to\PP^1$ over
a general point of~$\PP^1$ in the neighbourhood of a general point
of $\pi^{-1}(C)$ after a change of coordinates
$\omega=x$, \ $x_1=y$, \ $x_2=z$, \ $x_3=1$ (see Example~3.10).
The same equation locally determines $V$ if we regard it as an
equation in~$t$,~$x$,~$y$,~$z$. The corresponding point appears to be
non-singular in the cases $H_{5}$, $H_{12}$ and~$H_{17}$.
In the other cases we attribute new weights
$\wt(x)=\wwt_x$, \ $\wt(y)=\wwt_y$, \ $\wt(z)=\wwt_z$ (listed in the
sixth column) to the variables  $x,y,z$ and find out that the terms
of the lowest weight determine an isolated Du~Val singularity.
We notice that the weights $\wt(x)$, $\wt(y)$, $\wt(z)$ coincide
with the weights of a Du~Val singularity indicated in the
seventh column. Hence the singularity of~$Z$ in the chosen
neighbourhood is Du~Val of this type (by Theorem~2.10), and the
singularity of~$V$ is locally isomorphic to the product of~$\CC$
and the corresponding Du~Val singularity. In any case, $V$ has canonical
singularities in this neighbourhood and, by Corollary~2.7, 
all singularities of $V$ are 
\linebreak
canonical.

\medskip

{

\raggedbottom

$$
\vbox{
\let\PLUS=+
\binoppenalty10000
\relpenalty10000
\catcode`\+=\active
\raggedright
\def+{\allowbreak\PLUS\nobreak\unskip}
\define\mystrut{\vphantom{$\Bigl($}}
\def\ensp{\kern0.41em}
\tabskip=0pt \offinterlineskip
\halign{\vrule\mystrut
\ensp\hfill#\hfill\ensp&\vrule# &
\ensp\hfill#\hfill\ensp&\vrule# &
\ensp\hfill#\hfill\ensp&\vrule# &
\ensp\hfill\tabboxone{#}\hfill\ensp&\vrule# &
\ensp\hfill\tabboxtwo{#}\hfill\ensp&\vrule# &
\ensp\hfill\tabboxthree{#}\hfill\ensp&\vrule# &
\ensp\hfill\tabboxthree{\hsize 1.2cm {#}}\hfill\ensp\vrule
\cr
\multispan{13}
\hfill{\smc Table 1}
\cr
\noalign{\medskip}
\noalign{\hrule}
\lw{2.5}{$H_i$}&&\lw{2.5}{~$(d_1{,}d_2{,}d_3)$}&&\lw{2.5}{~~$b$}
&&
\centerline{\lw{2.5}{Equation of~$D$}} &&
Local equation of~$V$
&&\centerline{\lw{2.5}{Weights}}&&Singu\-larity
\cr
\noalign{\hrule}
$H_{4}$ && $(1,1,1)$ && $-2$ && \centerline{--} &&
\centerline{--} && \centerline{--} && \centerline{--}
\cr
\noalign{\hrule}
$H_{5}$ && $(2,1,0)$ && $-2$ &&
$\alpha_6x_1^4+\alpha_5x_1^3x_2+\alpha^1_4x_1^3x_3
+\alpha^2_4x_1^2x_2^2+\alpha^1_3x_1^2x_2x_3+\alpha^1_2x_1^2x_3^2+
\alpha^2_3x_1x_2^3+\alpha^2_2x_1x_2^2x_3
+\alpha^1_1x_1x_2x_3^2+\alpha^1_0x_1x_3^3+\alpha^3_2x_2^4
+\alpha^2_1x_2^3x_3+\alpha^2_0x_2^2x_3^2=\nomathbreak0$
&&
$x^2+y^3z+y^3+y^2z^2+y^2z+y^2+yz^3+yz^2+yz+y+z^4+z^3+z^2=\nomathbreak0$
&& \centerline{--} && Non-sin\-gu\-lar point
\cr
\noalign{\hrule}
$H_{6}$ && $(2,1,1)$ && $-4$ && \centerline{--} &&
\centerline{--} && \centerline{--} && \centerline{--}
\cr
\noalign{\hrule}
$H_{7}$ && $(2,2,0)$ && $-4$ &&
$\alpha^1_4x_1^4+\alpha^2_4x_1^3x_2+\alpha^3_4x_1^2x_2^2
+\alpha^4_4x_1x_2^3+\alpha^5_4x_2^4+\alpha^1_2x_1^3x_3+
\!\!\alpha^2_2x_1^2x_2x_3+\alpha^3_2x_1x_2^2x_3+\alpha^4_2x_2^3x_3
+\alpha^1_0x_1^2x_3^2+\alpha^2_0x_1x_2x_3^2
+\alpha^3_0x_2^2x_3^2=\nomathbreak0$
&&
$x^2+P_2(y,z)+P_3(y,z)+P_4(y,z)=\nomathbreak0$
\smash{($P_i$} is a homogeneous polynomial of degree~$i$) &&
$\wwt_x=\nomathbreak1$ \ $\wwt_y=\nomathbreak1$ \ $\wwt_z=\nomathbreak1$ &&
\centerline{$\Bbb{A}_1$}
\cr
\noalign{\hrule}
}}
$$
$$
\vbox{
\let\PLUS=+
\binoppenalty10000
\relpenalty10000
\catcode`\+=\active
\raggedright
\def+{\allowbreak\PLUS\nobreak\unskip}
\define\mystrut{\vphantom{$\Bigl($}}
\def\ensp{\kern0.41em}
\tabskip=0pt \offinterlineskip
\halign{\vrule\mystrut
\ensp\hfill#\hfill\ensp&\vrule# &
\ensp\hfill#\hfill\ensp&\vrule# &
\ensp\hfill#\hfill\ensp&\vrule# &
\ensp\hfill\tabboxone{#}\hfill\ensp&\vrule# &
\ensp\hfill\tabboxtwo{#}\hfill\ensp&\vrule# &
\ensp\hfill\tabboxthree{#}\hfill\ensp&\vrule# &
\ensp\hfill\tabboxthree{\hsize 1.2cm #}\hfill\ensp\vrule
\cr
\multispan{13}
\hfill{\smc Table 1 continued}
\cr
\noalign{\medskip}
\noalign{\hrule}
$H_{8}$ && $(2,2,1)$ && $-6$ &&
$\alpha^1_2x_1^4+\alpha^2_2x_1^3x_2+\alpha^3_2x_1^2x_2^2+
\alpha^2_4x_1x_2^3+\alpha^5_2x_2^4+\alpha^1_1x_1^3x_3+
\!\!\alpha^2_1x_1^2x_2x_3+\alpha^3_1x_1x_2^2x_3+\alpha^4_1x_2^3x_3+
\alpha^1_0x_1^2x_3^2+\alpha^2_0x_1x_2x_3^2+\alpha^3_0x_2^2x_3^2=\nomathbreak0$
&&
$x^2+P_2(y,z)+P_3(y,z)+P_4(y,z)=\nomathbreak0$ \vp\smash{($P_i$} is a
homogeneous polynomial of degree~$i$) &&
$\wwt_x=\nomathbreak1$ $\wwt_y=\nomathbreak1$ $\wwt_z=\nomathbreak1$ &&
\centerline{$\Bbb{A}_1$}
\cr
\noalign{\hrule}
$H_{9}$ && $(2,2,2)$ && $-8$ && \centerline{--} &&
\centerline{--} && \centerline{--} && \centerline{--}
\cr
\noalign{\hrule}
$H_{10}$ && $(3,0,0)$ && $-2$ &&
$\alpha_{10}x_1^4+\alpha^1_7x_1^3x_2+\alpha^2_7x_1^3x_3+
\alpha^1_4x_1^2x_2^2+\alpha^2_4x_1^2x_2x_3+
\alpha^3_4x_1^2x_3^2+\alpha^1_1x_1x_2^3+
\alpha^2_1x_1x_2^2x_3+
\!\!\alpha^3_1x_1x_2x_3^2+\alpha^4_1x_1x_3^3=\nomathbreak0$
&& $x^2+y^3z+y^3+y^2z^2+y^2z+y^2+yz^3+yz^2+yz+y=\nomathbreak0$ &&
\centerline{--} && Non-sin\-gu\-lar point
\cr
\noalign{\hrule}
$H_{11}$ && $(3,1,0)$ && $-4$ &&
$\alpha_8x_1^4+\alpha_6x_1^3x_2+\alpha_5x_1^3x_3+\alpha_4x_1^2x_2^2+
\alpha_3x_1^2x_2x_3+\alpha_2^1x_1^2x_3^2+
\!\!\alpha_0^1x_1x_2x_3^2\phantom{aa}+
\alpha_2^2x_1x_2^3+\alpha_1x_1x_2^2x_3+\alpha_0^2x_2^4=\nomathbreak0$ &&
$x^2+y^4+y^3z+y^3+y^2z^2+y^2z+y^2+yz^3+yz^2+yz+z^4=\nomathbreak0$ &&
$\wwt_x=\nomathbreak1$ $\wwt_y=\nomathbreak1$ $\wwt_z=\nomathbreak1$ &&
\centerline{$\Bbb{A}_1$}
\cr
\noalign{\hrule}
$H_{12}$ && $(3,1,1)$ && $-6$ &&
$\alpha_6x_1^4+\alpha^1_4x_1^3x_2+\alpha^2_4x_1^3x_3+
\alpha^1_2x_1^2x_2^2+\alpha^2_2x_1^2x_2x_3+
\alpha^3_2x_1^2x_3^2+\alpha^1_0x_1x_2^3+\alpha^2_0x_1x_2^2x_3+
\!\!\alpha^3_0x_1x_2x_3^2+\alpha^4_0x_1x_3^3=\nomathbreak0$ &&
$x^2+y^3z+y^3+y^2z^2+y^2z+y^2+yz^3+yz^2+yz+y=\nomathbreak0$ &&
\centerline{--} && Non-sin\-gu\-lar point
\cr
\noalign{\hrule}
$H_{13}$ && $(3,2,0)$ && $-6$ &&
$\alpha_6x_1^4+\alpha_5x_1^3x_2+\alpha_3^1x_1^3x_3+\alpha_4x_1^2x_2^2+
\alpha_2^1x_1^2x_2x_3+\alpha_0^1x_1^2x_3^2+\alpha_3^2x_1x_2^3+
\alpha_1x_1x_2^2x_3+\alpha_2^2x_2^4+\alpha_0^2x_2^3x_3=\nomathbreak0$ &&
$x^2+y^4+y^3z+y^3+y^2z^2+y^2z+y^2+yz^3+yz^2+z^4+z^3=\nomathbreak0$ &&
$\wwt_x=\nomathbreak3$ $\wwt_y=\nomathbreak3$ $\wwt_z=\nomathbreak2$ &&
\centerline{$\Bbb{A}_2$}
\cr
\noalign{\hrule}
}}
$$
$$
\vbox{
\let\PLUS=+
\binoppenalty10000
\relpenalty10000
\catcode`\+=\active
\raggedright
\def+{\allowbreak\PLUS\nobreak\unskip}
\define\mystrut{\vphantom{$\Bigl($}}
\def\ensp{\kern0.41em}
\tabskip=0pt \offinterlineskip
\halign{\vrule\mystrut
\ensp\hfill#\hfill\ensp&\vrule# &
\ensp\hfill#\hfill\ensp&\vrule# &
\ensp\hfill#\hfill\ensp&\vrule# &
\ensp\hfill\tabboxone{#}\hfill\ensp&\vrule# &
\ensp\hfill\tabboxtwo{#}\hfill\ensp&\vrule# &
\ensp\hfill\tabboxthree{#}\hfill\ensp&\vrule# &
\ensp\hfill\tabboxthree{\hsize 1.2cm #}\hfill\ensp\vrule
\cr
\multispan{13}
\hfill{\smc Table 1 continued}
\cr
\noalign{\medskip}
\noalign{\hrule}
$H_{14}$ && $(3,2,1)$ && $-8$ &&
$\alpha_4x_1^4+\alpha_3x_1^3x_2+\alpha_2^1x_1^3x_3+\alpha_2^2x_1^2x_2^2+
\alpha_1^1x_1^2x_2x_3+\alpha_0^1x_1^2x_3^2+\alpha_1^2x_1x_2^3+
\alpha_0^2x_1x_2^2x_3+\alpha_0^3x_2^4=\nomathbreak0$ &&
$x^2+y^4+y^3z+y^3+y^2z^2+y^2z+y^2+yz^3+yz^2+z^4=\nomathbreak0$ &&
$\wwt_x=\nomathbreak2$ $\wwt_y=\nomathbreak2$ $\wwt_z=\nomathbreak1$ &&
\centerline{$\Bbb{A}_3$}
\cr
\noalign{\hrule}
$H_{15}$ && $(3,3,0)$ && $-8$ &&
$\alpha^1_4x_1^4+\alpha^2_4x_1^3x_2+\alpha^3_4x_1^2x_2^2+
\alpha^4_4x_1x_2^3+\alpha^5_4x_2^4+\alpha^1_1x_1^3x_3+
\!\!\alpha^2_1x_1^2x_2x_3+
\!\!\alpha^3_1x_1x_2^2x_3+\alpha^4_1x_2^3x_3=\nomathbreak0$
&&
$x^2+P_3(y,z)+P_4(y,z)=\nomathbreak0$ \vp\smash{($P_i$}
is a homogeneous polynomial of degree~$i$) &&
$\wwt_x=\nomathbreak3$ $\wwt_y=\nomathbreak2$ $\wwt_z=\nomathbreak2$ &&
\centerline{$\Bbb{D}_4$}
\cr
\noalign{\hrule}
$H_{16}$ && $(3,3,1)$ && $-10$ &&
$\alpha^1_2x_1^4+\alpha^2_2x_1^3x_2+\alpha^3_2x_1^2x_2^2+
\alpha^4_2x_1x_2^3+\alpha^5_2x_2^4+
\alpha^1_0x_1^3x_3+
\!\!\alpha^2_0x_1^2x_2x_3+
\!\!\alpha^3_0x_1x_2^2x_3+\alpha^4_0x_2^3x_3=\nomathbreak0$ &&
$x^2+P_3(y,z)+P_4(y,z)=\nomathbreak0$ \vp\smash{($P_i$}
is a homogeneous polynomial of degree~$i$) &&
$\wwt_x=\nomathbreak3$ $\wwt_y=\nomathbreak2$ $\wwt_z=\nomathbreak2$ &&
\centerline{$\Bbb{D}_4$}
\cr
\noalign{\hrule}
$H_{17}$ && $(4,0,0)$ && $-4$ &&
$\alpha_{12}x_1^4+\alpha^1_8x_1^3x_2+\alpha^2_8x_1^3x_3+
\alpha^1_4x_1^2x_2^2+\alpha^2_4x_1^2x_2x_3+
\alpha^3_4x_1^2x_3^2+
\alpha^1_0x_1x_2^3+\alpha^2_0x_1x_2^2x_3+
\!\!\alpha^3_0x_1x_2x_3^2+\alpha^4_0x_1x_3^3=\nomathbreak0$ &&
$x^2+y^3z+y^3+y^2z^2+y^2z+y^2+yz^3+yz^2+yz+y=\nomathbreak0$ &&
\centerline{--} && Non-sin\-gu\-lar point
\cr
\noalign{\hrule}
$H_{18}$ && $(4,1,0)$ && $-6$ &&
$\alpha_{10}x_1^4+\alpha_8x_1^3x_2+\alpha_6x_1^3x_3+\alpha_4x_1^2x_2^2+
\!\!\alpha_3x_1^2x_2x_3+\alpha_2x_1^2x_3^2+\alpha_1x_1x_2^3+
\alpha_0x_1x_2^2x_3=\nomathbreak0$ &&
$x^2+y^4+y^3z+y^3+y^2z^2+y^2z+y^2+yz^3+yz^2=\nomathbreak 0$ &&
$\wwt_x=\nomathbreak2$ $\wwt_y=\nomathbreak2$ $\wwt_z=\nomathbreak1$ &&
\centerline{$\Bbb{A}_3$}
\cr
\noalign{\hrule}
$H_{19}$ && $(4,2,0)$ && $-8$ &&
$\alpha_8x_1^4+\alpha_6x_1^3x_2+\alpha_4x_1^3x_3+\alpha_4x_1^2x_2^2+
\alpha_2x_1^2x_2x_3+\alpha_0x_1^2x_3^2+\alpha_2x_1x_2^3+
\alpha_0x_1x_2^2x_3+\alpha_0x_2^4=\nomathbreak0$ &&
$x^2+y^4+y^3z+y^3+y^2z^2+y^2z+y^2+yz^3+yz^2+z^4=\nomathbreak0$&&
$\wwt_x=\nomathbreak2$ $\wwt_y=\nomathbreak2$ $\wwt_z=\nomathbreak1$ &&
\centerline{$\Bbb{A}_3$}
\cr
\noalign{\hrule}
}}
$$
$$
\vbox{
\let\PLUS=+
\binoppenalty10000
\relpenalty10000
\catcode`\+=\active
\raggedright
\def+{\allowbreak\PLUS\nobreak\unskip}
\define\mystrut{\vphantom{$\Bigl($}}
\def\ensp{\kern0.41em}
\tabskip=0pt \offinterlineskip
\halign{\vrule\mystrut
\ensp\hfill#\hfill\ensp&\vrule# &
\ensp\hfill#\hfill\ensp&\vrule# &
\ensp\hfill#\hfill\ensp&\vrule# &
\ensp\hfill\tabboxone{#}\hfill\ensp&\vrule# &
\ensp\hfill\tabboxtwo{#}\hfill\ensp&\vrule# &
\ensp\hfill\tabboxthree{#}\hfill\ensp&\vrule# &
\ensp\hfill\tabboxthree{\hsize 1.2cm #}\hfill\ensp\vrule
\cr
\multispan{13}
\hfill{\smc Table 1 continued}
\cr
\noalign{\medskip}
\noalign{\hrule}
$H_{20}$ && $(4,2,1)$ && $-10$ &&
$\alpha_6x_1^4+\alpha_4x_1^3x_2+\alpha_3x_1^3x_3+\alpha_2x_1^2x_2^2+
\alpha_1x_1^2x_2x_3+\alpha^1_0x_1x_2^3+\alpha^2_0x_1^2x_3^2=\nomathbreak0$ &&
$x^2+y^4+y^3z+y^3+y^2z^2+y^2z+yz^3+y^2=\nomathbreak0$ &&
$\wwt_x=\nomathbreak3$ $\wwt_y=\nomathbreak3$ $\wwt_z=\nomathbreak1$ &&
\centerline{$\Bbb{A}_5$}
\cr
\noalign{\hrule}
$H_{21}$ && $(4,3,0)$ && $-10$ &&
$\alpha_6x_1^4+\alpha_5x_1^3x_2+\alpha_4x_1^2x_2^2+
\alpha_3x_1x_2^3+\alpha^1_2x_2^4+
\alpha^2_2x_1^3x_3+
\!\!\alpha_1x_1^2x_2x_3+\alpha_0x_1x_2^2x_3=\nomathbreak0$ &&
$x^2+y^4+y^3z+y^2z^2+yz^3+z^4+y^3+y^2z+yz^2=\nomathbreak 0$ &&
$\wwt_x=\nomathbreak3$ $\wwt_y=\nomathbreak2$ $\wwt_z=\nomathbreak2$ &&
\centerline{$\Bbb{D}_4$}
\cr
\noalign{\hrule}
$H_{22}$ && $(4,3,1)$ && $-12$ &&
$\alpha_4x_1^4+\alpha_3^1x_1^3x_2+\alpha^1_1x_1^3x_3+\alpha_2x_1^2x_2^2+
\alpha^1_0x_1^2x_2x_3+\alpha^2_1x_1x_2^3+\alpha^2_0x_2^4=\nomathbreak0$ &&
$x^2+y^4+y^3z+y^3+y^2z^2+y^2z+yz^3+z^4=\nomathbreak0$ &&
$\wwt_x=\nomathbreak4$ $\wwt_y=\nomathbreak3$ $\wwt_z=\nomathbreak2$ &&
\centerline{$\Bbb{D}_5$}
\cr
\noalign{\hrule}
$H_{23}$ && $(4,4,0)$ && $-12$ &&
$\alpha^1_4x_1^4+\alpha^2_4x_1^3x_2+\alpha^3_4x_1^2x_2^2+
\alpha^4_4x_1x_2^3+\alpha^5_4x_2^4+\alpha^1_0x_1^3x_3+
\!\!\alpha^2_0x_1^2x_2x_3+
\!\!\alpha^3_0x_1x_2^2x_3+
\alpha^4_0x_2^3x_3=\nomathbreak0$ &&
$x^2+P_3(y,z)+P_4(y,z)=\nomathbreak0$ \vp\smash{($P_i$}
is a homogeneous polynomial of degree~$i$) &&
$\wwt_x=\nomathbreak3$ $\wwt_y=\nomathbreak2$ $\wwt_z=\nomathbreak2$ &&
\centerline{$\Bbb{D}_4$}
\cr
\noalign{\hrule}
$H_{24}$ && $(5,1,0)$ && $-8$ &&
$\alpha_{12}x_1^4+\alpha_8x_1^3x_2+\alpha_7x_1^3x_3+\alpha_4x_1^2x_2^2+
\alpha_3x_1^2x_2x_3+\alpha_0x_1x_2^3+\alpha_2x_1^2x_3^2=\nomathbreak0$ &&
$x^2+y^4+y^3z+y^3+y^2z^2+y^2z+yz^3+y^2=\nomathbreak0$ &&
$\wwt_x=\nomathbreak3$ $\wwt_y=\nomathbreak3$ $\wwt_z=\nomathbreak1$ &&
\centerline{$\Bbb{A}_5$}
\cr
\noalign{\hrule}
$H_{25}$ && $(5,2,0)$ && $-10$ &&
$\alpha_{10}x_1^4+\alpha_7x_1^3x_2+\alpha_5x_1^3x_3+\alpha_4x_1^2x_2^2+
\alpha^1_0x_1^2x_2x_3+\alpha_1x_1x_2^3+\alpha^2_0x_1^2x_3^2=\nomathbreak0$ &&
$x^2+y^4+y^3z+y^3+y^2z^2+y^2z+yz^3+y^2=\nomathbreak0$ &&
$\wwt_x=\nomathbreak3$ $\wwt_y=\nomathbreak3$ $\wwt_z=\nomathbreak1$ &&
\centerline{$\Bbb{A}_5$}
\cr
\noalign{\hrule}
$H_{26}$ && $(5,3,0)$ && $-12$ &&
$\alpha_8x_1^4+\alpha_6x_1^3x_2+\alpha_3x_1^3x_3+\alpha_4x_1^2x_2^2+
\alpha^1_1x_1^2x_2x_3+\alpha^2_1x_1x_2^3+\alpha_0x_2^4=\nomathbreak0$ &&
$x^2+y^4+y^3z+y^3+y^2z^2+y^2z+yz^3+z^4=\nomathbreak0$ &&
$\wwt_x=\nomathbreak4$ $\wwt_y=\nomathbreak3$ $\wwt_z=\nomathbreak2$ &&
\centerline{$\Bbb{D}_5$}
\cr
\noalign{\hrule}
}}
$$
$$
\vbox{
\let\PLUS=+
\binoppenalty10000
\relpenalty10000
\catcode`\+=\active
\raggedright
\def+{\allowbreak\PLUS\nobreak\unskip}
\define\mystrut{\vphantom{$\Bigl($}}
\def\ensp{\kern0.41em}
\tabskip=0pt \offinterlineskip
\halign{\vrule\mystrut
\ensp\hfill#\hfill\ensp&\vrule# &
\ensp\hfill#\hfill\ensp&\vrule# &
\ensp\hfill#\hfill\ensp&\vrule# &
\ensp\hfill\tabboxone{#}\hfill\ensp&\vrule# &
\ensp\hfill\tabboxtwo{#}\hfill\ensp&\vrule# &
\ensp\hfill\tabboxthree{#}\hfill\ensp&\vrule# &
\ensp\hfill\tabboxthree{\hsize 1.2cm #}\hfill\ensp\vrule
\cr
\multispan{13}
\hfill{\smc Table 1 continued}
\cr
\noalign{\medskip}
\noalign{\hrule}
$H_{27}$ && $(5,3,1)$ && $-14$ &&
$\alpha_6x_1^4+\alpha_4x_1^3x_2+\alpha^1_2x_1^3x_3+
\alpha^2_2x_1^2x_2^2+
\!\!\alpha^1_0x_1^2x_2x_3+\alpha^2_0x_1x_2^3=\nomathbreak0$ &&
$x^2+y^4+y^3z+y^3+y^2z^2+y^2z+yz^3=\nomathbreak 0$ &&
$\wwt_x=\nomathbreak5$ $\wwt_y=\nomathbreak4$ $\wwt_z=\nomathbreak2$ &&
\centerline{$\Bbb{D}_6$}
\cr
\noalign{\hrule}
$H_{28}$ && $(5,4,0)$ && $-14$ &&
$\alpha_6x_1^4+\alpha_5x_1^3x_2+\alpha_1x_1^3x_3+\alpha_4x_1^2x_2^2+
\alpha_0x_1^2x_2x_3+\alpha_3x_1x_2^3+\alpha_2x_2^4=\nomathbreak0$ &&
$x^2+y^4+y^3z+y^3+y^2z^2+y^2z+yz^3+z^4=\nomathbreak0$ &&
$\wwt_x=\nomathbreak4$ $\wwt_y=\nomathbreak3$ $\wwt_z=\nomathbreak2$ &&
\centerline{$\Bbb{D}_5$}
\cr
\noalign{\hrule}
$H_{29}$ && $(5,4,1)$ && $-16$ &&
$\alpha_4x_1^4+\alpha_3x_1^3x_2+\alpha^1_0x_1^3x_3+\alpha_2x_1^2x_2^2+
\alpha_1x_1x_2^3+\alpha^2_0x_2^4=\nomathbreak 0$ &&
$x^2+y^4+y^3z+y^3+y^2z^2+yz^3+z^4=\nomathbreak 0$ &&
$\wwt_x=\nomathbreak6$ $\wwt_y=\nomathbreak4$ $\wwt_z=\nomathbreak3$ &&
\centerline{$\Bbb{E}_6$}
\cr
\noalign{\hrule}
$H_{30}$ && $(6,2,0)$ && $-12$ &&
$\alpha_{12}x_1^4+\alpha_8x_1^3x_2+\alpha_6x_1^3x_3+
\alpha_4x_1^2x_2^2+\alpha_2x_1^2x_2x_3+
\alpha_0^1x_1^2x_3^2+\alpha_0^2x_1x_2^3=\nomathbreak0$ &&
$x^2+y^4+y^3z+y^3+y^2z^2+y^2z+y^2+yz^3=\nomathbreak0$ &&
$\wwt_x=\nomathbreak3$ $\wwt_y=\nomathbreak3$ $\wwt_z=\nomathbreak1$ &&
\centerline{$\Bbb{A}_5$}
\cr
\noalign{\hrule}
$H_{31}$ && $(6,3,0)$ && $-14$ &&
$\alpha_{10}x_1^4+\alpha_7x_1^3x_2+\alpha_4x_1^3x_3+
\alpha_2x_1^2x_2^2+
\!\!\alpha_0x_1^2x_2x_3+\alpha_1x_1x_2^3=\nomathbreak0$ &&
$x^2+y^4+y^3z+y^3+y^2z^2+y^2z+yz^3=\nomathbreak 0$ &&
$\wwt_x=\nomathbreak5$ $\wwt_y=\nomathbreak4$ $\wwt_z=\nomathbreak2$ &&
\centerline{$\Bbb{D}_6$}
\cr
\noalign{\hrule}
$H_{32}$ && $(6,4,0)$ && $-16$ &&
$\alpha_{8}x_1^4+\alpha_6x_1^3x_2+\alpha_2x_1^3x_3+
\alpha_4x_1^2x_2^2+
\!\!\alpha_0^1x_1^2x_2x_3+\alpha_0^2x_2^4=\nomathbreak0$ &&
$x^2+y^4+y^3z+y^3+y^2z^2+y^2z+z^4=\nomathbreak 0$ &&
$\wwt_x=\nomathbreak4$ $\wwt_y=\nomathbreak3$ $\wwt_z=\nomathbreak2$ &&
\centerline{$\Bbb{D}_5$}
\cr
\noalign{\hrule}
$H_{33}$ && $(6,4,1)$ && $-18$ &&
$\alpha_{6}x_1^4+\alpha_4x_1^3x_2+\alpha_0^1x_1^3x_3+
\alpha_2x_1^2x_2^2+\alpha_0^2x_1x_2^3=\nomathbreak0$ &&
$x^2+y^4+y^3z+y^3+y^2z^2+yz^3=\nomathbreak0$ &&
$\wwt_x=\nomathbreak9$ $\wwt_y=\nomathbreak6$ $\wwt_z=\nomathbreak4$ &&
\centerline{$\Bbb{E}_7$}
\cr
\noalign{\hrule}
$H_{34}$ && $(6,5,0)$ && $-18$ &&
$\alpha_{6}x_1^4+\alpha_5x_1^3x_2+\alpha_0x_1^3x_3+
\alpha_4x_1^2x_2^2+\alpha_3x_1x_2^3+\alpha_2x_2^4=\nomathbreak 0$
&&
$x^2+y^4+y^3z+y^3+y^2z^2+yz^3+z^4=\nomathbreak 0 $ &&
$\wwt_x=\nomathbreak6$ $\wwt_y=\nomathbreak4$ $\wwt_z=\nomathbreak3$ &&
\centerline{$\Bbb{E}_6$}
\cr
\noalign{\hrule}
$H_{35}$ && $(7,3,0)$ && $-16$ &&
$\alpha_{12}x_1^4+\alpha_8x_1^3x_2+\alpha_5x_1^3x_3+\alpha_4x_1^2x_2^2+
\!\!\alpha_1x_1^2x_2x_3+\alpha_0x_1x_2^3=\nomathbreak0$ &&
$x^2+y^4+y^3z+y^3+y^2z^2+y^2z+yz^3=\nomathbreak 0$ &&
$\wwt_x=\nomathbreak5$ $\wwt_y=\nomathbreak4$ $\wwt_z=\nomathbreak2$ &&
\centerline{$\Bbb{D}_6$}
\cr
\noalign{\hrule}
}}
$$
$$
\vbox{
\let\PLUS=+
\binoppenalty10000
\relpenalty10000
\catcode`\+=\active
\raggedright
\def+{\allowbreak\PLUS\nobreak\unskip}
\define\mystrut{\vphantom{$\Bigl($}}
\def\ensp{\kern0.40em}
\tabskip=0pt \offinterlineskip
\halign{\vrule\mystrut
\ensp\hfill#\hfill\ensp&\vrule# &
\ensp\hfill#\hfill\ensp&\vrule# &
\ensp\hfill#\hfill\ensp&\vrule# &
\ensp\hfill\tabboxone{#}\hfill\ensp&\vrule# &
\ensp\hfill\tabboxtwo{#}\hfill\ensp&\vrule# &
\ensp\hfill\tabboxthree{#}\hfill\ensp&\vrule# &
\ensp\hfill\tabboxthree{\hsize 1.2cm #}\hfill\ensp\vrule
\cr
\multispan{13}
\hfill{\smc Table 1 continued}
\cr
\noalign{\medskip}
\noalign{\hrule}
$H_{36}$ && $(7,4,0)$ && $-18$ &&
$\alpha_{10}x_1^4+\alpha_7x_1^3x_2+\alpha_3x_1^3x_3+\alpha_4x_1^2x_2^2+
\!\!\alpha_0x_1^2x_2x_3+\alpha_1x_1x_2^3=\nomathbreak0$ &&
$x^2+y^4+y^3z+y^3+y^2z^2+y^2z+yz^3=\nomathbreak 0$ &&
$\wwt_x=\nomathbreak5$ $\wwt_y=\nomathbreak4$ $\wwt_z=\nomathbreak2$ &&
\centerline{$\Bbb{D}_6$}
\cr
\noalign{\hrule}
$H_{37}$ && $(7,5,0)$ && $-20$ &&
$\alpha_{8}x_1^4+\alpha_6x_1^3x_2+\alpha_1x_1^3x_3+\alpha_4x_1^2x_2^2+
\alpha_2x_1x_2^3+\alpha_{0}x_{2}^{4}=\nomathbreak 0$ &&
$x^2+y^4+y^3z+y^3+y^2z^2+yz^3+z^{4}=\nomathbreak 0$&&
$\wwt_x=\nomathbreak6$ $\wwt_y=\nomathbreak4$ $\wwt_z=\nomathbreak3$ &&
\centerline{$\Bbb{E}_6$}
\cr
\noalign{\hrule}
$H_{38}$ && $(7,5,1)$ && $-22$ &&
$\alpha_{6}x_1^4+\alpha_4x_1^3x_2+\alpha_0^1x_1^3x_3+
\alpha_2x_1^2x_2^2+\alpha_0^2x_1x_2^3=\nomathbreak0$ &&
$x^2+y^4+y^3z+y^3+y^2z^2+yz^3=\nomathbreak0$ &&
$\wwt_x=\nomathbreak9$ $\wwt_y=\nomathbreak6$ $\wwt_z=\nomathbreak4$ &&
\centerline{$\Bbb{E}_7$}
\cr
\noalign{\hrule}
$H_{39}$ && $(8,4,0)$ && $-20$ &&
$\alpha_{12}x_1^4+\alpha_8x_1^3x_2+\alpha_4^1x_1^3x_3+\alpha_4^2x_1^2x_2^2+
\!\!\alpha_0^1x_1^2x_2x_3+\alpha_0^2x_1x_2^3=\nomathbreak0$ &&
$x^2+y^4+y^3z+y^3+y^2z^2+y^2z+yz^3=\nomathbreak 0$ &&
$\wwt_x=\nomathbreak5$ $\wwt_y=\nomathbreak4$ $\wwt_z=\nomathbreak2$ &&
\centerline{$\Bbb{D}_6$}
\cr
\noalign{\hrule}
$H_{40}$ && $(8,5,0)$ && $-22$ &&
$\alpha_{10}x_1^4+\alpha_7x_1^3x_2+\alpha_2x_1^3x_3+
\alpha_4x_1^2x_2^2+\alpha_1x_1x_2^3=\nomathbreak0$ &&
$x^2+y^4+y^3z+y^3+y^2z^2+yz^3=\nomathbreak0$ &&
$\wwt_x=\nomathbreak9$ $\wwt_y=\nomathbreak6$ $\wwt_z=\nomathbreak4$ &&
\centerline{$\Bbb{E}_7$}
\cr
\noalign{\hrule}
$H_{41}$ && $(8,6,0)$ && $-24$ &&
$\alpha_{8}x_1^4+\alpha_6x_1^3x_2+\alpha_0^1x_1^3x_3+
\alpha_4x_1^2x_2^2+\alpha_2x_1x_2^3+\alpha_0^2x_2^4=\nomathbreak0$ &&
$x^2+y^4+y^3z+y^3+y^2z^2+yz^3+z^4=\nomathbreak 0$ &&
$\wwt_x=\nomathbreak6$ $\wwt_y=\nomathbreak4$ $\wwt_z=\nomathbreak3$ &&
\centerline{$\Bbb{E}_6$}
\cr
\noalign{\hrule}
$H_{42}$ && $(9,5,0)$ && $-24$ &&
$\alpha_{12}x_1^4+\alpha_8x_1^3x_2+\alpha_3x_1^3x_3+
\alpha_4x_1^2x_2^2+\alpha_0x_1x_2^3=\nomathbreak0$ &&
$x^2+y^4+y^3z+y^3+y^2z^2+yz^3=\nomathbreak0$ &&
$\wwt_x=\nomathbreak9$ $\wwt_y=\nomathbreak6$ $\wwt_z=\nomathbreak4$ &&
\centerline{$\Bbb{E}_7$}
\cr
\noalign{\hrule}
$H_{43}$ && $(9,6,0)$ && $-26$ &&
$\alpha_{10}x_1^4+\alpha_7x_1^3x_2+\alpha_1^1x_1^3x_3+
\alpha_4x_1^2x_2^2+\alpha_1^2x_1x_2^3=\nomathbreak0$ &&
$x^2+y^4+y^3z+y^3+y^2z^2+yz^3=\nomathbreak0$ &&
$\wwt_x=\nomathbreak9$ $\wwt_y=\nomathbreak6$ $\wwt_z=\nomathbreak4$ &&
\centerline{$\Bbb{E}_7$}
\cr
\noalign{\hrule}
$H_{44}$ && $(10,6,0)$ && $-28$ &&
$\alpha_{12}x_1^4+\alpha_8x_1^3x_2+\alpha_2x_1^3x_3+
\alpha_4x_1^2x_2^2+\alpha_0x_1x_2^3=\nomathbreak0$ &&
$x^2+y^4+y^3z+y^3+y^2z^2+yz^3=\nomathbreak0$ &&
$\wwt_x=\nomathbreak9$ $\wwt_y=\nomathbreak6$ $\wwt_z=\nomathbreak4$ &&
\centerline{$\Bbb{E}_7$}
\cr
\noalign{\hrule}
$H_{45}$ && $(10,7,0)$ && $-30$ &&
$\alpha_{10}x_1^4+\alpha_7x_1^3x_2+\alpha_0x_1^3x_3+
\alpha_4x_1^2x_2^2+\alpha_1x_1x_2^3=\nomathbreak0$ &&
$x^2+y^4+y^3z+y^3+y^2z^2+yz^3=\nomathbreak0$ &&
$\wwt_x=\nomathbreak9$ $\wwt_y=\nomathbreak6$ $\wwt_z=\nomathbreak4$ &&
\centerline{$\Bbb{E}_7$}
\cr
\noalign{\hrule}
}}
$$
$$
\vbox{
\let\PLUS=+
\binoppenalty10000
\relpenalty10000
\catcode`\+=\active
\raggedright
\def+{\allowbreak\PLUS\nobreak\unskip}
\define\mystrut{\vphantom{$\Bigl($}}
\def\ensp{\kern0.40em}
\tabskip=0pt \offinterlineskip
\halign{\vrule\mystrut
\ensp\hfill#\hfill\ensp&\vrule# &
\ensp\hfill#\hfill\ensp&\vrule# &
\ensp\hfill#\hfill\ensp&\vrule# &
\ensp\hfill\tabboxone{#}\hfill\ensp&\vrule# &
\ensp\hfill\tabboxtwo{#}\hfill\ensp&\vrule# &
\ensp\hfill\tabboxthree{#}\hfill\ensp&\vrule# &
\ensp\hfill\tabboxthree{\hsize 1.2cm #}\hfill\ensp\vrule
\cr
\multispan{13}
\hfill{\smc Table 1 continued}
\cr
\noalign{\medskip}
\noalign{\hrule}
$H_{46}$ && $(11,7,0)$ && $-32$ &&
$\alpha_{12}x_1^4+\alpha_8x_1^3x_2+\alpha_1x_1^3x_3+
\alpha_4x_1^2x_2^2+\alpha_0x_1x_2^3=\nomathbreak0$ &&
$x^2+y^4+y^3z+y^3+y^2z^2+yz^3=\nomathbreak0$ &&
$\wwt_x=\nomathbreak9$ $\wwt_y=\nomathbreak6$ $\wwt_z=\nomathbreak4$ &&
\centerline{$\Bbb{E}_7$}
\cr
\noalign{\hrule}
$H_{47}$ && $(12,8,0)$ && $-36$ &&
$\alpha_{12}x_1^4+\alpha_8x_1^3x_2+\alpha_0x_1^3x_3+
\alpha_4x_1^2x_2^2+\alpha_0x_1x_2^3=\nomathbreak0$ &&
$x^2+y^4+y^3z+y^3+y^2z^2+yz^3=\nomathbreak0$ &&
$\wwt_x=\nomathbreak9$ $\wwt_y=\nomathbreak6$ $\wwt_z=\nomathbreak4$ &&
\centerline{$\Bbb{E}_7$}
\cr
\noalign{\hrule}
}}
$$

}

Thus Theorem~1.5 is proved.

\remark{Remark 3.11}
The proof of Theorem~1.5 can also be used to describe the singularities of
all 3-folds $H_{i}$. For example, sufficiently general 3-folds $H_{i}$ are
smooth for $i\in\{1,2,3,4,6,9\}$. A sufficiently general 3-fold $H_{5}$
has a single isolated ordinary double point and is not
$\Bbb{Q}$\kkk-factorial. The singularities of~$H_{i}$ are non-isolated for
$i\in\{8,10,12,14,16,17,18,20,22,27,29,33,38\}$. In all other cases, a
sufficiently general 3-fold $H_{i}$ has a single isolated
non-$\operatorname{cDV}$ singular point.
\endremark

\remark{Remark 3.12}
One can simplify the proof of Theorem~1.5 arguing as follows. If $X$ is
a del~Pezzo surface of degree~$2$ over some field with a non-Du~Val
singular point defined over this field, then the ramification divisor of
the double covering $X\to\Bbb{P}^2$ is a union of four lines. The latter
condition can easily be checked in terms of the numbers~$d_i$. However,
the authors did not use this approach, having in mind the future applications
to the rationality questions for~$H_i$: it is sometimes useful to know
the type of singularity or even the explicit local equations of~$X$
(see~\S\,5).
\endremark

\head
\S\,4. Proof of Theorem~1.6
\endhead

Let $X$ be a Fano 3-fold with canonical Gorenstein singularities such that
the linear system $|-K_{X}|$ has no base points and the induced morphism
$\varphi_{|-K_{X}|}$ is an embedding, but the anticanonical image
$\varphi_{|-K_{X}|}(X)\subset\Bbb{P}^{n}$ is not an intersection of quadrics.
$\bigl($Here $n=-{\frac{K_{X}^{3}}{2}}+2.\bigr)$

\remark{Remark 4.1}
If $-K_{X}^{3}=4$, then the 3-fold $X$ is a quartic (possibly singular)
in~$\Bbb{P}^{4}$. The birational geometry of such 3-folds was studied in
\cite{22}, \cite{54}, \cite{37}, \cite{16},
\cite{28}, \cite{63}, \cite{64}, \cite{102},~\cite{58}.
\endremark

\remark{Remark 4.2}
If $-K_{X}^{3}=6$, then the 3-fold $X$ is a complete intersection (possibly
singular) of a quadric and a cubic in~$\Bbb{P}^{5}$. This easily follows
from either Theorem~2.14 or Proposition~2.21. The birational geometry of
such 3-folds was studied in~\cite{54},
\cite{37}, \cite{16}, \cite{29}, \cite{89},~\cite{63}.
\endremark

\medskip
Thus we may assume that $-K_{X}^{3}\geqslant 8$.
Hence Theorem~2.14 implies that $X$ is projectively normal in~$\Bbb{P}^{n}$
and the quadrics through $X$ in~$\Bbb{P}^{n}$ cut out a 4-fold
$Y\subset\Bbb{P}^{n}$ of degree $n-3$. Moreover, if
$\Pi\subset\Bbb{P}^{n}$ is a general linear subspace of codimension
$2$, then the curve $X\cap\Pi$ is either a canonically embedded smooth
trigonal curve or a canonically embedded smooth plane quintic curve, and
$\operatorname{deg}(X\subset\Bbb{P}^{n})=10$. In the former case, the 4-fold
$Y$ is the image of a rational scroll
$\operatorname{Proj}\bigl(\,
\bopl_{i=1}^{4}\Cal{O}_{\Bbb{P}^{1}}(d_{i})\bigr)$
under the map given by the tautological line bundle over~$\PP^1$,
where $0\ne d_{1}\geqslant\dots\geqslant d_{4}\geqslant 0$ and
$-K_{X}^{3}=2(d_{1}+d_{2}+d_{3}+d_{4})+2$. In the latter case we have
$n=7$ and $Y$ is a cone over the Veronese surface $v_2(\PP^2)$.

\remark{Remark 4.3}
The cone over a Veronese surface in~$\Bbb{P}^{7}$ is isomorphic to
$\Bbb{P}(1^{3},2,2)$. Therefore, if~$Y$ is a cone over a Veronese
surface, then the 3-fold $X$ is a hypersurface in
$\Bbb{P}(1^{3},2^{2})$ of degree~$5$ because $-K_{X}^{3}=10$.
\endremark

\proclaim{Lemma 4.4}
Let $Y$ be a cone over a Veronese surface in~$\Bbb{P}^{7}$ whose vertex is
a line $L\subset\Bbb{P}^{7}$. Take a resolution of singularities
$f\:U\to Y$, where
$U=\operatorname{Proj}(\Cal{O}_{\Bbb{P}^{2}}(2)\oplus
\Cal{O}_{\Bbb{P}^{2}}\oplus\Cal{O}_{\Bbb{P}^{2}})$.
Put $T=f^{*}(\Cal{O}_{\Bbb{P}^{7}}(1)|_{Y})$ and let $F$
be the pull back of~$\Cal{O}_{\Bbb{P}^{2}}(1)$ under the
natural projection of~$U$ to~$\Bbb{P}^{2}$. Put
$V=f^{-1}(X)\subset U$. Then $V$ has canonical Gorenstein singularities,
$-K_{V}$ is big and numerically effective, and we have $V\sim 2T+F$
on~$U$ and $X=\varphi_{|-rK_{V}|}(V)$ for $r\gg 0$. In particular, the
birational morphism $f|_{V}\:V\to X$ is crepant.
\endproclaim

\demo{Proof}
The line $L$ is contained in~$X$ by Remark~4.3. On the other hand,
the 3-fold $X$ is singular along $L\subset X$. Indeed, let
$O\in L$ be a point, $\Pi\subset\Bbb{P}^{7}$ a sufficiently general linear
subspace of codimension~$2$ through~$O$, and~$C=\Pi\cap X$.
Suppose that $O$ is smooth on~$X$. Then $C$ is a smooth anticanonically
embedded plane quintic curve. Hence quadrics through $C$ in
$\Pi\cong\Bbb{P}^{5}$ cut out a smooth Veronese surface by Theorem~2.13.
On the other hand, quadrics through $C$ in~$\Pi\cong\Bbb{P}^{5}$ cut out
$Y\cap\Pi$. However, the surface $Y\cap\Pi$ must be singular because
$Y$ is singular at~$O$ by the hypothesis.

Therefore the morphism $f|_{V}\:V\to X$ is crepant at the general point of
the line $L\subset X$ by Proposition~2.15. It follows that
$V$ contains no fibres of~$f$ and $V\sim 2T+F$ on~$U$. Hence the 3-fold
$V$ is normal (see~\cite{83}, Proposition~8.23).
Therefore $V$ has canonical Gorenstein singularities and
$-K_{V}$ is a crepant pull back of~$-K_{X}$.
\enddemo

\remark{Remark 4.5}
Lemma~4.4 does not {\it a priori\/} imply the existence of the corresponding
Fano 3-fold $X\subset Y$. However, this existence is easily seen. In the
notation of Lemma~4.4, the linear system $|2T+L|$ on the 4-fold
$U=\operatorname{Proj}(\Cal{O}_{\Bbb{P}^{2}}(2)\oplus
\Cal{O}_{\Bbb{P}^{2}}\oplus\Cal{O}_{\Bbb{P}^{2}})$ is free.
In particular, sufficiently general divisors in this system are smooth.
Let $D$ be a divisor in~$|2T+L|$ with canonical singularities. Then the
adjunction formula implies that $-K_{D}\sim T$. Therefore $D$ is a weak
Fano 3-fold, that is, the divisor $-K_{D}$ is numerically effective and
big. The vanishing theorem (see~\cite{91}, \cite{120}) implies that
$\varphi_{|-K_{D}|}=\varphi_{|T|}|_{D}$. In particular, the 3-fold
$\varphi_{|-K_{D}|}(D)$ is a Fano 3-fold with canonical Gorenstein
singularities.
\endremark

\medskip
In what follows, we may thus assume that
$X\cap\Pi$ is a canonically embedded smooth trigonal curve for any
general linear subspace $\Pi\subset\Bbb{P}^{n}$ of codimension~$2$.
Therefore quadrics through $X$ in~$\Bbb{P}^{n}$ cut out a 4-fold $Y$
which is the image of a rational scroll
$\operatorname{Proj}\bigl(\,
\bopl_{i=1}^{4}\Cal{O}_{\Bbb{P}^{1}}(d_{i})\bigr)$
under the map given by the tautological line bundle, where
$0\ne d_{1}\geqslant\dots\geqslant d_{4}\geqslant 0$ and
$-K_{X}^{3}=2(d_{1}+d_{2}+d_{3}+d_{4})+2$.

\proclaim{Lemma 4.6}
The inclusion $\operatorname{Sing}(Y)\cap X\subset\operatorname{Sing}(X)$
holds.
\endproclaim

\demo{Proof}
Let $O$ be a singular point on~$Y$ such that $O\in X$ and the 3-fold
$X$ is non-singular at~$O$. Take a sufficiently general linear subspace
$\Pi\subset\Bbb{P}^{n}$ of codimension~$2$ passing through~$O$. Put
$C=\Pi\cap X$. Then the curve $C\subset\Pi\cong\Bbb{P}^{n-2}$ is a smooth
anticanonically embedded trigonal curve. Therefore quadrics through $C$
in~$\Pi\cong\Bbb{P}^{5}$ cut out a smooth surface by Theorem~2.13. On
the other hand, quadrics through $C$ in~$\Pi\cong\Bbb{P}^{n-2}$ cut out the
surface $Y\cap\Pi$, which must be singular at~$O$ because the 4-fold~$Y$
is singular at~$O$ by the assumption.
\enddemo

Let $f\:{\kern-.8pt}U{\kern-.8pt}\!\to\! Y$ be the birational morphism 
$f{\kern-.6pt}\!=\!\varphi_{|{\kern-.4pt}M{\kern-.4pt}|}$, where
$U{\kern-.6pt}\!=\!\operatorname{Proj}\bigl(
\bopl_{i={\kern-.4pt}1}^{4}\!\Cal{O}_{\Bbb{P}^{1}}{\kern-.4pt}(d_{i})\bigr)$, 
$0\ne d_{1}\geqslant\dots\geqslant d_{4}\geqslant 0$, \
$-K_{X}^{3}=2(d_{1}+d_{2}+d_{3}+d_{4})+2\geqslant 8$ and
$M$ is the tautological line bundle on~$U$. We put
$V=f^{-1}(X)\subset U$ and $h=f|_{V}\:V\to X$.

\proclaim{Lemma 4.7}
The 3-fold $V$ has canonical Gorenstein singularities, the anticanonical
divisor $-K_{V}$ is numerically effective and big, and
$K_{V}\sim h^{*}(K_{X})$, that is, the morphism $h$ is crepant.
\endproclaim

\demo{Proof}
The 3-fold $V$ is normal if and only if it is smooth in codimension~$1$
(see~\cite{83}, Proposition~8.23). On the other hand, if~$d_{3}\ne 0$, then
the 3-folds $X$ and $V$ are isomorphic in codimension~$2$. This immediately
yields the claim (compare~\cite{92}).

We may thus assume that $d_{3}=d_{4}=0$. Put
$Z=\operatorname{Sing}(Y)$. Then $\operatorname{dim}(Z)\leqslant 2$,
and the equation $\operatorname{dim}(Z)=2$ holds if and only if
$d_{2}=0$. Moreover, if~$\operatorname{dim}(Z\cap X)=0$, then the 3-folds
$X$ and $V$ are isomorphic in codimension~$2$, which implies the claim.
On the other hand, we have $Z\cap X\subset\operatorname{Sing}(X)$ by
Lemma~4.6. Since $X$ is normal, it follows that
$\operatorname{dim}(Z\cap X)\leqslant
\operatorname{dim}(\operatorname{Sing}(X))\leqslant 1$. We may thus assume
that the intersection $Z\cap X$ consists of finitely many curves and
$X$ is singular along every curve in~$Z\cap X$.

The canonicity of singularities of~$X$ and Proposition~2.15 imply that
$g$ is crepant at the general point of every curve in~$Z\cap X$ and the
singularities of~$V$ are canonical Gorenstein over the general point of
every curve in~$Z\cap X$. This proves the assertion of the lemma for
the complement to a subset of codimension 2 in~$V$. It follows that
$V$ is normal (see~\cite{83}, Proposition~8.23).
Hence $V\subset U$ is a divisor on a smooth 4-fold, $V$ is a normal 3-fold,
and $V$ has canonical Gorenstein singularities in codimension~2. It follows
that singularities of~$V$ are canonical Gorenstein, and
$K_{V}\sim h^{*}(K_{X})$.
\enddemo

Let $M$ be the tautological line bundle on~$U$, and let $L$ be a fibre of
the natural projection of~$U$ to~$\Bbb{P}^{1}$. Then $-K_{V}\sim M|_{V}$
by the construction.

\remark{Remark 4.8}
The contraction theorem (see~\cite{93}) does not {\it a priori\/} imply that
$X$ is an anticanonical model of~$V$. However, the vanishing theorem
(see~\cite{91}, \cite{120}) implies that $|-K_{V}|=|M|_{V}|$.
Therefore we see {\it a posteriori\/} that $X$ is an anticanonical image
of~$V$, that is, $X=\varphi_{|-K_{V}|}(V)$.
\endremark

\medskip
The adjunction formula implies that $V\sim 3M-(d_{1}+d_{2}+d_{3}+d_{4}-2)L$
on~$U$. Let $Y_{j}\subset\nomathbreak V$ be the subscroll induced by the
natural projection
$$
\bopl_{i=1}^{4}\Cal{O}_{\Bbb{P}^1}(d_i)\to
\bopl_{i=j}^{4}\Cal{O}_{\Bbb{P}^1}(d_i).
$$
In particular, $Y_{4}$ is a curve, $Y_{3}$ is a surface, and
$Y_{2}$ is a 3-fold.

\proclaim{Lemma 4.9} The following inequalities hold:
$$
\operatorname{mult}_{Y_{2}}(V)=0,\qquad
\operatorname{mult}_{Y_{3}}(V)\leqslant 1,\qquad
\operatorname{mult}_{Y_{4}}(V)\leqslant 2.
$$
\endproclaim

\demo{Proof}
The first inequality is
obvious, the second one follows since $V$ is normal,
and the last one follows from the canonicity of~$V$ at a general point of
$Y_{4}$ by Theorem~2.9.
\enddemo

Therefore Corollary 2.20 yield the inequalities
$$
2d_{2}-d_{1}-d_{3}-d_{4}+2\geqslant 0,\qquad
d_3-d_2-d_4+2\geqslant 0,\qquad
2-d_2-d_3+d_1\geqslant 0.
$$
We also have
$0\ne d_{1}\geqslant\dots\geqslant d_{4}\geqslant 0$
and $d_{1}+d_{2}+d_{3}+d_{4}\geqslant 3$ by the assumption.
These inequalities determine exactly $66$ different rational scrolls
$$
\Bbb{F}(d_{1},d_{2},d_{3},d_{4})=\operatorname{Proj}\biggl(\,
\bopl_{i=1}^{4}\Cal{O}_{\Bbb{P}^{1}}(d_{i})\biggr).
$$

\remark{Remark 4.10}
The 3-fold $X$ is an anticanonical model of the 3-fold $V$, that is,
$X=\varphi_{|-rK_{V}|}(V)$ for $r\gg 0$. Thus the contraction theorem
(see~\cite{93}) implies that $X$ is a Fano 3-fold with canonical Gorenstein
singularities if and only if~$V$ has canonical Gorenstein singularities and
$-K_{V}$ is numerically effective and big. On the other hand, the 3-fold
$V$ is uniquely determined by the rational scroll
$\Bbb{F}(d_{1},d_{2},d_{3},d_{4})$ and the class
$3M-(d_{1}+d_{2}+d_{3}+d_{4}-2)L$ in the Picard group of
$\Bbb{F}(d_{1},d_{2},d_{3},d_{4})$. However, the linear system
$|3M-(d_{1}+d_{2}+d_{3}+d_{4}-2)L|$ may {\it a priori\/} contain no divisor
with canonical singularities.
\endremark

\medskip
In the rest of the section we explicitly show that, for each of the
$66$ possible cases, the linear system $|3M-(d_{1}+d_{2}+d_{3}+d_{4}-2)L|$
contains a divisor with canonical singularities. To do this, we use
Corollary~2.5 and Proposition~2.19. This will complete the proof of
Theorem~1.6.

\remark{Remark 4.11}
The same idea was used in the classification of smooth trigonal Fano 3-folds
(see~\cite{15}, \cite{17}) and in~\cite{41} to prove the effective
boundedness of degree for the trigonal Fano 3-folds with canonical
singularities. The maximal value of the degree is attained by
Theorem~1.6.
\endremark

\medskip
We start by considering two different possible cases in full detail.

\example{Example 4.12}
Let $V\subset\FF=\FF(6,5,3,0)$ be a sufficiently general divisor in the
linear system $|3M-12L|$. Let us show that $V$ has canonical singularities.

\goodbreak
By Proposition~2.19, $V$ is given by the zeros of a bihomogeneous polynomial
$$
\align
&\alpha_6(t_1,t_2)\,x_1^3+\alpha_5(t_1,t_2)\,x_1^2x_2
+\alpha^1_3(t_1,t_2)\,x_1^2x_3+\alpha^1_0(t_1,t_2)\,x_1^2x_4
+\alpha_4(t_1,t_2)\,x_1x_2^2
\\
&\qquad\qquad
+\alpha_2(t_1,t_2)\,x_1x_2x_3+\alpha^2_0(t_1,t_2)\,x_1x_3^2
+\alpha^2_3(t_1,t_2)\,x_2^3+\alpha_1(t_1,t_2)\,x_2^2x_3,
\endalign
$$
where $\alpha_d(t_1,t_2)$ (or $\alpha^i_d(t_1,t_2)$) is a form of
degree~$d$. Let $E$ be the surface $x_1=x_2=0$, and let $C$ be the curve
$x_1=x_2=x_3=0$. We note that the base locus of~$|3M-12L|$ is equal to~$E$.
Since the automorphism group of~$E\cong\FF(3,0)$ acts transitively on
$E\setminus C$ and $V$ has multiplicity~$1$ at a general point of~$E$,
we see from Corollary~2.5 that it suffices to prove that for any point
$P$ on~$C$ there is a divisor $V$ with canonical singularities in a
neighbourhood of~$P$ (compare Example~3.10).

Let $Y$ be the fibre of~$V\subset\FF$ over a sufficiently general point
$P\in\PP^1$. We put $O=C\cap Y$. As above, it suffices to prove that
$V$ has at most canonical singularities in a neighbourhood of~$O$. Hence
it suffices to prove that $O$ is a Du~Val point on~$Y$.

Let $(\gamma:\delta)$ be the homogeneous coordinates of the point
$P\in\PP^1$. Then $Y$ may be presented as the hypersurface
$$
\align
&\alpha_6x_1^3+\alpha_5x_1^2x_2+\alpha^1_3x_1^2x_3+
\alpha^1_0x_1^2x_4+\alpha_4x_1x_2^2
\\
&\qquad\qquad
+\alpha_2x_1x_2x_3+\alpha^2_0x_1x_3^2+\alpha^2_3x_2^3+\alpha_1x_2^2x_3=0
\endalign
$$
in $\PP^4$, where $\alpha_i=\alpha_i(\gamma,\delta)$
(resp\.~$\alpha_i^j=\alpha_i^j(\gamma,\delta)$) and $O=(0:0:0:1)$.
Since $P$ is general, we have $\alpha_i\ne 0$
(resp\.~$\alpha_i^j\ne 0$) for all~$i$.
Hence we may assume for convenience that $\alpha_i=1$
(resp\.~$\alpha_i^j=1$) for all~$i$ and~$j$.

Let $x_1=x$, \ $x_2=y$, \ $x_3=z$ and $x_4=1$. Then the local equation of
$Y$ in a neighbourhood of~$O$ is
$$
x^3+x^2y+x^2z+x^2+xy^2+xyz+xz^2+y^3+y^2z=0.
$$
We put $\wt(x)=4$, \ $\wt(y)=3$ and $\wt(z)=2$. Then
$\wt(x^2+xz^2+y^2z)=8$, \ $\wt(x^3)=12$, \ $\wt(x^2y)=11$, \
$\wt(xy^2)=10$, \ $\wt(y^3)=9$, \ $\wt(x^2z)=10$ and $\wt(xyz)=9$.
Moreover, the singularity given by the equation $x^2+xz^2+y^2z=0$
is isolated. Therefore the singularity of~$Y$ at~$O$ is locally isomorphic
to a Du~Val singularity of type~$\Bbb{D}_5$. In particular, $V$
has a singularity of type $\Bbb{D}_5\times\CC$ at a general
point of~$C$. Using the generality in the choice of~$V$, we may actually
assume that $P$ is an arbitrary point of the curve~$C$ (compare
Example~3.10). Thus, for any given point $P$ of~$C$, the linear system
$|3M-12L|$ contains a divisor with at most canonical singularities in
the neighbourhood of~$P$. Hence the singularities of~$V$ are canonical
by Corollary~2.5.
\endexample

\example{Example 4.13}
Let $V\subset\FF=\FF(7,3,1,0)$ be a general divisor in the linear
system $|3M-\nomathbreak 9L|$. Let us show that $V$ has canonical
singularities.

By Proposition~2.19, $V$ is given by the zeros of a bihomogeneous polynomial
$$
\align
&\alpha_{12}(t_1,t_2)\,x_1^3+\alpha_8(t_1,t_2)\,x_1^2x_2
+\alpha_6(t_1,t_2)\,x_1^2x_3+\alpha_5(t_1,t_2)\,x_1^2x_4
\\
&\qquad\qquad
+\alpha_4(t_1,t_2)\,x_1x_2^2+\alpha_2(t_1,t_2)\,x_1x_2x_3
+\alpha_1(t_1,t_2)\,x_1x_2x_4
\\
&\qquad\qquad
+\alpha^1_0(t_1,t_2)\,x_1x_3^2+\alpha^2_0(t_1,t_2)\,x_2^3,
\endalign
$$
\goodbreak
\noindent
where $\alpha_d(t_1,t_2)$ (or $\alpha^i_d(t_1,t_2)$) is a form of
degree~$d$. Let $E$ be the surface $x_1=x_2=0$, and let $C$ be the curve
$x_1=x_2=x_3=0$. We note that the base locus of~$|3M-9L|$ is equal to~$E$.
Since the automorphism group of~$E\cong\FF(1,0)$ acts transitively on
$E\setminus C$ and $V$ has multiplicity~$1$ at a general point of~$E$,
Corollary~2.5 implies that $V$ has canonical singularities on
$E\setminus C$, and it remains to verify that $V$ has canonical
singularities at points of~$C$ (compare Example~3.10).

Let $Y$ be the fibre of~$V\subset\FF$ over a general point $P\in\PP^1$.
We put $O=C\cap Y$. As above, it suffices to prove that $V$ has canonical
singularities at~$O$. Hence it suffices to prove that $O$ is a Du~Val
point on~$Y$.

Let $(\gamma:\delta)$ be the homogeneous coordinates of~$P\in\PP^1$.
Then $Y$ may be presented as a hypersurface corresponding to the
polynomial
$$
\align
&\alpha_{12}x_1^3+\alpha_8x_1^2x_2+\alpha_6x_1^2x_3+
\alpha_5x_1^2x_4+\alpha_4x_1x_2^2\\
&\qquad\qquad
+\alpha_2x_1x_2x_3+\alpha_1x_1x_2x_4+\alpha^1_0x_1x_3^2+\alpha^2_0x_2^3
\endalign
$$
on $\PP^4$, where $\alpha_i=\alpha_i(\gamma,\delta)$
(resp\.~$\alpha_i^j=\alpha_i^j(\gamma, \delta)$) and
$O=(0:0:0:1)$. Since $P$ is general, we have
$\alpha_i\ne 0$ (resp\.~$\alpha_i^j\ne 0$) for all~$i$.
Therefore we may assume for convenience that $\alpha_i=1$
(resp\.~$\alpha_i^j=1$) for all~$i$ and~$j$.

We put $x_1=x$, \ $x_2=y$, \ $x_3=z$ and $x_4=1$. Then the local equation
of $Y$ in a neighbourhood of~$O$ is
$$
x^3+x^2y+x^2z+x^2+xy^2+xyz+xy+xz^2+y^3=0.
$$
It is easy to prove that one cannot choose ``good'' weights for this
polynomial. Hence we cannot use Theorem~2.10 as before and we must
explicitly resolve singularities in a neighbourhood of~$O$.
It is easy to see that, resolving the singularity at~$O$, we can skip
monomials whose weight is larger than at least one of the others.
(This will be clear from the forthcoming blow ups.) Hence it suffices to
consider the polynomial
$$
x^2+xy+xz^2+y^3.
$$

This polynomial determines an isolated singular point at~$O$.
We blow up this point. The formulae in the local charts are as follows.

1) $x\neq 0$: the change of coordinates $x=x$, \ $y=xy$, \ $z=xz$
brings the local equation (after dividing by~$x^2$) to the form
$$
1+y+xz^2+xy^3=0,
$$
so our surface is smooth in this chart.

2) $y\neq 0$: the change of coordinates $x=xy$, \ $y=y$, \ $z=yz$
yields (after dividing by~$y^2$)
$$
x^2+x+xyz^2+y=0,
$$
and our surface is smooth in this chart.

3) $z\neq 0$: the change of coordinates $x=xz$, \ $y=yz$, \ $z=z$
yields (after dividing by~$z^2$)
$$
x^2+xy+xz+y^3z=0,
$$
and we have two extremal
$(-2)$-curves with the only singular point $(0,0,0)$
near $z=0$.

Therefore more blow ups are necessary. We must study the
singularities given by the local equation
$$
x^2+xy+xz+y^3z=0.
$$
We blow up the point $(0,0,0)$. Here are the formulae in the local charts.

1) $x\neq 0$: the change of coordinates $x=x$, \ $y=xy$, \ $z=xz$
yields (after dividing by~$x^2$)
$$
1+y+z+x^2y^3z=0,
$$
and the surface is smooth in this chart.

2) $y\neq 0$: the change of coordinates $x=xy$, \ $y=y$, \ $z=yz$
yields (after dividing by~$y^2$)
$$
x^2+x+xz+y^2z=x^2+x(z+1)+y^2(z+1)-y^2=x^2+y^{\prime\,2}+z^{\prime\,2}-y^{\prime\,2}z',
$$
where $y'=iy$ and $z'=z+1$. Thus the point $(0,0,-1)$ is Du~Val of
type~$\Bbb{A}_1$.

3) $z\neq 0$: the change of coordinates $x=xz$, \ $y=yz$, \ $z=z$ yields
(after dividing by~$z^2$)
$$
x^2+xy+x+y^3z^2=0.
$$
We can see that the point $(0,-1,0)$ is singular. It coincides with the
singular point in the chart $y\neq 0$.

Summarizing, we have two $(-2)$-curves after the first blow up, two
$(-2)$-curves after the second blow up and one point of type
$\Bbb{A}_1$ on one of these curves. Hence the graph of resolution of the
original singularity corresponds to a Du~Val singularity of type
$\Bbb{A}_5$. In particular, the 3-fold $V$ has a singularity of type
$\Bbb{A}_5\times\CC$ at the general point of~$C$. As in Example~4.12,
it follows that the singularities of~$V$ are canonical by Corollary~2.5.
\endexample

We state the results of some easy calculations, which will be used in
the remaining part of the proof of Theorem~1.6.

\proclaim{Lemma 4.14} \
\rom{1.}
A surface singularity given by
$$
x^3+x^2y+x^2z+x^2+xy^2+xyz+xy+xz^2+xz+y^3=0
$$
is Du~Val of type
$\Bbb{A}_2$.

\rom{2.}
A surface singularity given by
$$
x^3+x^2y+x^2z+x^2+xy^2+xyz+xy+xz^2+xz+y^3+y^2z,
$$
is Du~Val of type $\Bbb{A}_3$.

\rom{3.}
A surface singularity given by
$$
x^3+x^2y+x^2z+x^2+xy^2+xyz+xy+xz^2+xz+y^3+y^2z+yz^2+z^3,
$$
is Du~Val of type $\Bbb{A}_3$.

\rom{4.}
A surface singularity given by
$$
x^3+x^2y+x^2z+x^2+xy^2+xyz+xy+xz^2+y^3+y^2z,
$$
is Du~Val of type $\Bbb{A}_4$.

\rom{5.}
A surface singularity given by
$$
x^3+x^2y+x^2z+x^2+xy^2+xyz+xy+xz^2+y^3,
$$
is Du~Val of type $\Bbb{A}_5$.
\endproclaim

\demo{Proof}
In case~1, changing the coordinates by~$x'=x$, \
$y'=y$, \ $z'=x+y+z$, we get the equation
$$
x'z'+y^{\prime\,3}+Q(x', y',z')=0,
$$
where $Q$ consists of the terms whose weight
(with respect to any choice of weights)
is larger than that of either $x'z'$ or $y^{\prime\,3}$.
Hence this singularity is Du~Val of 
\linebreak
type $\Bbb{A}_2$.

\looseness=1
In the remaining cases, it is easy to see that the corresponding equation
describes an isolated singularity. Since it is impossible to use Theorem~2.5,
we shall explicitly resolve these singularities. Case~5 has already been
discussed in Example~4.13, along with the details of calculations.
In cases 2 and~3, a single blow up yields two exceptional $(-2)$-curves
and one point of type $\Bbb{A}_1$ on one of them. This means that the
original singularities are Du~Val of type $\Bbb{A}_3$. In case~4, two blow
ups yield a smooth surface, and the exceptional curves form a Dynkin diagram
of type~$\Bbb{A}_4$.
\enddemo

In the rest of the section we consider all possible cases, following the
pattern of Examples~4.12 and~4.13 and using Lemma~4.14 when necessary.
The differences appear only in numerical characteristics, equations, types
of singularities etc. They are surveyed in Table~2.

\looseness=1
Table 2 is organized as follows. The first column contains the labels of
varieties $V$ in the notation of Theorem~1.6. The second column yields a
quadruple $(d_1, d_2, d_3, d_4)$ such that $V$ is a divisor on
$\FF(d_1, d_2, d_3, d_4)$. The third column yields the number~$b$
such that $|V|=|3M+bL|$. In cases $T_{5}$, $T_{8}$ and~$T_{14}$, this
linear system appears to be base point free, whence $V$ is smooth by
the Bertini theorem. Then we do not need the information contained in
the other columns. In the other cases,
$\operatorname{Bs}|V|$ is either the curve $C=Y_4$ given by
$x_1=x_2=x_3=0$ (then $V$ is non-singular outside~$C$, so it suffices to
verify that $V$ has only canonical singularities at points of~$C$) or
the surface $E=Y_3$ given by~$x_1=x_2=0$ (then
$V$ has multiplicity~$1$ at a general point of~$E$ and, since the
automorphism group of~$E\cong\FF(d_3, d_4)$ acts transitively on
$E\setminus C$ (compare Example~3.10), we see from Corollary~2.5 that it
suffices to study the singularities of~$V$ at points of~$C$ only).

\looseness=1
The fourth column yields an equation of a general divisor~$V$ in the linear
system $|3M+bL|$. The fifth column yields the equation of the fibre $Y$ of
the projection $V\to\PP^1$ over a general point of~$\PP^1$ in the
neighbourhood of a general point of~$C$ after the change of coordinates
$x_1=x$, \ $x_2=y$, \ $x_3=z$, \ $x_4=1$ (see Example~4.12). The same
equation locally describes $V$ as an equation in~$t$,~$x$,~$y$,~$z$. For
the 3-folds $T_j$ with
$j\in\{4,\,6,\,7,\,9,\,11,\,15,\,16,\,17,\,19,\,25,\,26,\,28,\,36,\,45\}$, the
corresponding point appears to be non-singular. For
$T_j$ with $j\in\{10, 12, 13, 18, 
\mathbreak
22, \,23, \,24, \,30, \,31, \,33, \,34, \,35, \,41, \,43, \,44, \,49, \,50, \,
51, \,52, \,55, \,56, \,57, \,58, \,60, \hfil\dots,\hfil 69\}$,
\linebreak
we attribute new weights
$\wt(x)=\wwt_x$, \ $\wt(y)=\wwt_y$, \ $\wt(z)=\wwt_z$
(listed in the sixth column) to the variables $x,y,z$ and find out that the
terms of the lowest weight describe an isolated Du~Val singularity.
We note that the weights $\wt(x)$, $\wt(y)$, $\wt(z)$ coincide with those
of a Du~Val singularity of type given in the seventh column. Hence the
singularity of~$Y$ in the chosen neighbourhood is Du~Val of this type by
Theorem~2.10. Unfortunately, it is impossible to find such weights in
the remaining cases. But the corresponding equations have already been
considered in Lemma~4.14: the case of~$T_{37}$ is exactly case~1,
the case of~$T_{27}$ is case~2, the cases of~$T_{20}$, $T_{21}$,~$T_{29}$
fall into case~3, the cases of
$T_{32}$, $T_{38}$, $T_{39}$, $T_{42}$ and~$T_{48}$ correspond to case~4,
and the cases of~$T_{40}$, $T_{46}$,
$T_{47}$, $T_{53}$, $T_{54}$,~$T_{59}$ correspond to case~5 of Lemma~4.14.
In~either case, the singularity of~$V$ is Du~Val of type given in the seventh
\linebreak
column.

The singularity of~$V$ is locally isomorphic to the product of~$\CC$ and
the corresponding Du~Val singularity. Hence $V$ has canonical singularities
in the chosen neighbourhood and, by Corollary~2.5, $V$ has canonical
singularities.

\medskip

{

\raggedbottom

$$
\vbox{
\let\PLUS=+
\binoppenalty10000
\relpenalty10000
\catcode`\+=\active
\raggedright
\def+{\allowbreak\PLUS\nobreak\unskip}
\define\mystrut{\vphantom{$\Bigl($}}
\def\ensp{\kern0.35em}
\tabskip=0pt \offinterlineskip
\halign{\vrule\mystrut
\ensp\hfill#\hfill\ensp&\vrule# &
\ensp\hfill#\hfill\ensp&\vrule# &
\ensp\hfill#\hfill\ensp&\vrule# &
\ensp\hfill\tabboxone{#}\hfill\ensp&\vrule# &
\ensp\hfill\tabboxtwo{#}\hfill\ensp&\vrule# &
\ensp\hfill\tabboxthree{#}\hfill\ensp&\vrule# &
\ensp\hfill\tabboxthree{\hsize 1.2cm #}\hfill\ensp\vrule
\cr
\multispan{13}
\hfill{\smc Table 2}
\cr
\noalign{\medskip}
\noalign{\hrule}
\lw{2.5}{$T_i$}&&\lw{2.5}{$(d_1{,}d_2{,}d_3{,}d_4)$} &&\lw{2.5}{$b$}&&
\centerline{\lw{2.5}{Equation of~$V$}} &&
Local equation of~$V$
&& \centerline{\lw{2.5}{Weights}} && Singu\-larity
\cr
\noalign{\hrule}
$T_{4}$ && $(1,1,1,0)$ && $-1$ &&
$\alpha^1_2x_1^3+\alpha^2_2x_1^2x_2+\alpha^3_2x_1^2x_3+\alpha^4_2x_1^2x_2+
\alpha^5_2x_1x_2x_3+\alpha^6_2x_1x_3^2+\alpha^7_2x_2^3+
\alpha^8_2x_2^2x_3+\alpha^9_2x_2x_3^2+\alpha^{10}_2x_3^3+\alpha^1_1x_1^2x_4+
\alpha^2_1x_1x_2x_4+
\alpha^3_1x_1x_3x_4+\alpha^4_1x_2^2x_4+\alpha^5_1x_2x_3x_4+
\alpha^5_1x_3^2x_4+\alpha^1_0x_1x_4^2+\alpha^2_0x_2x_4^2+\alpha^3_0x_3x_4^2=0$
&& $P_1(x,y,z)+P_2(x,y,z)+P_3(x{,}\kkk y{,}\kkk z){=}\,0$ \vp\smash{($P_i$
is a} (general) homogeneous polynomial of degree~$i$) &&
\centerline{--} && Non-sin\-gu\-lar point
\cr
\noalign{\hrule}
$T_{5}$ && $(1,1,1,1)$ && $-2$ && \centerline{--} && \centerline{--} &&
\centerline{--} && \centerline{--}
\cr
\noalign{\hrule}
}}
$$
$$
\vbox{
\let\PLUS=+
\binoppenalty10000
\relpenalty10000
\catcode`\+=\active
\raggedright
\def+{\allowbreak\PLUS\nobreak\unskip}
\define\mystrut{\vphantom{$\Bigl($}}
\def\ensp{\kern0.41em}
\tabskip=0pt \offinterlineskip
\halign{\vrule\mystrut
\ensp\hfill#\hfill\ensp&\vrule# &
\ensp\hfill#\hfill\ensp&\vrule# &
\ensp\hfill#\hfill\ensp&\vrule# &
\ensp\hfill\tabboxone{#}\hfill\ensp&\vrule# &
\ensp\hfill\tabboxtwo{#}\hfill\ensp&\vrule# &
\ensp\hfill\tabboxthree{#}\hfill\ensp&\vrule# &
\ensp\hfill\tabboxthree{\hsize 1.2cm #}\hfill\ensp\vrule
\cr
\multispan{13}
\hfill{\smc Table 2 continued}
\cr
\noalign{\medskip}
\noalign{\hrule}
$T_{6}$ && $(2,1,0,0)$ && $-1$ &&
$\alpha_5x_1^3+\alpha_4x_1^2x_2+\alpha^1_3x_1^2x_3+
\alpha^2_3x_1^2x_4+ \alpha^3_3x_1x_2^2+\alpha^1_2x_1x_2x_3+
\alpha^2_2x_1x_2x_4+\alpha^1_1x_1x_3^2+
\!\alpha^2_1x_1x_3x_4\phantom{aaa}+
\alpha^3_1x_1x_4^2+\alpha^3_2x_2^3+\alpha^4_1x_2^2x_3+
\alpha^5_1x_2^2x_4+\alpha^1_0x_2x_3^2+\alpha^2_0x_2x_3x_4+
\alpha^3_0x_2x_4^2=0$ &&
$x^3+x^2y+x^2z+x^2+xy^2+xyz+xy+xz^2+xz+x+y^3+y^2z+y^2+yz^2+yz+y=
\nomathbreak 0$ &&
\centerline{--} && Non-sin\-gu\-lar point
\cr
\noalign{\hrule}
$T_{7}$ && $(2,1,1,0)$ && $-2$ &&
$\alpha_4x_1^3+\alpha^1_3x_1^2x_2+\alpha^2_3x_1^2x_3+
\alpha^1_2x_1^2x_4+\alpha^2_2x_1x_2^2+\alpha^3_2x_1x_2x_3+
\alpha^1_1x_1x_2x_4+\alpha^4_2x_1x_3^2+\alpha^2_1x_1x_3x_4+
\alpha^1_0x_1x_4^2+\alpha^3_1x_2^3+\alpha^4_1x_2^2x_3+\alpha^2_0x_2^2x_4+
\alpha^5_1x_2x_3^2+\alpha^3_0x_2x_3x_4+\alpha^6_1x_3^3+\alpha^4_0x_3^2x_4=0$
&& $x^3+x^2y+x^2z+x^2+xy^2+xyz+xy+xz^2+xz+x+y^3+y^2z+y^2+yz^2+yz+z^3+z^2=0$
&& \centerline{--} && Non-sin\-gu\-lar point
\cr
\noalign{\hrule}
$T_{8}$ && $(2,1,1,1)$ && $-3$ && \centerline{--} && \centerline{--} && \centerline{--} && \centerline{--}
\cr
\noalign{\hrule}
$T_{9}$ && $(2,2,0,0)$ && $-2$ &&
$\alpha^1_4x_1^3+\alpha^2_4x_1^2x_2+\alpha^1_2x_1^2x_3+
\alpha^2_2x_1^2x_4+\alpha^3_4x_1x_2^2+\alpha^3_2x_1x_2x_3+
\alpha^4_2x_1x_2x_4+\alpha^1_0x_1x_3^2+
\!\alpha^2_0x_1x_3x_4\phantom{aaa}+
\alpha^3_0x_1x_4^2+\alpha^4_4x_2^3+\alpha^5_2x_2^2x_3+
\alpha^6_2x_2^2x_4+\alpha^4_0x_2x_3^2+\alpha^5_0x_2x_3x_4+
\alpha^6_0x_2x_4^2=0$ &&
$x^3+x^2y+x^2z+x^2+xy^2+xyz+xy+xz^2+xz+x+y^3+y^2z+y^2+yz^2+yz+y=
\nomathbreak 0$ &&
\centerline{--} && Non-sin\-gu\-lar point
\cr
\noalign{\hrule}
}}
$$
$$
\vbox{
\let\PLUS=+
\binoppenalty10000
\relpenalty10000
\catcode`\+=\active
\raggedright
\def+{\allowbreak\PLUS\nobreak\unskip}
\define\mystrut{\vphantom{$\Bigl($}}
\def\ensp{\kern0.41em}
\tabskip=0pt \offinterlineskip
\halign{\vrule\mystrut
\ensp\hfill#\hfill\ensp&\vrule# &
\ensp\hfill#\hfill\ensp&\vrule# &
\ensp\hfill#\hfill\ensp&\vrule# &
\ensp\hfill\tabboxone{#}\hfill\ensp&\vrule# &
\ensp\hfill\tabboxtwo{#}\hfill\ensp&\vrule# &
\ensp\hfill\tabboxthree{#}\hfill\ensp&\vrule# &
\ensp\hfill\tabboxthree{\hsize 1.2cm #}\hfill\ensp\vrule
\cr
\multispan{13}
\hfill{\smc Table 2 continued}
\cr
\noalign{\medskip}
\noalign{\hrule}
$T_{10}$ && $(2,2,1,0)$ && $-3$ &&
$\alpha^1_3x_1^3+\alpha^2_3x_1^2x_2+\alpha^1_2x_1^2x_3+
\alpha^1_1x_1^2x_4+\alpha^3_3x_1x_2^2+\alpha^2_2x_1x_2x_3+
\alpha^2_1x_1x_2x_4+\alpha^3_1x_1x_3^2+\alpha^1_0x_1x_3x_4+
\alpha^4_3x_2^3+\alpha^3_2x_2^2x_3+\alpha^4_1x_2^2x_4+
\alpha^5_1x_2x_3^2+\alpha^2_0x_2x_3x_4+\alpha^3_0x_3^3=0$ &&
$P_3(x,y)+P^1_2(x,y)z+P^2_2(x,y)+P_1(x{,}\kkk y)z{=}\, 0$
\vp($P_i$, \smash{$P^j_i$} are homogeneous polynomials of degree~$i$) &&
$\wwt_x=1$ $\wwt_y=1$ $\wwt_z=1$ && \centerline{$\Bbb{A}_1$}
\cr
\noalign{\hrule}
$T_{11}$ && $(2,2,1,1)$ && $-4$ &&
$\alpha^1_2x_1^3+\alpha^2_2x_1^2x_2+\alpha^1_1x_1^2x_3+
\alpha^2_1x_1^2x_4+\alpha^3_2x_1x_2^2+\alpha^3_1x_1x_2x_3+
\alpha^4_1x_1x_2x_4+\alpha^1_0x_1x_3^2+
\!\alpha^2_0x_1x_3x_4\phantom{aaa}+
\alpha^3_0x_1x_4^2+\alpha^4_2x_2^3+\alpha^5_1x_2^2x_3+
\alpha^6_1x_2^2x_4+\alpha^4_0x_2x_3^2+\alpha^5_0x_2x_3x_4+
\alpha^6_0x_2x_4^2=0$ &&
$x^3+x^2y+x^2z+x^2+xy^2+xyz+xy+xz^2+xz+x+y^3+y^2z+y^2+yz^2+yz+y=
\nomathbreak 0$ &&
\centerline{--} && Non-sin\-gu\-lar point
\cr
\noalign{\hrule}
$T_{12}$ && $(2,2,2,0)$ && $-4$ &&
$\alpha^1_2x_1^3+\alpha^2_2x_1^2x_2+\alpha^3_2x_1^2x_3+
\alpha^1_0x_1^2x_4+\alpha^4_2x_1x_2^2+\alpha^5_2x_1x_2x_3+
\alpha^2_0x_1x_2x_4+\alpha^6_2x_1x_3^2+\alpha^3_0x_1x_3x_4+
\alpha^7_2x_2^3+\alpha^8_2x_2^2x_3+\alpha^4_0x_2^2x_4+
\alpha^9_2x_2x_3^2+\alpha^5_0x_2x_3x_4+\alpha^{10}_2x_3^3+
\alpha^6_0x_3^2x_4{=}\,0$ && $P_3(x,y,z)+P_2(x{,}\kkk y{,}\kkk z){=}\,0$
\vp\smash{($P_i$} is a homogeneous polynomial of degree~$i$) &&
$\wwt_x=1$ $\wwt_y=1$ $\wwt_z=1$ && \centerline{$\Bbb{A}_1$}
\cr
\noalign{\hrule}
$T_{13}$ && $(2,2,2,1)$ && $-5$ &&
$\alpha^1_1x_1^3+\alpha^2_1x_1^2x_2+\alpha^3_1x_1^2x_3+
\alpha^1_0x_1^2x_4+\alpha^4_1x_1x_2^2+\alpha^5_1x_1x_2x_3+
\alpha^2_0x_1x_2x_4+\alpha^6_1x_1x_3^2+\alpha^3_0x_1x_3x_4+
\alpha^7_1x_2^3+\alpha^8_1x_2^2x_3+\alpha^4_0x_2^2x_4+
\alpha^9_1x_2x_3^2+\alpha^5_0x_2x_3x_4+\alpha^{10}_1x_3^3+
\alpha^6_0x_3^2x_4{=}\,0$ && $P_3(x,y,z)+P_2(x{,}\kkk y{,}\kkk z){=}\,0$
\vp\smash{($P_i$} is a homogeneous polynomial of degree~$i$) &&
$\wwt_x=1$ $\wwt_y=1$ $\wwt_z=1$ && \centerline{$\Bbb{A}_1$}
\cr
\noalign{\hrule}
}}
$$
$$
\vbox{
\let\PLUS=+
\binoppenalty10000
\relpenalty10000
\catcode`\+=\active
\raggedright
\def+{\allowbreak\PLUS\nobreak\unskip}
\define\mystrut{\vphantom{$\Bigl($}}
\def\ensp{\kern0.41em}
\tabskip=0pt \offinterlineskip
\halign{\vrule\mystrut
\ensp\hfill#\hfill\ensp&\vrule# &
\ensp\hfill#\hfill\ensp&\vrule# &
\ensp\hfill#\hfill\ensp&\vrule# &
\ensp\hfill\tabboxone{#}\hfill\ensp&\vrule# &
\ensp\hfill\tabboxtwo{#}\hfill\ensp&\vrule# &
\ensp\hfill\tabboxthree{#}\hfill\ensp&\vrule# &
\ensp\hfill\tabboxthree{\hsize 1.2cm #}\hfill\ensp\vrule
\cr
\multispan{13}
\hfill{\smc Table 2 continued}
\cr
\noalign{\medskip}
\noalign{\hrule}
$T_{14}$ && $(2,2,2,2)$ && $-6$ && \centerline{--} && \centerline{--} &&
\centerline{--} && \centerline{--}
\cr
\noalign{\hrule}
$T_{15}$ && $(3,1,0,0)$ && $-2$ &&
$\alpha_7x_1^3+\alpha_5x_1^2x_2+\alpha^1_4x_1^2x_3+
\alpha^2_4x_1^2x_4+\alpha_3x_1x_2^2+\alpha^1_2x_1x_2x_3+
\alpha^2_2x_1x_2x_4+\alpha^1_1x_1x_3^2+
\alpha^2_1x_1x_3x_4+\alpha^3_1x_1x_4^2+\alpha^4_1x_2^3+\alpha^1_0x_2^2x_3+
\alpha^2_0x_2^2x_4=0$ &&
$x^3+x^2y+x^2z+x^2+xy^2+xyz+xy+xz^2+xz+x+y^3+y^2z+y^2=\nomathbreak
0$ &&
\centerline{--} && Non-sin\-gu\-lar point
\cr
\noalign{\hrule}
$T_{16}$ && $(3,1,1,0)$ && $-3$ &&
$\alpha_6x_1^3+\alpha^1_4x_1^2x_2+\alpha^2_4x_1^2x_3+\alpha_3x_1^2x_4+
\alpha_2^1x_1x_2^2+\alpha_2^2x_1x_2x_3+\alpha_1^1x_1x_2x_4+
\alpha_2^3x_1x_3^2+\alpha_1^2x_1x_3x_4+\alpha_0^1x_1x_4^2+
\alpha_0^2x_2^3+\alpha_0^3x_2^2x_3+\alpha_0^4x_2x_3^2+
\alpha_0^5x_3^3=0$ && $xQ(x,y,z)+P_3(y,z)=0$ \vp($Q(0)\neq 0$,
\smash{$P_3$} is a homogeneous polynomial of degree~$3$)&&
\centerline{--} && Non-sin\-gu\-lar point
\cr
\noalign{\hrule}
$T_{17}$ && $(3,2,0,0)$ && $-3$ &&
$\alpha_6x_1^3+\alpha_5x_1^2x_2+\alpha_3^1x_1^2x_3+\alpha_3^2x_1^2x_4+
\alpha_4x_1x_2^2+\alpha_2^1x_1x_2x_3+\alpha_2^2x_1x_2x_4+
\alpha_0^1x_1x_3^2+\alpha_0^2x_1x_3x_4+\alpha_0^3x_1x_4^2+
\alpha_3^3x_2^3+\alpha_1^1x_2^2x_3+\alpha_1^2x_2^2x_4=0$ &&
$xQ_1(x,y,z)\!+y^2Q_2\kern-.5mm(y{,}\kkk z){=}\,0$ ($Q_1(0)\neq 0$)&&
\centerline{--} && Non-sin\-gu\-lar point
\cr
\noalign{\hrule}
$T_{18}$ && $(3,2,1,0)$ && $-4$ &&
$\alpha_5x_1^3+\alpha_4x_1^2x_2+\alpha_3^1x_1^2x_3+\alpha^1_2x_1^2x_4+
\alpha^2_3x_1x_2^2+\alpha^2_2x_1x_2x_3+\alpha^1_1x_1x_2x_4+
\alpha^2_1x_1x_3^2+\alpha^1_0x_1x_3x_4+\alpha_2^3x_2^3+
\alpha_1^3x_2^2x_3+\alpha_0^2x_2^2x_4+\alpha_0^3x_2x_3^2=0$ &&
$x^2Q_1(x,y,z)+xyQ_2(x,y,z)+xzQ_3(x,y,z)+y^2Q_4(x,y,z)+yz^2=0$
($Q_i(0)\neq 0$) &&
$\wwt_x=1$ $\wwt_y=1$ $\wwt_z=1$ && \centerline{$\Bbb{A}_1$}
\cr
\noalign{\hrule}
}}
$$
$$
\vbox{
\let\PLUS=+
\binoppenalty10000
\relpenalty10000
\catcode`\+=\active
\raggedright
\def+{\allowbreak\PLUS\nobreak\unskip}
\define\mystrut{\vphantom{$\Bigl($}}
\def\ensp{\kern0.41em}
\tabskip=0pt \offinterlineskip
\halign{\vrule\mystrut
\ensp\hfill#\hfill\ensp&\vrule# &
\ensp\hfill#\hfill\ensp&\vrule# &
\ensp\hfill#\hfill\ensp&\vrule# &
\ensp\hfill\tabboxone{#}\hfill\ensp&\vrule# &
\ensp\hfill\tabboxtwo{#}\hfill\ensp&\vrule# &
\ensp\hfill\tabboxthree{#}\hfill\ensp&\vrule# &
\ensp\hfill\tabboxthree{\hsize 1.2cm #}\hfill\ensp\vrule
\cr
\multispan{13}
\hfill{\smc Table 2 continued}
\cr
\noalign{\medskip}
\noalign{\hrule}
$T_{19}$ && $(3,2,1,1)$ && $-5$ &&
$\alpha_4x_1^3+\alpha_3x_1^2x_2+\alpha^1_2x_1^2x_3+\alpha^2_2x_1^2x_4+
\alpha^3_2x_1x_2^2+\alpha^1_1x_1x_2x_3+\alpha^2_1x_1x_2x_4+\alpha^1_0x_1x_3^2+
\alpha^2_0x_1x_3x_4+\alpha^3_0x_1x_4^2+\alpha^3_1x_2^3+
\alpha^4_0x_2^2x_3+\alpha^5_0x_2^2x_4=0$ &&
$xQ_1(x,y,z)+Q_2(y,z)=0$ ($Q_1(0)\neq 0$) &&
\centerline{--} && Non-sin\-gu\-lar point
\cr
\noalign{\hrule}
$T_{20}$ && $(3,2,2,0)$ && $-5$ &&
$\alpha_4x_1^3+\alpha^1_3x_1^2x_2+\alpha^2_3x_1^2x_3+
\alpha^1_1x_1^2x_4+\alpha^1_2x_1x_2^2+\alpha^2_2x_1x_2x_3+
\alpha^1_0x_1x_2x_4+\alpha^3_2x_1x_3^2+\alpha^2_0x_1x_3x_4+
\alpha^2_1x_2^3+\alpha^3_1x_2^2x_3+\alpha^4_1x_2x_3^2+
\alpha^5_1x_3^3=0$ &&
$x^3+x^2y+x^2z+x^2+xy^2+xyz+xy+xz^2+xz+y^3+y^2z+yz^2+z^3=0$ &&
\centerline{--} && \centerline{$\Bbb{A}_3$}
\cr
\noalign{\hrule}
$T_{21}$ && $(3,2,2,1)$ && $-6$ &&
$\alpha_3x_1^3+\alpha^1_2x_1^2x_2+\alpha^2_2x_1^2x_3+
\alpha^1_1x_1^2x_4+\alpha^2_1x_1x_2^2+\alpha^3_1x_1x_2x_3+
\alpha^1_0x_1x_2x_4+\alpha^4_1x_1x_3^2+
\alpha^2_0x_1x_3x_4+\alpha^3_0x_2^3+\alpha^4_0x_2^2x_3+\alpha^5_0x_2x_3^2+
\alpha^6_0x_3^3=0$ &&
$x^3+x^2y+x^2z+x^2+xy^2+xyz+xy+xz^2+xz+y^3+y^2z+yz^2+z^3=0$ &&
\centerline{--} && \centerline{$\Bbb{A}_3$}
\cr
\noalign{\hrule}
$T_{22}$ && $(3,3,1,0)$ && $-5$ &&
$\alpha^1_4x_1^3+\alpha^2_4x_1^2x_2+\alpha^1_2x_1^2x_3+
\alpha^1_1x_1^2x_4+\alpha^3_4x_1x_2^2+\alpha^2_2x_1x_2x_3+
\!\alpha^2_1x_1x_2x_4\phantom{aaa}
+\alpha^1_0x_1x_3^2+\alpha^4_4x_2^3+
\alpha^3_2x_2^2x_3+\alpha^3_1x_2^2x_4+\alpha^2_0x_2x_3^2{=}\,0$ &&
$x^3+x^2y+x^2z+x^2+xy^2+xyz+xy+xz^2+y^3+y^2z+y^2+yz^2=0$ &&
$\wwt_x=2$ $\wwt_y=2$ $\wwt_z=1$ && \centerline{$\Bbb{A}_3$}
\cr
\noalign{\hrule}
}}
$$
$$
\vbox{
\let\PLUS=+
\binoppenalty10000
\relpenalty10000
\catcode`\+=\active
\raggedright
\def+{\allowbreak\PLUS\nobreak\unskip}
\define\mystrut{\vphantom{$\Bigl($}}
\def\ensp{\kern0.41em}
\tabskip=0pt \offinterlineskip
\halign{\vrule\mystrut
\ensp\hfill#\hfill\ensp&\vrule# &
\ensp\hfill#\hfill\ensp&\vrule# &
\ensp\hfill#\hfill\ensp&\vrule# &
\ensp\hfill\tabboxone{#}\hfill\ensp&\vrule# &
\ensp\hfill\tabboxtwo{#}\hfill\ensp&\vrule# &
\ensp\hfill\tabboxthree{#}\hfill\ensp&\vrule# &
\ensp\hfill\tabboxthree{\hsize 1.2cm #}\hfill\ensp\vrule
\cr
\multispan{13}
\hfill{\smc Table 2 continued}
\cr
\noalign{\medskip}
\noalign{\hrule}
$T_{23}$ && $(3,3,2,0)$ && $-6$ &&
$\alpha^1_3x_1^3+\alpha^2_3x_1^2x_2+\alpha^1_2x_1^2x_3+
\alpha^1_0x_1^2x_4+\alpha^3_3x_1x_2^2+\alpha^2_2x_1x_2x_3+
\alpha^2_0x_1x_2x_4+\alpha^1_1x_1x_3^2+\alpha^4_3x_2^3+
\alpha^3_2x_2^2x_3+\alpha^3_0x_2^2x_4+\alpha^2_1x_2x_3^2+
\alpha^4_0x_3^3=\nomathbreak 0$ &&
$x^3+x^2y+x^2z+x^2+xy^2+xyz+xy+xz^2+y^3+y^2z+y^2+yz^2+z^3=\nomathbreak
0$ &&
$\wwt_x=3$ $\wwt_y=3$ $\wwt_z=2$ && \centerline{$\Bbb{A}_2$}
\cr
\noalign{\hrule}
$T_{24}$ && $(3,3,2,1)$ && $-7$ &&
$\alpha^1_2x_1^3+\alpha^2_2x_1^2x_2+\alpha^1_1x_1^2x_3+
\alpha^1_0x_1^2x_4+\alpha^3_2x_1x_2^2+\alpha^2_1x_1x_2x_3+
\alpha^2_0x_1x_2x_4+\alpha^3_0x_1x_3^2+\alpha^4_2x_2^3+
\alpha^3_1x_2^2x_3+\alpha^4_0x_2^2x_4+\alpha^5_0x_2x_3^2{=}\,0$ &&
$x^3+x^2y+x^2z+x^2+xy^2+xyz+xy +xz^2+y^3+y^2z+ y^2+yz^2=0$ &&
$\wwt_x=2$ $\wwt_y=2$ $\wwt_z=1$ && \centerline{$\Bbb{A}_3$}
\cr
\noalign{\hrule}
$T_{25}$ && $(4,1,0,0)$ && $-3$ &&
$\alpha_9x_1^3+\alpha_6x_1^2x_2+\alpha^1_5x_1^2x_3+
\alpha^2_5x_1^2x_4+\alpha_3x_1x_2^2+
\alpha^1_2x_1x_2x_3+\alpha^2_2x_1x_2x_4+\alpha^1_1x_1x_3^2+
\alpha^2_1x_1x_3x_4+\alpha^3_1x_1x_4^2+\alpha_0x_2^3=0$ &&
$xQ(x,y,z)+y^3=0$ ($Q(0)\neq 0$) &&
\centerline{--} && Non-sin\-gular point
\cr
\noalign{\hrule}
$T_{26}$ && $(4,2,0,0)$ && $-4$ &&
$\alpha_8x_1^3+\alpha_6x_1^2x_2+\alpha^1_4x_1^2x_3+
\alpha^2_4x_1^2x_4+\alpha^3_4x_1x_2^2+\alpha^1_2x_1x_2x_3+
\alpha^2_2x_1x_2x_4+\alpha^1_0x_1x_3^2+\alpha^2_0x_1x_3x_4+
\alpha^3_0x_1x_4^2+\alpha^3_2x_2^3+\alpha^4_0x_2^2x_3+\alpha^5_0x_2^2x_4=0$
&& $xQ(x,y,z)+y^2+y^3=0$ ($Q(0)\neq 0$) &&
\centerline{--} && Non-sin\-gu\-lar point
\cr
\noalign{\hrule}
$T_{27}$ && $(4,2,1,0)$ && $-5$ &&
$\alpha_7x_1^3+\alpha_5x_1^2x_2+\alpha_4x_1^2x_3+
\alpha^1_3x_1^2x_4+\alpha^2_3x_1x_2^2+\alpha_2x_1x_2x_3+
\alpha^1_1x_1x_2x_4+\alpha^2_1x_1x_3^2+\alpha^1_0x_1x_3x_4+
\alpha^3_1x_2^3+\alpha^2_0x_2^2x_3=0$ &&
$x^3+x^2y+x^2z+x^2+xy^2+xyz+xy+xz^2+xz+y^3+y^2z=0$ &&
\centerline{--} && \centerline{$\Bbb{A}_3$}
\cr
\noalign{\hrule}
}}
$$
$$
\vbox{
\let\PLUS=+
\binoppenalty10000
\relpenalty10000
\catcode`\+=\active
\raggedright
\def+{\allowbreak\PLUS\nobreak\unskip}
\define\mystrut{\vphantom{$\Bigl($}}
\def\ensp{\kern0.41em}
\tabskip=0pt \offinterlineskip
\halign{\vrule\mystrut
\ensp\hfill#\hfill\ensp&\vrule# &
\ensp\hfill#\hfill\ensp&\vrule# &
\ensp\hfill#\hfill\ensp&\vrule# &
\ensp\hfill\tabboxone{#}\hfill\ensp&\vrule# &
\ensp\hfill\tabboxtwo{#}\hfill\ensp&\vrule# &
\ensp\hfill\tabboxthree{#}\hfill\ensp&\vrule# &
\ensp\hfill\tabboxthree{\hsize 1.2cm #}\hfill\ensp\vrule
\cr
\multispan{13}
\hfill{\smc Table 2 continued}
\cr
\noalign{\medskip}
\noalign{\hrule}
$T_{28}$ && $(4,2,1,1)$ && $-6$ &&
$\alpha_6x_1^3+\alpha_4x_1^2x_2+\alpha^1_3x_1^2x_3+
\alpha^2_3x_1^2x_4+\alpha_2x_1x_2^2+\alpha^1_1x_1x_2x_3+
\alpha^2_1x_1x_2x_4+\alpha^1_0x_1x_3^2+\alpha^2_0x_1x_3x_4+
\alpha^3_0x_1x_4^2+\alpha^4_0x_2^3=0$ && $xQ(x,y,z)+y^3=0$
($Q(0)\neq 0$) &&
\centerline{--} && Non-sin\-gu\-lar point
\cr
\noalign{\hrule}
$T_{29}$ && $(4,2,2,0)$ && $-6$ &&
$\alpha_6x_1^3+\alpha^1_4x_1^2x_2+\alpha^2_4x_1^2x_3+
\alpha^1_2x_1^2x_4+\alpha^2_2x_1x_2^2+\alpha^3_2x_1x_2x_3+
\alpha^1_0x_1x_2x_4+\alpha^4_2x_1x_3^2+\alpha^2_0x_1x_3x_4+
\alpha^3_0x_2^3+\alpha^4_0x_2^2x_3+\alpha^5_0x_2x_3^2+
\alpha^6_0x_3^3=0$ &&
$x^3+x^2y+x^2z+x^2+xy^2+xyz+xy+xz^2+xz+y^3+y^2z+yz^2+z^3=0$ &&
\centerline{--} && \centerline{$\Bbb{A}_3$}
\cr
\noalign{\hrule}
$T_{30}$ && $(4,3,1,0)$ && $-6$ &&
$\alpha_6x_1^3+\alpha_5x_1^2x_2+\alpha^1_3x_1^2x_3+
\alpha^1_2x_1^2x_4+\alpha_4x_1x_2^2+\alpha^2_2x_1x_2x_3+
\alpha^1_1x_1x_2x_4+\alpha^1_0x_1x_3^2+\alpha^2_3x_2^3+
\alpha^2_1x_2^2x_3+\alpha^2_0x_2^2x_4=0$ &&
$x^3+x^2y+x^2z+x^2+xy^2+xyz+xy+xz^2+y^3+y^2z+y^2=0$ &&
$\wwt_x=2$ $\wwt_y=2$ $\wwt_z=1$ && \centerline{$\Bbb{A}_3$}
\cr
\noalign{\hrule}
$T_{31}$ && $(4,3,2,0)$ && $-7$ &&
$\alpha_5x_1^3+\alpha_4x_1^2x_2+\alpha^1_3x_1^2x_3+
\alpha^1_1x_1^2x_4+\alpha^2_3x_1x_2^2+\alpha^1_2x_1x_2x_3+
\alpha^1_0x_1x_2x_4+
\alpha^2_1x_1x_3^2+\alpha^2_2x_2^3+\alpha^3_1x_2^2x_3+\alpha^2_0x_2x_3^2=0$
&& $x^3+x^2y+x^2z+x^2+xy^2+xyz+xy+xz^2+y^3+y^2z+yz^2\mathbreak =0$ &&
$\wwt_x=2$ $\wwt_y=2$ $\wwt_z=1$ && \centerline{$\Bbb{A}_3$}
\cr
\noalign{\hrule}
$T_{32}$ && $(4,3,2,1)$ && $-8$ &&
$\alpha_4x_1^3+\alpha_3x_1^2x_2+\alpha^1_2x_1^2x_3+
\alpha^1_1x_1^2x_4+\alpha^2_2x_1x_2^2+\alpha^2_1x_1x_2x_3+
\alpha^1_0x_1x_2x_4+\alpha^2_0x_1x_3^2+\alpha^3_1x_2^3+
\alpha^3_0x_2^2x_3=\nomathbreak 0$ &&
$x^3+x^2y+x^2z+x^2+xy^2+xyz+xy+xz^2+y^3+y^2z=0$ &&
\centerline{--} && \centerline{$\Bbb{A}_4$}
\cr
\noalign{\hrule}
}}
$$
$$
\vbox{
\let\PLUS=+
\binoppenalty10000
\relpenalty10000
\catcode`\+=\active
\raggedright
\def+{\allowbreak\PLUS\nobreak\unskip}
\define\mystrut{\vphantom{$\Bigl($}}
\def\ensp{\kern0.41em}
\tabskip=0pt \offinterlineskip
\halign{\vrule\mystrut
\ensp\hfill#\hfill\ensp&\vrule# &
\ensp\hfill#\hfill\ensp&\vrule# &
\ensp\hfill#\hfill\ensp&\vrule# &
\ensp\hfill\tabboxone{#}\hfill\ensp&\vrule# &
\ensp\hfill\tabboxtwo{#}\hfill\ensp&\vrule# &
\ensp\hfill\tabboxthree{#}\hfill\ensp&\vrule# &
\ensp\hfill\tabboxthree{\hsize 1.2cm #}\hfill\ensp\vrule
\cr
\multispan{13}
\hfill{\smc Table 2 continued}
\cr
\noalign{\medskip}
\noalign{\hrule}
$T_{33}$ && $(4,3,3,0)$ && $-8$ &&
$\alpha_4x_1^3+\alpha^1_3x_1^2x_2+\alpha^2_3x_1^2x_3+
\alpha_0x_1^2x_4+\alpha^1_2x_1x_2^2+\alpha^2_2x_1x_2x_3+
\alpha^3_2x_1x_3^2+\alpha^1_1x_2^3+\alpha^2_1x_2^2x_3+
\alpha^3_1x_2x_3^2+\alpha^4_1x_3^3=0$ &&
$x^3+x^2y+x^2z+x^2+xy^2+xyz+xz^2+y^3+y^2z+yz^2+z^3=0$ &&
$\wwt_x=3$ $\wwt_y=2$ $\wwt_z=2$ && \centerline{$\Bbb{D}_4$}
\cr
\noalign{\hrule}
$T_{34}$ && $(4,3,3,1)$ && $-9$ &&
$\alpha_3x_1^3+\alpha^1_2x_1^2x_2+\alpha^2_2x_1^2x_3+
\alpha^1_0x_1^2x_4+\alpha^1_1x_1x_2^2+\alpha^2_1x_1x_2x_3+
\!\alpha^3_1x_1x_3^2+\alpha^2_0x_2^3+\alpha^3_0x_2^2x_3+
\alpha^4_0x_2x_3^2+\alpha^5_0x_3^3=0$ &&
$x^3+x^2y+x^2z+x^2+xy^2+xyz+xz^2+y^3+y^2z+yz^2+z^3=0$ &&
$\wwt_x=3$ $\wwt_y=2$ $\wwt_z=2$ && \centerline{$\Bbb{D}_4$}
\cr
\noalign{\hrule}
$T_{35}$ && $(4,4,2,0)$ && $-8$ &&
$\alpha^1_4x_1^3+\alpha^2_4x_1^2x_2+\alpha^1_2x_1^2x_3+
\alpha^1_0x_1^2x_4+\alpha^3_4x_1x_2^2+\alpha^2_2x_1x_2x_3+
\alpha^2_0x_1x_2x_4+\alpha^3_0x_1x_3^2+\alpha^4_4x_2^3+
\alpha^3_2x_2^2x_3+\alpha^4_0x_2^2x_4+\alpha^5_0x_2x_3^2=\nomathbreak 0$ &&
$x^3+x^2y+x^2z+x^2+xy^2+xyz+xy+xz^2+y^3+y^2z+y^2+yz^2=0$ &&
$\wwt_x=2$ $\wwt_y=2$ $\wwt_z=1$ && \centerline{$\Bbb{A}_3$}
\cr
\noalign{\hrule}
$T_{36}$ && $(5,2,0,0)$ && $-5$ &&
$\alpha_{10}x_1^3+\alpha_7x_1^2x_2+\alpha^1_5x_1^2x_3+
\alpha^2_5x_1^2x_4+\alpha_4x_1x_2^2+\alpha^1_2x_1x_2x_3+
\alpha^2_2x_1x_2x_4+\alpha^1_0x_1x_3^2+\alpha^2_0x_1x_3x_4+
\alpha^3_0x_1x_4^2+\alpha_1x_2^3=0$ && $xQ(x,y,z)+y^3=0$ ($Q(0)\neq
0$) &&
\centerline{--} && Non-sin\-gu\-lar point
\cr
\noalign{\hrule}
$T_{37}$ && $(5,2,1,0)$ && $-6$ &&
$\alpha_9x_1^3+\alpha_6x_1^2x_2+\alpha_5x_1^2x_3+
\alpha_4x_1^2x_4+\alpha_3x_1x_2^2+
\alpha_2x_1x_2x_3+
\!\alpha^1_1x_1x_2x_4\phantom{aaa}+
\alpha^2_1x_1x_3^2+
\alpha^1_0x_1x_3x_4+\alpha^2_0x_2^3=\nomathbreak 0$ &&
$x^3+x^2y+x^2z+x^2+xy^2+xyz+xy+xz^2+xz+y^3=\nomathbreak 0$ &&
\centerline{--} && \centerline{$\Bbb{A}_2$}
\cr
\noalign{\hrule}
}}
$$
$$
\vbox{
\let\PLUS=+
\binoppenalty10000
\relpenalty10000
\catcode`\+=\active
\raggedright
\def+{\allowbreak\PLUS\nobreak\unskip}
\define\mystrut{\vphantom{$\Bigl($}}
\def\ensp{\kern0.39em}
\tabskip=0pt \offinterlineskip
\halign{\vrule\mystrut
\ensp\hfill#\hfill\ensp&\vrule# &
\ensp\hfill#\hfill\ensp&\vrule# &
\ensp\hfill#\hfill\ensp&\vrule# &
\ensp\hfill\tabboxone{#}\hfill\ensp&\vrule# &
\ensp\hfill\tabboxtwo{#}\hfill\ensp&\vrule# &
\ensp\hfill\tabboxthree{#}\hfill\ensp&\vrule# &
\ensp\hfill\tabboxthree{\hsize 1.2cm #}\hfill\ensp\vrule
\cr
\multispan{13}
\hfill{\smc Table 2 continued}
\cr
\noalign{\medskip}
\noalign{\hrule}
$T_{38}$ && $(5,3,1,0)$ && $-7$ &&
$\alpha_8x_1^3+\alpha_6x_1^2x_2+\alpha^1_4x_1^2x_3+
\alpha_3x_1^2x_4+\alpha^2_4x_1x_2^2+\alpha^1_2x_1x_2x_3+
\alpha_1x_1x_2x_4+\alpha^1_0x_1x_3^2+\alpha^2_2x_2^3+
\alpha^2_0x_2^2x_3=\nomathbreak 0$ &&
$x^3+x^2y+x^2z+x^2+xy^2+xyz+xy+xz^2+y^3+y^2z=0$ &&
\centerline{--} && \centerline{$\Bbb{A}_4$}
\cr
\noalign{\hrule}
$T_{39}$ && $(5,3,2,0)$ && $-8$ &&
$\alpha_7x_1^3+\alpha_5x_1^2x_2+\alpha_4x_1^2x_3+
\alpha^1_2x_1^2x_4+\alpha_3x_1x_2^2+\alpha^2_2x_1x_2x_3+\alpha^1_0x_1x_2x_4+
\alpha^1_1x_1x_3^2+\alpha^2_1x_2^3+\alpha^2_0x_2^2x_3=\nomathbreak
0$ && $x^3+x^2y+x^2z+x^2+xy^2+xyz+xy+xz^2+y^3+y^2z=0$ &&
\centerline{--} && \centerline{$\Bbb{A}_4$}
\cr
\noalign{\hrule}
$T_{40}$ && $(5,3,2,1)$ && $-9$ &&
$\alpha_6x_1^3+\alpha_4x_1^2x_2+\alpha_3x_1^2x_3+
\alpha^1_2x_1^2x_4+\alpha^2_2x_1x_2^2+
\alpha_1x_1x_2x_3+\alpha^1_0x_1x_2x_4+\alpha^2_0x_1x_3^2+\alpha^3_0x_2^3=0$
&& $x^3+x^2y+x^2z+x^2+xy^2+xyz+xy+xz^2+y^3\,{=}\,0$ &&
\centerline{--} && \centerline{$\Bbb{A}_5$}
\cr
\noalign{\hrule}
$T_{41}$ && $(5,3,3,0)$ && $-9$ &&
$\alpha_6x_1^3+\alpha^1_4x_1^2x_2+\alpha^2_4x_1^2x_3+
\alpha_1x_1^2x_4+\alpha^1_2x_1x_2^2+
\alpha^2_2x_1x_2x_3+\alpha^3_2x_1x_3^2+\alpha^1_0x_2^3+\alpha^2_0x_2^2x_3+
\alpha^3_0x_2x_3^2+\alpha^4_0x_3^3=0$ &&
$x^3+x^2y+x^2z+x^2+xy^2+xyz+xz^2+y^3+y^2z+yz^2+z^3=0$ &&
$\wwt_x=3$ $\wwt_y=2$ $\wwt_z=2$ && \centerline{$\Bbb{D}_4$}
\cr
\noalign{\hrule}
$T_{42}$ && $(5,4,2,0)$ && $-9$ &&
$\alpha_6x_1^3+\alpha_5x_1^2x_2+\alpha^1_3x_1^2x_3+
\alpha^1_1x_1^2x_4+\alpha_4x_1x_2^2+\alpha_2x_1x_2x_3+
\alpha^1_0x_1x_2x_4+\alpha^2_0x_1x_3^2+\alpha^2_3x_2^3+
\alpha^2_1x_2^2x_3=\nomathbreak 0$ &&
$x^3+x^2y+x^2z+x^2+xy^2+xyz+xy+xz^2+y^3+y^2z=0$ &&
\centerline{--} && \centerline{$\Bbb{A}_4$}
\cr
\noalign{\hrule}
$T_{43}$ && $(5,4,3,0)$ && $-10$ &&
$\alpha_5x_1^3+\alpha_4x_1^2x_2+\alpha^1_3x_1^2x_3+
\alpha^1_0x_1^2x_4+\alpha^2_3x_1x_2^2+\alpha^1_2x_1x_2x_3+
\alpha^1_1x_1x_3^2+\alpha^2_2x_2^3+\alpha^2_1x_2^2x_3+
\alpha^2_0x_2x_3^2=\nomathbreak 0$ &&
$x^3+x^2y+x^2z+x^2+xy^2+xyz+xz^2+y^3+y^2z+yz^2=0$ &&
$\wwt_x=3$ $\wwt_y=2$ $\wwt_z=2$ && \centerline{$\Bbb{D}_4$}
\cr
\noalign{\hrule}
}}
$$
$$
\vbox{
\let\PLUS=+
\binoppenalty10000
\relpenalty10000
\catcode`\+=\active
\raggedright
\def+{\allowbreak\PLUS\nobreak\unskip}
\define\mystrut{\vphantom{$\Bigl($}}
\def\ensp{\kern0.39em}
\tabskip=0pt \offinterlineskip
\halign{\vrule\mystrut
\ensp\hfill#\hfill\ensp&\vrule# &
\ensp\hfill#\hfill\ensp&\vrule# &
\ensp\hfill#\hfill\ensp&\vrule# &
\ensp\hfill\tabboxone{#}\hfill\ensp&\vrule# &
\ensp\hfill\tabboxtwo{#}\hfill\ensp&\vrule# &
\ensp\hfill\tabboxthree{#}\hfill\ensp&\vrule# &
\ensp\hfill\tabboxthree{\hsize 1.2cm #}\hfill\ensp\vrule
\cr
\multispan{13}
\hfill{\smc Table 2 continued}
\cr
\noalign{\medskip}
\noalign{\hrule}
$T_{44}$ && $(5,4,3,1)$ && $-11$ &&
$\alpha_4x_1^3+\alpha_3x_1^2x_2+\alpha^1_2x_1^2x_3+
\alpha^1_0x_1^2x_4+\alpha^2_2x_1x_2^2+
\alpha^1_1x_1x_2x_3+\alpha^2_0x_1x_3^2+\alpha^2_1x_2^3+\alpha^3_0x_2^2x_3=0$
&& $x^3+x^2y+x^2z+x^2+xy^2+xyz+xz^2+y^3+y^2z=0$ &&
$\wwt_x=4$ $\wwt_y=3$ $\wwt_z=2$ && \centerline{$\Bbb{D}_5$}
\cr
\noalign{\hrule}
$T_{45}$ && $(6,2,0,0)$ && $-6$ &&
$\alpha_{12}x_1^3+\alpha_8x_1^2x_2+\alpha^1_6x_1^2x_3+
\alpha^2_6x_1^2x_4+\alpha_4x_1x_2^2+\alpha^1_2x_1x_2x_3+
\alpha^2_2x_1x_2x_4+\alpha^1_0x_1x_3^2+\alpha^2_0x_1x_3x_4+
\alpha^3_0x_1x_4^2+\alpha^4_0x_2^3=0$ && $xQ(x,y,z)+y^3=0$
($Q(0)\neq 0$) &&
\centerline{--} && Non-sin\-gu\-lar point
\cr
\noalign{\hrule}
$T_{46}$ && $(6,3,1,0)$ && $-8$ &&
$\alpha_{10}x_1^3+\alpha_7x_1^2x_2+\alpha_5x_1^2x_3+
\alpha^1_4x_1^2x_4+\alpha^2_4x_1x_2^2+
\alpha_2x_1x_2x_3+\alpha^1_1x_1x_2x_4+\alpha_0x_1x_3^2+\alpha^2_1x_2^3=0$
&& $x^3+x^2y+x^2z+x^2+xy^2+xyz+xy+xz^2+y^3{=}\,0$ &&
\centerline{--} && \centerline{$\Bbb{A}_5$}
\cr
\noalign{\hrule}
$T_{47}$ && $(6,3,2,0)$ && $-9$ &&
$\alpha_9x_1^3+\alpha_6x_1^2x_2+\alpha_5x_1^2x_3+
\alpha^1_3x_1^2x_4+\alpha^2_3x_1x_2^2+
\alpha_2x_1x_2x_3+\alpha^1_0x_1x_2x_4+\alpha_1x_1x_3^2+\alpha^2_0x_2^3=0$
&& $x^3+x^2y+x^2z+x^2+xy^2+xyz+xy+xz^2+y^3{=}\,0$ &&
\centerline{--} && \centerline{$\Bbb{A}_5$}
\cr
\noalign{\hrule}
$T_{48}$ && $(6,4,2,0)$ && $-10$ &&
$\alpha_8x_1^3+\alpha_6x_1^2x_2+\alpha^1_4x_1^2x_3+
\alpha^1_2x_1^2x_4+\alpha^2_4x_1x_2^2+\alpha^2_2x_1x_2x_3+
\alpha^1_0x_1x_2x_4+\alpha^2_0x_1x_3^2+\alpha^3_2x_2^3+
\alpha^3_0x_2^2x_3=\nomathbreak 0$ &&
$x^3+x^2y+x^2z+x^2+xy^2+xyz+xy+xz^2+y^3+y^2z=0$ &&
\centerline{--} && \centerline{$\Bbb{A}_4$}
\cr
\noalign{\hrule}
$T_{49}$ && $(6,4,3,0)$ && $-11$ &&
$\alpha_7x_1^3+\alpha_5x_1^2x_2+\alpha_4x_1^2x_3+
\alpha^1_1x_1^2x_4+\alpha_3x_1x_2^2+
\alpha_2x_1x_2x_3+\alpha^2_1x_1x_3^2+\alpha^3_1x_2^3+\alpha_0x_2^2x_3=0$
&& $x^3+x^2y+x^2z+x^2+xy^2+xyz+xz^2+y^3+y^2z=0$ &&
$\wwt_x=4$ $\wwt_y=3$ $\wwt_z=2$ && \centerline{$\Bbb{D}_5$}
\cr
\noalign{\hrule}
}}
$$
$$
\vbox{
\let\PLUS=+
\binoppenalty10000
\relpenalty10000
\catcode`\+=\active
\raggedright
\def+{\allowbreak\PLUS\nobreak\unskip}
\define\mystrut{\vphantom{$\Bigl($}}
\def\ensp{\kern0.39em}
\tabskip=0pt \offinterlineskip
\halign{\vrule\mystrut
\ensp\hfill#\hfill\ensp&\vrule# &
\ensp\hfill#\hfill\ensp&\vrule# &
\ensp\hfill#\hfill\ensp&\vrule# &
\ensp\hfill\tabboxone{#}\hfill\ensp&\vrule# &
\ensp\hfill\tabboxtwo{#}\hfill\ensp&\vrule# &
\ensp\hfill\tabboxthree{#}\hfill\ensp&\vrule# &
\ensp\hfill\tabboxthree{\hsize 1.2cm #}\hfill\ensp\vrule
\cr
\multispan{13}
\hfill{\smc Table 2 continued}
\cr
\noalign{\medskip}
\noalign{\hrule}
$T_{50}$ && $(6,4,3,1)$ && $-12$ &&
$\alpha_6x_1^3+\alpha_4x_1^2x_2+\alpha_3x_1^2x_3+
\alpha^1_1x_1^2x_4+\alpha_2x_1x_2^2+\alpha^2_1x_1x_2x_3+
\alpha^1_0x_1x_3^2+\alpha^2_0x_2^3=0$ &&
$x^3+x^2y+x^2z+x^2+xy^2+xyz+xz^2+y^3=0$ &&
$\wwt_x=6$ $\wwt_y=4$ $\wwt_z=3$ && \centerline{$\Bbb{E}_6$}
\cr
\noalign{\hrule}
$T_{51}$ && $(6,4,4,0)$ && $-12$ &&
$\alpha_6x_1^3+\alpha^1_4x_1^2x_2+\alpha^2_4x_1^2x_3+\alpha^1_0x_1^2x_4+
\alpha^1_2x_1x_2^2+\alpha^2_2x_1x_2x_3+\alpha^3_2x_1x_3^2+
\alpha^1_0x_2^3+\alpha^2_0x_2^2x_3+\alpha^3_0x_2x_3^2+\alpha^4_0x_3^3=0$
&& $x^2+P_3(y,z)+xP_2(y,z){=}\,0$ \vp\smash{($P_i$} is a homogeneous
polynomial of degree~$i$) &&
$\wwt_x=3$ $\wwt_y=2$ $\wwt_z=2$ && \centerline{$\Bbb{D}_4$}
\cr
\noalign{\hrule}
$T_{52}$ && $(6,5,3,0)$ && $-12$ &&
$\alpha_6x_1^3+\alpha_5x_1^2x_2+\alpha^1_3x_1^2x_3+
\alpha^1_0x_1^2x_4+\alpha_4x_1x_2^2+\alpha_2x_1x_2x_3+
\alpha^2_0x_1x_3^2+\alpha^2_3x_2^3+\alpha_1x_2^2x_3=0$ &&
$x^3+x^2y+x^2z+x^2+xy^2+xyz+xz^2+y^3+y^2z=0$ &&
$\wwt_x=4$ $\wwt_y=3$ $\wwt_z=2$ && \centerline{$\Bbb{D}_5$}
\cr
\noalign{\hrule}
$T_{53}$ && $(7,3,1,0)$ && $-9$ &&
$\alpha_{12}x_1^3+\alpha_8x_1^2x_2+\alpha_6x_1^2x_3+\alpha_5x_1^2x_4+
\alpha_4x_1x_2^2+\alpha_2x_1x_2x_3+\alpha_1x_1x_2x_4+
\alpha^1_0x_1x_3^2+\alpha^2_0x_2^3=0$ &&
$x^3+x^2y+x^2z+x^2+xy^2+xyz+xy+xz^2+y^3{=}\,0$ &&
\centerline{--} && \centerline{$\Bbb{A}_5$}
\cr
\noalign{\hrule}
$T_{54}$ && $(7,4,2,0)$ && $-11$ &&
$\alpha_{10}x_1^3+\alpha_7x_1^2x_2+\alpha_5x_1^2x_3+\alpha_3x_1^2x_4+
\alpha_4x_1x_2^2+\alpha_2x_1x_2x_3+\alpha^1_0x_1x_2x_4+\alpha^2_0x_1x_3^2+
\alpha_1x_2^3=0$ && $x^3+x^2y+x^2z+x^2+xy^2+xyz+xy+xz^2+y^3{=}\,0$ &&
\centerline{--} && \centerline{$\Bbb{A}_5$}
\cr
\noalign{\hrule}
$T_{55}$ && $(7,4,3,0)$ && $-12$ &&
$\alpha_9x_1^3+\alpha_6x_1^2x_2+\alpha_5x_1^2x_3+\alpha^1_2x_1^2x_4+
\alpha_3x_1x_2^2+\alpha^2_2x_1x_2x_3+\alpha_1x_1x_3^2+\alpha_0x_2^3=0$
&& $x^3+x^2y+x^2z+x^2+xy^2+xyz+xz^2+y^3=0$ &&
$\wwt_x=6$ $\wwt_y=4$ $\wwt_z=3$ && \centerline{$\Bbb{E}_6$}
\cr
\noalign{\hrule}
$T_{56}$ && $(7,5,3,0)$ && $-13$ &&
$\alpha_8x_1^3+\alpha_6x_1^2x_2+\alpha^1_4x_1^2x_3+\alpha^1_1x_1^2x_4+
\alpha^2_4x_1x_2^2+\alpha^1_2x_1x_2x_3+\alpha_0x_1x_3^2+
\alpha^2_2x_2^3+\alpha^2_1x_2^2x_3=0$ &&
$x^3+x^2y+x^2z+x^2+xy^2+xyz+xz^2+y^3+y^2z=0$ &&
$\wwt_x=4$ $\wwt_y=3$ $\wwt_z=2$ && \centerline{$\Bbb{D}_5$}
\cr
\noalign{\hrule}
}}
$$
$$
\vbox{
\let\PLUS=+
\binoppenalty10000
\relpenalty10000
\catcode`\+=\active
\raggedright
\def+{\allowbreak\PLUS\nobreak\unskip}
\define\mystrut{\vphantom{$\Bigl($}}
\def\ensp{\kern0.39em}
\tabskip=0pt \offinterlineskip
\halign{\vrule\mystrut
\ensp\hfill#\hfill\ensp&\vrule# &
\ensp\hfill#\hfill\ensp&\vrule# &
\ensp\hfill#\hfill\ensp&\vrule# &
\ensp\hfill\tabboxone{#}\hfill\ensp&\vrule# &
\ensp\hfill\tabboxtwo{#}\hfill\ensp&\vrule# &
\ensp\hfill\tabboxthree{#}\hfill\ensp&\vrule# &
\ensp\hfill\tabboxthree{\hsize 1.2cm #}\hfill\ensp\vrule
\cr
\multispan{13}
\hfill{\smc Table 2 continued}
\cr
\noalign{\medskip}
\noalign{\hrule}
$T_{57}$ && $(7,5,4,0)$ && $-14$ &&
$\alpha_7x_1^3+\alpha_5x_1^2x_2+\alpha_4x_1^2x_3+
\alpha^1_0x_1^2x_4+\alpha_3x_1x_2^2+\alpha_2x_1x_2x_3+
\alpha^1_1x_1x_3^2+\alpha^2_1x_2^3+\alpha^2_0x_2^2x_3=0$ &&
$x^3+x^2y+x^2z+x^2+xy^2+xyz+xz^2+y^3+y^2z=0$ &&
$\wwt_x=4$ $\wwt_y=3$ $\wwt_z=2$ && \centerline{$\Bbb{D}_5$}
\cr
\noalign{\hrule}
$T_{58}$ && $(7,5,4,1)$ && $-15$ &&
$\alpha_6x_1^3+\alpha_4x_1^2x_2+\alpha_3x_1^2x_3+\alpha^1_0x_1^2x_4+
\alpha_2x_1x_2^2+\alpha_1x_1x_2x_3+\alpha^2_0x_1x_3^2+\alpha^3_0x_2^3=0$
&& $x^3+x^2y+x^2z+x^2+xy^2+xyz+xz^2+y^3=0$ &&
$\wwt_x=6$ $\wwt_y=4$ $\wwt_z=3$ && \centerline{$\Bbb{E}_6$}
\cr
\noalign{\hrule}
$T_{59}$ && $(8,4,2,0)$ && $-12$ &&
$\alpha_{12}x_1^3+\alpha_8x_1^2x_2+\alpha_6x_1^2x_3+
\alpha^1_4x_1^2x_4+\alpha^2_4x_1x_2^2+\alpha^1_2x_1x_2x_3+
\alpha^2_2x_1x_2x_4+\alpha^1_0x_1x_3^2+\alpha^2_0x_2^3=0$ &&
$x^3+x^2y+x^2z+x^2+xy^2+xyz+xy+xz^2+y^3{=}\,0$ &&
\centerline{--} && \centerline{$\Bbb{A}_5$}
\cr
\noalign{\hrule}
$T_{60}$ && $(8,5,3,0)$ && $-14$ &&
$\alpha_{10}x_1^3+\alpha_7x_1^2x_2+\alpha_5x_1^2x_3+\alpha^1_2x_1^2x_4+
\alpha_4x_1x_2^2+\alpha^2_2x_1x_2x_3+\alpha_0x_1x_3^2+\alpha_1x_2^3=0$
&& $x^3+x^2y+x^2z+x^2+xy^2+xyz+xz^2+y^3=0$ &&
$\wwt_x=6$ $\wwt_y=4$ $\wwt_z=3$ && \centerline{$\Bbb{E}_6$}
\cr
\noalign{\hrule}
$T_{61}$ && $(8,5,4,0)$ && $-15$ &&
$\alpha_9x_1^3+\alpha_6x_1^2x_2+\alpha_5x_1^2x_3+\alpha^1_1x_1^2x_4+
\alpha_3x_1x_2^2+\alpha_2x_1x_2x_3+\alpha^2_1x_1x_3^2+\alpha_0x_2^3=0$
&& $x^3+x^2y+x^2z+x^2+xy^2+xyz+xz^2+y^3=0$ &&
$\wwt_x=6$ $\wwt_y=4$ $\wwt_z=3$ && \centerline{$\Bbb{E}_6$}
\cr
\noalign{\hrule}
$T_{62}$ && $(8,6,4,0)$ && $-16$ &&
$\alpha_8x_1^3+\alpha_6x_1^2x_2+\alpha^1_4x_1^2x_3+
\alpha^1_0x_1^2x_4+\alpha^2_4x_1x_2^2+\alpha^1_2x_1x_2x_3+
\alpha^2_0x_1x_3^2+\alpha^2_2x_2^3+\alpha^3_0x_2^2x_3=0$ &&
$x^3+x^2y+x^2z+x^2+xy^2+xyz+xz^2+y^3+y^2z=0$ &&
$\wwt_x=4$ $\wwt_y=3$ $\wwt_z=2$ && \centerline{$\Bbb{D}_5$}
\cr
\noalign{\hrule}
$T_{63}$ && $(9,5,3,0)$ && $-15$ &&
$\alpha_{12}x_1^3+\alpha_8x_1^2x_2+\alpha_6x_1^2x_3+\alpha_3x_1^2x_4+
\alpha_4x_1x_2^2+\alpha_2x_1x_2x_3+\alpha^1_0x_1x_3^2+\alpha^2_0x_2^3=0$
&& $x^3+x^2y+x^2z+x^2+xy^2+xyz+xz^2+y^3=0$ &&
$\wwt_x=6$ $\wwt_y=4$ $\wwt_z=3$ && \centerline{$\Bbb{E}_6$}
\cr
\noalign{\hrule}
}}
$$
$$
\vbox{
\let\PLUS=+
\binoppenalty10000
\relpenalty10000
\catcode`\+=\active
\raggedright
\def+{\allowbreak\PLUS\nobreak\unskip}
\define\mystrut{\vphantom{$\Bigl($}}
\def\ensp{\kern0.35em}
\tabskip=0pt \offinterlineskip
\halign{\vrule\mystrut
\ensp\hfill#\hfill\ensp&\vrule# &
\ensp\hfill#\hfill\ensp&\vrule# &
\ensp\hfill#\hfill\ensp&\vrule# &
\ensp\hfill\tabboxone{#}\hfill\ensp&\vrule# &
\ensp\hfill\tabboxtwo{#}\hfill\ensp&\vrule# &
\ensp\hfill\tabboxthree{#}\hfill\ensp&\vrule# &
\ensp\hfill\tabboxthree{\hsize 1.2cm #}\hfill\ensp\vrule
\cr
\multispan{13}
\hfill{\smc Table 2 continued}
\cr
\noalign{\medskip}
\noalign{\hrule}
$T_{64}$ && $(9,6,4,0)$ && $-17$ &&
$\alpha_{10}x_1^3+\alpha_7x_1^2x_2+\alpha_5x_1^2x_3+\alpha^1_1x_1^2x_4+
\alpha_4x_1x_2^2+\alpha_2x_1x_2x_3+\alpha_0x_1x_3^2+\alpha^2_1x_2^3=0$
&& $x^3+x^2y+x^2z+x^2+xy^2+xyz+xz^2+y^3=0$ &&
$\wwt_x=6$ $\wwt_y=4$ $\wwt_z=3$ && \centerline{$\Bbb{E}_6$}
\cr
\noalign{\hrule}
$T_{65}$ && $(9,6,5,0)$ && $-18$ &&
$\alpha_9x_1^3+\alpha_6x_1^2x_2+\alpha_5x_1^2x_3+\alpha^1_0x_1^2x_4+
\alpha_3x_1x_2^2+\alpha_2x_1x_2x_3+\alpha_1x_1x_3^2+\alpha^2_0x_2^3=0$
&& $x^3+x^2y+x^2z+x^2+xy^2+xyz+xz^2+y^3=0$ &&
$\wwt_x=6$ $\wwt_y=4$ $\wwt_z=3$ && \centerline{$\Bbb{E}_6$}
\cr
\noalign{\hrule}
$T_{66}$ && $(10,6,4,0)$ && $-18$ &&
$\alpha_{12}x_1^3+\alpha_8x_1^2x_2+\alpha_6x_1^2x_3+\alpha^1_2x_1^2x_4+
\alpha_4x_1x_2^2+\alpha^2_2x_1x_2x_3+\alpha^1_0x_1x_3^2+\alpha^2_0x_2^3=0$
&& $x^3+x^2y+x^2z+x^2+xy^2+xyz+xz^2+y^3=0$ &&
$\wwt_x=6$ $\wwt_y=4$ $\wwt_z=3$ && \centerline{$\Bbb{E}_6$}
\cr
\noalign{\hrule}
$T_{67}$ && $(10,7,5,0)$ && $-20$ &&
$\alpha_{10}x_1^3+\alpha_7x_1^2x_2+\alpha_5x_1^2x_3+\alpha^1_0x_1^2x_4+
\alpha_4x_1x_2^2+\alpha_2x_1x_2x_3+\alpha^2_0x_1x_3^2+\alpha_1x_2^3=0$
&& $x^3+x^2y+x^2z+x^2+xy^2+xyz+xz^2+y^3=0$ &&
$\wwt_x=6$ $\wwt_y=4$ $\wwt_z=3$ && \centerline{$\Bbb{E}_6$}
\cr
\noalign{\hrule}
$T_{68}$ && $(11,7,5,0)$ && $-21$ &&
$\alpha_{12}x_1^3+\alpha_8x_1^2x_2+\alpha_6x_1^2x_3+\alpha^1_1x_1^2x_4+
\alpha_4x_1x_2^2+\alpha_2x_1x_2x_3+\alpha_0x_1x_3^2+\alpha^2_1x_2^3=0$
&& $x^3+x^2y+x^2z+x^2+xy^2+xyz+xz^2+y^3=0$ &&
$\wwt_x=6$ $\wwt_y=4$ $\wwt_z=3$ && \centerline{$\Bbb{E}_6$}
\cr
\noalign{\hrule}
$T_{69}$ && $(12,8,6,0)$ && $-24$ &&
$\alpha_{12}x_1^3+\alpha_8x_1^2x_2+\alpha_6x_1^2x_3+\alpha^1_0x_1^2x_4+
\alpha_4x_1x_2^2+\alpha_2x_1x_2x_3+\alpha^2_0x_1x_3^2+\alpha^3_0x_2^3=0$
&& $x^3+x^2y+x^2z+x^2+xy^2+xyz+xz^2+y^3=0$ && $\wwt_x=6$ $\wwt_y=4$
$\wwt_z=3$ && \centerline{$\Bbb{E}_6$}
\cr
\noalign{\hrule}
}}
$$

}

Thus Theorem~1.6 is proved.

\remark{Remark 4.15}
The proof of Theorem~1.6 gives a description of possible singularities
of the 3-folds~$T_{j}$. For example, sufficiently general 3-folds $T_{j}$
are smooth for $j\in\{1,2,5,8,14\}$ and have only isolated ordinary
double points for $j\in\{4,7,11,16,19\}$. The smooth trigonal 3-folds
$T_{j}$ are well known (see~\cite{15}, \cite{88}). On the other hand,
the 3-fold $T_{j}$ always has non-isolated singularities for
$j\in\{6,9,13,15,17,21,24,25,26,28,32,\allowmathbreak
34,\allowmathbreak
36,40,41,44,45,50,58\}$.
In all other cases, the 3-fold $T_{j}$ has at least one
non-$\operatorname{cDV}$-point.
\endremark

\remark{Remark 4.16}
In the case of~$T_{15}$ the variety $V$ is always singular along the curve
$x_{1}=x_{2}=\alpha^1_1x_3^2+\alpha^2_1x_3x_4+\alpha^3_1x_4^2=0$.
In the case of~$T_{17}$ it is singular along the curve
$x_{1}=x_{2}=\alpha_0^1x_3^2+\alpha_0^2x_3x_4+\alpha_0^3x_4^2=\nomathbreak0$.
In the case of~$T_{19}$ it is singular along the curve
$x_{1}=x_{2}=\alpha^1_0x_3^2+\alpha^2_0x_3x_4+\alpha^3_0x_4^2=0$.
In the case of~$T_{25}$ it is singular along the curve
$x_{1}=x_{2}=\alpha^1_1x_3^2+\alpha^2_1x_3x_4+\alpha^3_1x_4^2=0$.
In the case of~$T_{26}$ it is singular along the curve
$x_{1}=x_{2}=\alpha^1_0x_3^2+\alpha^2_0x_3x_4+\alpha^3_0x_4^2=0$.
In the case of~$T_{28}$ it is singular along the curve
$x_{1}=x_{2}=\alpha^1_0x_3^2+\alpha^2_0x_3x_4+ \alpha^3_0x_4^2=0$.
In the case of~$T_{36}$ it is singular along the curve
$x_{1}=x_{2}=\alpha^1_0x_3^2+\alpha^2_0x_3x_4+\alpha^3_0x_4^2=0$.
In the case of~$T_{45}$ it is singular along the curve
$x_{1}=x_{2}=\alpha^1_0x_3^2+\alpha^2_0x_3x_4+\alpha^3_0x_4^2=0$.
All of these curves are bisections of the corresponding projections
$\varphi\:\FF(d_1, d_2, d_3, d_4)\to\PP^1$. They are reducible in the
cases of~$T_{17}$, $T_{19}$, $T_{26}$, $T_{28}$, $T_{36}$ and~$T_{45}$.
This simple observation will enable us to apply Lemma~5.2 to these
varieties and prove their rationality.
\endremark

\remark{Remark 4.17}
In the case of~$T_{7}$, the linear system $|M-L|$ determines a birational
map $\psi\:V\dasharrow\Bbb{P}^{3}$, which may be factorized as
$\psi=\omega\circ\gamma\circ\beta$. Here $\beta$ flops the curve~$C$,
$\gamma$ contracts the strict transform of the surface with equation
$x_{1}=0$ (on~$V$) onto a smooth rational curve whose image on
$\Bbb{P}^{3}$ is a line, and $\omega$ is a double covering of
$\Bbb{P}^{3}$ branched over a non-sin\-gu\-lar quartic surface.
In particular, $V$ is birationally isomorphic to a hypersurface of
degree~$4$ in~$\Bbb{P}(1^{4},2)$. The latter variety is also known as
a double space of index two. It was studied in
\cite{54}, \cite{31}, \cite{32}, \cite{37},
\cite{33}, \cite{60} and~\cite{34}.
\endremark

\remark{Remark 4.18}
In the cases of~$T_{4}$ and~$T_{6}$, the variety $X\subset \PP^6$ is an
anticanonically embedded Fano variety with canonical Gorenstein
singularities  and with $(-K_X)^3=\nomathbreak8$ 
(compare \cite{88}, Statement~4.1.12).
\endremark

\remark{Remark 4.19}
One can simplify the proof of Theorem~1.6 by arguing as follows. If
$X$ is a del~Pezzo surface of degree 3 over some field with a non-Du~Val
singular point defined over this field, then $X$ is a cone. The authors
did not use this approach by the reasons pointed out in Remark~3.12.
\endremark

\head
\S\,5. Rationality and non-rationality
\endhead

In this section we prove Proposition~1.10. Let $H_{i}$ and $T_{j}$ be the
Fano 3-folds from Theorems~1.5 and~1.6 respectively. The non-rationality
of sufficiently general 3-folds
$H_{1}$, $H_{2}$, $H_{3}$, $H_{4}$, $H_{6}$, $T_{1}$,
$T_{2}$, $T_{7}$, $T_{8}$ certainly follows (see Remark~1.8 and
Example~1.11) from the results of
\cite{22}, \cite{62}, \cite{54}, \cite{37},
\cite{16}, \cite{31}--\cite{33},
\cite{60}, \cite{38}, \cite{52}, \cite{34},
\cite{28}, \cite{29}, \cite{61}, \cite{111}, \cite{5},
\cite{82}, \cite{63}, \cite{98}, \cite{64},
\cite{6}, \cite{102}, \cite{7},~\cite{59}. On the other hand,
it is clear that the 3-folds $H_{9}$, $T_{5}$ and $T_{14}$ are always
rational (see Remark~1.8).

We may thus assume that $i\not\in\{1,2,3,4,6,9\}$ and
$j\not\in\{1,2,5,7,8,14\}$. Then the 3-fold $H_{i}$
is naturally birationally equivalent to a del~Pezzo fibration
$\tau\:Y_{i}\to\Bbb{P}^{1}$ of degree~$2$ (see Theorem~1.5) with
canonical Gorenstein singularities, and the 3-fold $T_{j}$
is naturally birationally equivalent to a del~Pezzo fibration
$\psi\:V_{j}\to\Bbb{P}^{1}$ of degree~$3$ (see Theorem~1.6) with
canonical Gorenstein singularities. Let $\overline{Y}_{i}$
and $\overline{V}_{j}$ be generic fibres of~$\tau$ and~$\psi$
respectively. Then $\overline{Y}_{i}$ and $\overline{V}_{j}$ are
del~Pezzo surfaces with Du~Val singularities defined over the
field $\Bbb{C}(x)$.

\remark{Remark 5.1}
Rationality of the surfaces $\overline{Y}_{i}$ and $\overline{V}_{j}$
over $\Bbb{C}(x)$ implies the rationality of the 3-folds $Y_{i}$ and $V_{j}$
respectively.
\endremark

\medskip
The del~Pezzo surfaces $\overline{Y}_{i}$ and $\overline{V}_{j}$ always have
a $\Bbb{C}(x)$-point by Theorem~2.24. Moreover, the sets of their
$\Bbb{C}(x)$-points are huge by Theorem~2.25.

\proclaim{Lemma 5.2}
Let $S$ be a del~Pezzo surface of degree~$3$ with canonical singularities
defined over an arbitrary perfect field $\Bbb{F}$. Suppose that the set
$\operatorname{Sing}(S)$ contains an $\Bbb{F}$-point $O\in S$.
Then $S$ is rational over~$\Bbb{F}$.
\endproclaim

\demo{Proof}
The surface $S$ is a cubic hypersurface in~$\Bbb{P}^{3}$ (see
\cite{23}, \cite{25}, \cite{20},~\cite{95}). Thus the projection from~$O$
gives a birational map to~$\Bbb{P}^{2}$.
\enddemo

Therefore the proof of Theorem~1.6 along with Lemma~5.2 immediately yields
the rationality of the 3-fold $T_{j}$ for
$j\in\{10, 12, 13, 17, \dots, 24, 26, \dots, 69\}$.

\proclaim{Lemma 5.3}
Let $S$ be a del~Pezzo surface of degree~$2$ with Du~Val singularities
defined over an arbitrary perfect field~$\Bbb{F}$.
Suppose that the singularity set $\operatorname{Sing}(S)$ contains an
$\Bbb{F}$-point $O\in S$ which is locally isomorphic to one of the following
Du~Val points\rom: $\Bbb{E}_{6}$, $\Bbb{E}_{7}$, $\Bbb{D}_{n}$ or
$\Bbb{A}_{k}$ for $n\geqslant 5$ and~$k\geqslant 7$. Then $S$ is rational
over~$\Bbb{F}$.
\endproclaim

\demo{Proof}
Let $f\:W\to S$ be a minimal resolution of singularities of~$S$, and let
$E=f^{-1}(O)\subset W$ be a connected curve defined over~$\Bbb F$. Then
$K_{W}\sim f^{*}(K_{S})$. In particular, $W$ is a weak del~Pezzo surface
(see~\cite{67}) of degree~$2$, the curve $E$ is
$\operatorname{Gal}(\overline{\Bbb{F}}/\Bbb{F})$-invariant, and all
irreducible components of~$E$ that are defined over $\overline{\Bbb{F}}$
must split into disjoint
$\operatorname{Gal}(\overline{\Bbb{F}}/\Bbb{F})$-orbits. However,
irreducible components of~$E$ form a graph of type
$\Bbb{E}_{6}$, $\Bbb{E}_{7}$,
$\Bbb{D}_{n}$ or $\Bbb{A}_{k}$ for $n\geqslant 5$ and
$k\geqslant 7$. Therefore the curve $E$ splits into at least $4$ (possibly
reducible) curves defined over $\Bbb{F}$. Since the intersection form
of irreducible  components of~$E$ is negative (see~\cite{51}), it follows
that the rank of~$\operatorname{Pic}(W)$ is at least~5.

There is a birational morphism $g\:W\to U$ defined over~$\Bbb{F}$
such that the surface  $U$ is minimal (see~\cite{14},
\cite{103}, \cite{20}), that is, no curve on~$U$ can be contracted
to a smooth point. Moreover, the rank of~$\operatorname{Pic}(U)$
does not exceed~$2$ by Theorem~2.22. Therefore
$K_{U}^{2}\geqslant K_{W}^{2}+3=5$. Thus the surface $U$ is rational
over $\Bbb{F}$ by Theorem~2.23.
\enddemo

Lemma~5.3 and the proof of Theorem~1.5 imply that the hyperelliptic 3-folds
$H_{i}$ are rational for $i\in\{22, 26, 27, 28, 29, 31,\dots, 47\}$.

\remark{Remark 5.4}
Non-rationality of the surfaces $\overline{Y}_{i}$ and $\overline{V}_{j}$
over $\Bbb{C}(x)$ does not imply non-rationality of the 3-folds
$H_{i}$ and $T_{j}$ respectively. However we believe that the rough
method used above can be also applied to prove non-rationality of~$H_{i}$
in many of the
remaining cases. For example, one can try to use the proofs of
Theorems~1.5 and~1.6 to describe the geometry
of the surfaces $\overline{Y}_{i}$ and $\overline{V}_{j}$ in more detail
and then use the results of
\cite{26}, \cite{35} and~\cite{20}.
\endremark

\proclaim{Proposition 5.5}
Let $X$ be a sufficiently general\,\footnote{Here and in what follows
we always understand ``general'' as ``belonging to a Zariski open subset
of the moduli space'' unless otherwise specified.}
Fano $3$-fold $T_{3}$ from Theorem~$1.6$. Then $X$ is non-rational.
\endproclaim

\demo{Proof}
Suppose that $U=\operatorname{Proj}(\Cal{O}_{\Bbb{P}^{2}}(2)\oplus
\Cal{O}_{\Bbb{P}^{2}}\oplus\Cal{O}_{\Bbb{P}^{2}})$, \
$f\:U\to\Bbb{P}^{2}$ is the natural projection, $T$ is the
tautological line bundle on~$U$, and $F=f^{*}(\Cal{O}_{\Bbb{P}^{2}}(1))$.
Then $X$ is an anticanonical image of a sufficiently general divisor
$V\in|2T+F|$. The 3-fold $V$ is smooth by the Bertini theorem. Moreover,
the Lefschetz theorem (see~\cite{55}, \cite{50}) implies that
$\operatorname{Pic}(V)\cong\Bbb{Z}\oplus\Bbb{Z}$.

Let $g\:V\to\Bbb{P}^{2}$ be the restriction of the projection
$f\: U\to\nomathbreak \Bbb{P}^{2}$. Then $g$ is a conic bundle.
Let $\Delta$ be the degeneration divisor of
$g$, and let $Y$ be a sufficiently general surface in the linear system
$|g^{*}(\Cal{O}_{\Bbb{P}^{2}}(1))|$. Then $Y$ is smooth and
$K_{Y}^{2}=1$ by the adjunction formula. Therefore the conic bundle
$g|_{Y}$ has 7 reducible fibres. Thus the degree of the divisor
$\Delta\subset\Bbb{P}^{2}$ is equal to~$7$ (see~\cite{56}, \S\,3.5),
and $V$ is non-rational by Theorem~2.16.\enddemo

\proclaim{Proposition 5.6}
Let $X$ be a sufficiently general\,\footnote{The complement to a
countable union of Zariski closed subset in the moduli space.}
$3$-fold $H_{5}$ from Theorem~$1.5$. Then $X$ is non-rational.
\endproclaim

\demo{Proof}
The 3-fold $X$ is an anticanonical model of a smooth weak Fano 3-fold
$V$, which may be described as a double covering
$\pi\:V\to U=\operatorname{Proj}(\Cal{O}_{\Bbb{P}^{1}}(2)\oplus
\Cal{O}_{\Bbb{P}^{1}}(1)\oplus\Cal{O}_{\Bbb{P}^{1}})$ branched over
a divisor $D\in |4M-2L|$, where $M$ is the tautological line bundle
on~$U$ and $L$ is a fibre of the natural projection of
$U$ to~$\Bbb{P}^{1}$. The divisor $D$ may be given in
the bihomogeneous coordinates (see Proposition~2.19) by the zeros of
the bihomogeneous polynomial
$$
\align
&\alpha_6x_1^4+\alpha_5x_1^3x_2+\alpha^1_4x_1^3x_3+\alpha^2_4x_1^2x_2^2
+\alpha^1_3x_1^2x_2x_3+\alpha^1_2x_1^2x_3^2+\alpha^2_3x_1x_2^3
\\
&\qquad\qquad
+\alpha^2_2x_1x_2^2x_3+\alpha^1_1x_1x_2x_3^2+\alpha^1_0x_1x_3^3
+\alpha^3_2x_2^4+\alpha^2_1x_2^3x_3+\alpha^2_0x_2^2x_3^2,
\endalign
$$
where $\alpha_d^{i}=\alpha_d^{i}(t_1,t_2)$ is a homogeneous polynomial of
degree~$d$.

Consider a double covering $\chi\:Y\to U$ branched over a sufficiently
general divisor $\Delta\subset U$ which is given by the zeros of the same
bihomogeneous polynomial as~$D$ with the only exception
that $\alpha_0^{1}=0$. Then $Y$ is not smooth because $\Delta$ has
singularities along the curve $Y_{3}\subset U$ given by~$x_1=x_2=0$.
The curve $Y_{3}$ is the smallest negative subscroll of~$U$ (see
Proposition~2.19). We may assume that
$\Delta\subset U$ is a sufficiently general element of the
linear subsystem in the system $|4M-2L|$ consisting of all
divisors with singularities along $Y_{3}$. The divisor
$\Delta$ is smooth outside $Y_{3}$ by the Bertini theorem.

Put $C=\chi^{-1}(Y_{3})$. Then the 3-fold $Y$ has singularities of type
$\Bbb{A}_{1}\times\Bbb{C}$ at the general point of the curve~$C$.
Moreover, the singularities of~$Y$ at other points of~$C$ are locally
isomorphic to the singularity
$$
x^{2}+y^{2}+z^{2}t=0\subset\Bbb{C}^{4}\cong
\operatorname{Spec}(\Bbb{C}[x,y,z,t]),
$$
where the curve $C$ is locally given by~$x=y=z=0$. It follows that one can
resolve the singularities of~$Y$ by one blow up
$f\:\widetilde{Y}\to\nomathbreak Y$ of the curve~$C$.

Let $g\:\widetilde{U}\to U$ be the blow up of the curve $Y_{3}\subset U$.
Then the diagram
$$
\xymatrix{
\wtilde{Y}\ar@{->}[d]_{f}\ar@{->}[rr]^{\wtilde{\chi}}&&
\wtilde{U}\ar[d]^{g}
\\
Y\ar@{->}[rr]^{\chi}&&U}
$$
is commutative, where $\widetilde{\chi}\:\widetilde{Y}\to\widetilde{U}$
is a double covering. Let $E$ be the exceptional divisor of~$g$. Then
$\widetilde{\chi}$ is branched over the divisor
$g^{-1}(\Delta)\sim g^{*}(4M-2L)-2E$.

In the case when the divisor $g^{-1}(\Delta)$ is ample 
on~$\widetilde{U}$, the Lefschetz theorem (see~\cite{55}, \cite{50},
\cite{121}) implies that
$\operatorname{Pic}(\,\widetilde{Y})\cong
\operatorname{Pic}(\,\widetilde{U})\cong\Bbb{Z}^{3}$
(see \cite{60}, \cite{65},~\cite{66}). However, the divisor
$g^{-1}(\Delta)$ is not ample, although it is numerically effective
and big. Indeed, the linear system $|g^{*}(M-L)-E|$ is free
and the linear system $|g^{*}(M)-E|$ gives a $\Bbb{P}^{1}$-bundle
$$
\tau\:\widetilde{U}\to\operatorname{Proj}(\Cal{O}_{\Bbb{P}^{1}}(2)\oplus
\Cal{O}_{\Bbb{P}^{1}}(1))\cong\Bbb{F}_{1}.
$$
Therefore the divisor $g^{-1}(\Delta)\sim g^{*}(4M-2L)-2E$ is numerically
effective and big. Hence we can replace the Lefschetz theorem by the first
part of the proof of Proposition~32 in~\cite{56} to get
$\operatorname{Pic}(\,\widetilde{Y})\cong\operatorname{Pic}(\,\widetilde{U})
\cong\Bbb{Z}^{3}$.

Let $Y_{2}\subset U$ be the largest negative subscroll
(see Proposition~2.19). The surface $Y_{2}$ is given by the equation
$x_{1}=0$ in the bihomogeneous coordinates on~$U$. Moreover,
$Y_{2}\cong\operatorname{Proj}\bigl(\Cal{O}_{\Bbb{P}^{1}}(1)\oplus
\Cal{O}_{\Bbb{P}^{1}}\bigr)$.
Put $S=g^{-1}(Y_{2})$. Then $S\cong Y_{2}$ and the morphism $\tau$
contracts the surface $S$ to the exceptional section of~$\Bbb{F}_{1}$.

By construction, the $\Bbb{P}^{1}$-bundle $\tau$ induces a conic bundle
$\widetilde{\tau}=\tau\circ\widetilde{\chi}$: $\widetilde{Y}\to\Bbb{F}_{1}$.
Put $\widetilde{S}=\widetilde{\chi}^{-1}(S)$, and let
$Z\subset\widetilde{Y}$ be a general fibre of the natural projection of
$\widetilde{Y}$ to~$\Bbb{P}^{1}$. Then $Z$ is a smooth weak del~Pezzo
surface of degree~$2$, that is, $-K_{Z}$ is numerically effective and big
and~$K_{Z}^{2}=2$. Moreover, the morphism
$g\circ\widetilde{\chi}|_{\widetilde{S}}\:\widetilde{S}\to Y_{2}$ is a
double covering branched over a divisor with the following equation in the
bihomogeneous coordinates:
$$
\alpha^3_2(t_{0},t_{1})\,x_2^2+\alpha^2_1(t_{0},t_{1})\,x_2x_3
+\alpha^2_0(t_{0},t_{1})\,x_3^2=0,
$$
where $\alpha_d^{i}(t_1,t_2)$ is the homogeneous polynomial of degree~$d$
from the bihomogeneous equation of~$\Delta$.

Let $\Xi\subset\Bbb{F}_{1}$ be the degeneration divisor of the conic bundle
$\widetilde{\tau}$. Then $\Xi\sim 6s_{\infty}+al$, where
$s_{\infty}$ is the exceptional section of~$\Bbb{F}_{1}$, \  $l$
is a fibre of the projection of~$\Bbb{F}_{1}$ to~$\Bbb{P}^{1}$,
and $a\in\Bbb{Z}$. The structure of the morphism
$g\circ\widetilde{\chi}|_{\widetilde{S}}$ implies that
$s_{\infty}\not\subset\Xi$. Moreover, the intersection
$s_{\infty}\cdot\Xi$ is equal to the number of reducible fibres
of the induced conic bundle $\widetilde{\tau}|_{\widetilde{S}}$.
This number can easily be calculated from the bihomogeneous equation
of the ramification divisor of~$g\circ\widetilde{\chi}|_{\widetilde{S}}$.
More precisely, reducible fibres of~$\widetilde{\tau}|_{\widetilde{S}}$
correspond to zeros of the discriminant
$(\alpha^2_1)^{2}-4\alpha^2_0\alpha^3_2$, whence $s_{\infty}\cdot\Xi=2$.
Therefore $a=8$. Thus $Y$ is non-rational by Theorem~2.16.

The 3-fold $Y$ is rationally connected (see~\cite{95}).
Thus the non-rationality of~$Y$ implies that $Y$ is non-ruled as well.
Therefore the 3-fold $V$ is non-ruled by Theorem~2.18 because we assumed
$V$ to be sufficiently general. Hence $X$ is non-rational.
\enddemo

\proclaim{Proposition 5.7}
Let $X$ be a sufficiently general $3$-fold $H_{7}$ from Theorem~$1.5$.
Then $X$ is non-rational.
\endproclaim

\demo{Proof}
The 3-fold $X$ is an anticanonical model of a weak Fano 3-fold
$V$ such that there is a double covering
$$
\pi\:V\to U=\operatorname{Proj}\bigl(\Cal{O}_{\Bbb{P}^{1}}(2)\oplus
\Cal{O}_{\Bbb{P}^{1}}(2)\oplus\Cal{O}_{\Bbb{P}^{1}}\bigr),
$$
branched over a divisor $D\in |4M-4L|$, where $M$ is the tautological line
bundle on~$U$ and $L$ is a fibre of the natural projection of~$U$ to
$\Bbb{P}^{1}$. The divisor $D$ may be given in the bihomogeneous coordinates
(see Proposition~2.19) by the zeros of a bihomogeneous polynomial
$$
\align
&\alpha^1_4x_1^4+\alpha^2_4x_1^3x_2+\alpha^3_4x_1^2x_2^2+\alpha^4_4x_1x_2^3
+\alpha^5_4x_2^4+\alpha^1_2x_1^3x_3+\alpha^2_2x_1^2x_2x_3
\\
&\qquad\qquad
+\alpha^3_2x_1x_2^2x_3+\alpha^4_2x_2^3x_3+\alpha^1_0x_1^2x_3^2
+\alpha^2_0x_1x_2x_3^2+\alpha^3_0x_2^2x_3^2,
\endalign
$$
where $\alpha_d^{i}=\alpha_d^{i}(t_1,t_2)$ is a homogeneous polynomial of
degree~$d$.

The divisor $D$ has singularities along the curve $Y_{3}\subset U$
given by~$x_1=x_2=0$. Since $X$ is general, the divisor $D\subset U$
is a sufficiently general element of the linear system $|4M-4L|$.
In particular, $D$ is smooth outside $Y_{3}$ by the Bertini theorem.
The 3-fold $V$ has singularities of the type $\Bbb{A}_{1}\times\Bbb{C}$
at a general point of the curve $C=\chi^{-1}(Y_{3})$. Moreover, one can
resolve the singularities of~$V$ by one blow up $f\:\widetilde{V}\to V$
of the curve~$C$.

Let $g\:\widetilde{U}\to U$ be the blow up of the curve $Y_{3}\subset U$.
Then the diagram
$$
\xymatrix{
\wtilde{V}\ar@{->}[d]_{f}\ar@{->}[rr]^{\wtilde{\pi}}&&
\wtilde{U}\ar[d]^{g}
\\
V\ar@{->}[rr]^{\pi}&&U}
$$
is commutative, where the morphism
$\widetilde{\pi}\:\widetilde{V}\to\widetilde{U}$ is a double covering.
Let $E$ be the exceptional divisor of~$g$. Then $\widetilde{\pi}$ is
branched over the divisor $g^{-1}(D)\sim g^{*}(4M-4L)-2E$. On the other
hand, the linear system $|g^{*}(M-2L)-E|$ is a free pencil whose image
on~$U$ is generated by the divisors $x_{1}=0$ and~$x_{2}=0$.
In particular, the divisor $g^{-1}(D)\sim g^{*}(4M-4L)-2E$
is numerically effective and big on~$\widetilde{U}$.
Then the first part of the proof of Proposition~32 in~\cite{56} (a
stronger version of the Lefschetz theorem) implies that
$$
\operatorname{Pic}(\,\widetilde{V})\cong
\operatorname{Pic}(\,\widetilde{U})\cong\Bbb{Z}^{3}.
$$

The linear system $|g^{*}(M-L)-E|$ is also free and gives a
$\Bbb{P}^{1}$-bundle
$$
\tau\:\widetilde{U}\to\operatorname{Proj}\bigl(\Cal{O}_{\Bbb{P}^{1}}(2)\oplus
\Cal{O}_{\Bbb{P}^{1}}(2)\bigr)\cong\Bbb{F}_{0}.
$$
The rational map
$\tau\circ g^{-1}$ is given in the bihomogeneous coordinates by the
linear system on~$U$ spanned by
$\beta_{1}(t_{0},t_{2})\,x_{1}+\allowmathbreak
\beta_{2}(t_{0},t_{2})\,x_{2}$, where
$\beta_{i}(t_{0},t_{2})$ is a homogeneous polynomial of degree~$1$.

The $\Bbb{P}^{1}$-bundle $\tau$ induces a conic bundle $\widetilde{\tau}=
\tau\circ\widetilde{\pi}\:\widetilde{V}\to\Bbb{F}_{0}$.
Let $\Delta\subset\Bbb{F}_{0}$ be the degeneration divisor of
$\widetilde{\tau}$, and let $L_{1}$,~$L_{2}$ be fibres of the two
projections of~$\Bbb{F}_{0}$ to~$\Bbb{P}^{1}$ such that
$\tau^{*}(L_{1})\sim g^{*}(L)$ and $\tau^{*}(L_{2})\sim g^{*}(M-2L)-E$.
Then $\Delta\sim nL_{1}+6L_{2}$ for some $n\in\Bbb{Z}$. Moreover, we have
$n=4$ by elementary calculations (see the proof of Proposition~5.6).
Hence $V$ is non-rational by Theorem~2.16.
\enddemo

\proclaim{Proposition 5.8}
Let $X$ be a $3$-fold $H_{8}$ from Theorem~$1.5$. Then $X$ is rational.
\endproclaim

\demo{Proof} \
Arguing as in the proof of Proposition~5.7, we get a conic bundle
$\widetilde{\tau}=\tau\circ\widetilde{\pi}\:
\widetilde{V}\to\Bbb{P}^{1}\times\Bbb{P}^{1}$,
where $\widetilde{V}$ is birationally isomorphic to~$X$.
Moreover, this case is simpler since the proof of rationality of
$\widetilde{V}$ does not require to prove that the conic bundle
$\widetilde{\tau}$ is standard, that is, that
$\operatorname{Pic}(\,\widetilde{V})\cong\Bbb{Z}^{3}$. Simple
calculations show that the degeneration divisor
$\Delta\subset\Bbb{P}^{1}\times\Bbb{P}^{1}$ of~$\widetilde{\tau}$ has
bidegree $(6,2)$. Now the rationality of~$X$ follows immediately from
Theorems~2.24, 2.22 and~2.23.
\enddemo

\proclaim{Proposition 5.9}
Let $X$ be a sufficiently general $3$-fold $T_{4}$ from Theorem~$1.6$.
Then $X$ is non-rational.
\endproclaim

\demo{Proof}
The 3-fold $X$ is an anticanonical image of a sufficiently general divisor
$$
V\subset
U=\operatorname{Proj}\bigl(\Cal{O}_{\Bbb{P}^{1}}(1)\oplus
\Cal{O}_{\Bbb{P}^{1}}(1)\oplus\Cal{O}_{\Bbb{P}^{1}}(1)\oplus
\Cal{O}_{\Bbb{P}^{1}}\bigr)
$$
belonging to the linear system $|3M-L|$, where $M$ is the tautological
line bundle on~$U$, and $L$ is a fibre of the natural projection of~$U$
to~$\Bbb{P}^{1}$. The divisor $D$ is given in the bihomogeneous coordinates
on~$U$ by
$$
\align
&\alpha^1_2x_1^3+\alpha^2_2x_1^2x_2+\alpha^3_2x_1^2x_3
+\alpha^4_2x_1^2x_2+\alpha^5_2x_1x_2x_3+\alpha^6_2x_1x_3^2
+\alpha^7_2x_2^3
\\
&\qquad\qquad
+\alpha^8_2x_2^2x_3+\alpha^9_2x_2x_3^2
+\alpha^{10}_2x_3^3+\alpha^1_1x_1^2x_4+\alpha^2_1x_1x_2x_4
+\alpha^3_1x_1x_3x_4
\\
&\qquad\qquad
+\alpha^4_1x_2^2x_4+\alpha^5_1x_2x_3x_4
+\alpha^5_1x_3^2x_4+\alpha^1_0x_1x_4^2+\alpha^2_0x_2x_4^2
+\alpha^3_0x_3x_4^2=0,
\endalign
$$
where $\alpha^{i}_{d}=\alpha^{i}_{d}(t_{0},t_{1})$
is a homogeneous polynomial of degree~$d$. Since $X$ is general, $V$ is
smooth. Moreover, the anticanonical morphism
$\varphi_{|-K_{V}|}$ contracts a single curve $C\subset V$
(given by~$x_{1}=x_{2}=x_{3}=0$) to an ordinary double point
$O$ on~$X$. The corresponding birational morphism $\varphi_{|M|}$ maps
the rational scroll $U$ to the cone $\overline{U}$ over
$\Bbb{P}^{1}\times\Bbb{P}^{2}$ with vertex~$O$.

The 3-fold $X$ and the 4-fold $\overline{U}$ are not $\Bbb{Q}$\kkk-factorial.
Moreover, the birational morphisms $\varphi_{|-K_{V}|}$ and $\varphi_{|M|}$
may be regarded as~$\Bbb{Q}$\kkk-factorializations of~$X$ and~$\overline{U}$
respectively (see~\cite{92}). There are also other ways to
$\Bbb{Q}$\kkk-factorialize $X$ and~$\overline{U}$. Namely, one can find a
scroll
$\widetilde{U}=\operatorname{Proj}\bigl(\Cal{O}_{\Bbb{P}^{2}}(1)\oplus
\Cal{O}_{\Bbb{P}^{2}}(1)\oplus\Cal{O}_{\Bbb{P}^{2}}\bigr)$
and a birational morphism
$\varphi_{|T|}\:\widetilde{U}\to\overline{U}$, where $T$ is the tautological
line bundle on~$\widetilde{U}$. Moreover, the birational map
$\varphi_{|T|}^{-1}\circ\varphi_{|M|}$ is an antiflip (see~\cite{93},
\cite{99}) in the  curve $C\subset U$.

Let $Y\subset\widetilde{U}$ be the proper transform of~$X$ on the 4-fold
$\widetilde{U}$. Then $Y$ is a smooth weak Fano 3-fold and
$Y\sim 2T+F$ for $F=f^{*}(\Cal{O}_{\Bbb{P}^{2}}(1))$, where $f$ is
the natural projection of~$\widetilde{U}$ to~$\Bbb{P}^{2}$.
The original 3-fold $X$ is an anticanonical image of the 3-fold~$Y$, and
the birational map $\varphi_{|-K_{Y}|}^{-1}\circ\varphi_{|-K_{V}|}$ is a
simple flop in the curve $C\subset V$ induced by the antiflip
$\varphi_{|T|}^{-1}\circ\nomathbreak\varphi_{|M|}$. The Lefschetz theorem
implies that $\operatorname{Pic}(Y)\cong\Bbb{Z}\oplus\Bbb{Z}$.

The restriction $g\:Y\to\Bbb{P}^{2}$ of the projection
$f\:\widetilde{U}\to\Bbb{P}^{2}$ is a conic bundle. Let $\Delta$
be the degeneration divisor of~$g$. Simple calculations (see the
proof of Proposition~5.5) imply that $\Delta\sim\Cal{O}_{\Bbb{P}^{2}}(7)$
(see~\cite{56}, \S\,4.4.1). Therefore the 3-fold $Y$ is non-rational by
Theorem~2.16 (see~\cite{54}, \cite{37}).
\enddemo

\proclaim{Proposition 5.10}
Let $X$ be a sufficiently general $3$-fold $T_{6}$ from Theorem~$1.6$.
Then $X$ is non-rational.
\endproclaim

\demo{Proof}
The 3-fold $X$ is an anticanonical image of a weak Fano 3-fold~$V$,
which may be regarded as a sufficiently general divisor on the rational
scroll
$U=\operatorname{Proj}(\Cal{O}_{\Bbb{P}^{1}}(2)\oplus
\Cal{O}_{\Bbb{P}^{1}}(1)\oplus\Cal{O}_{\Bbb{P}^{1}}\oplus
\Cal{O}_{\Bbb{P}^{1}})$
lying in the linear system $|3M-L|$, where $M$ is the tautological line
bundle on~$U$ and $L$ is a fibre of the natural projection of~$U$
to~$\Bbb{P}^{1}$. Thus the 3-fold $V$ is given in the bihomogeneous
coordinates on~$U$ by
$$
\align
&\alpha_5x_1^3+\alpha_4x_1^2x_2+\alpha^1_3x_1^2x_3
+\alpha^2_3x_1^2x_4+\alpha^3_3x_1x_2^2+\alpha^1_2x_1x_2x_3
\\
&\qquad\qquad
+\alpha^2_2x_1x_2x_4+\alpha^1_1x_1x_3^2
+\alpha^2_1x_1x_3x_4+\alpha^3_1x_1x_4^2+\alpha^3_2x_2^3
\\
&\qquad\qquad
+\alpha^4_1x_2^2x_3+\alpha^5_1x_2^2x_4+\alpha^1_0x_2x_3^2
+\alpha^2_0x_2x_3x_4+\alpha^3_0x_2x_4^2=0,
\endalign
$$
where $\alpha^{i}_{d}=\alpha^{i}_{d}(t_{0},t_{1})$ is a homogeneous
polynomial of degree~$d$. The 3-fold $V$ contains a surface
$Y_{3}\cong\Bbb{P}^{1}\times\Bbb{P}^{1}$ given by~$x_{1}=x_{2}=0$.
This surface is the base locus of the linear system~$|3M-L|$. However,
$V$ is smooth at the general point of~$Y_{3}$. On the  other hand,
$V$ is always singular at the points where
$$
x_{1}=x_{2}=\alpha^1_1x_3^2+\alpha^2_1x_3x_4+\alpha^3_1x_4^2
=\alpha^1_0x_3^2+\alpha^2_0x_3x_4+\alpha^3_0x_4^2=0.
$$
Since $V$ is general, it follows that these points are ordinary double
points on~$V$ and $V$ is smooth outside them.

Let $g\:\widetilde{U}\to U$ be the blow up of~$Y_{3}\subset U$, \
$E$ the exceptional divisor of~$g$, and
$\widetilde{V}=g^{-1}(V)\subset\widetilde{U}$. Then
$\widetilde{V}\sim g^{*}(3M-L)-E$, \ $\widetilde{V}$ is smooth, and
$g|_{\widetilde{V}}$ is a small resolution of the 3-fold~$V$. On the
other hand, the linear system $|g^{*}(M-L)-E|$ is free. Therefore the
divisor $\widetilde{V}$ is numerically effective and big on
$\widetilde{U}$. In the case when $\widetilde{V}$ is ample, the Lefschetz
theorem implies that $\operatorname{Pic}(\,\widetilde{V})\cong
\operatorname{Pic}(\,\widetilde{U})\cong\Bbb{Z}^{3}$.
However, the divisor $\widetilde{V}$ is not ample. Nevertheless, we can
replace the Lefschetz theorem by the arguments of the first part of the
proof of Proposition~32 in~\cite{56} to get
$\operatorname{Pic}(\,\widetilde{V})\cong
\operatorname{Pic}(\,\widetilde{U})\cong\Bbb{Z}^{3}$.

The linear system $|g^{*}(M)-E|$ is free and determines a
$\Bbb{P}^{2}$-bundle
$$
\tau\:\widetilde{U}\to\operatorname{Proj}\bigl(\Cal{O}_{\Bbb{P}^{1}}(2)\oplus
\Cal{O}_{\Bbb{P}^{1}}(1)\bigr)\cong\Bbb{F}_{1}.
$$
The rational map $\tau\circ g^{-1}$ is given in the bihomogeneous coordinates
by a linear system on~$U$ spanned by
$\beta_{1}(t_{0},t_{2})\,x_{1}+\allowmathbreak
\beta_{2}(t_{0},t_{2})\,x_{2}$, where $\beta_{i}(t_{0},t_{2})$
is a homogeneous polynomial of degree~$1$.

The $\Bbb{P}^{2}$-bundle $\tau$ induces a conic bundle
$\widetilde{\tau}=\tau|_{\widetilde{V}}\:\widetilde{V}\to\Bbb{F}_{1}$.
Let $\Delta\subset\Bbb{F}_{1}$ be the degeneration divisor of
$\widetilde{\tau}$. We see from the construction that
$\Delta\sim 5s_{\infty}+al$, where $s_{\infty}$ is the exceptional
section of~$\Bbb{F}_{1}$ and $l$ is a fibre of the natural projection
of $\Bbb{F}_{1}$ to~$\Bbb{P}^{1}$.

Let $s_{0}$ be a sufficiently general section of~$\Bbb{F}_{1}$ such that
$s_{0}\cap s_{\infty}=\varnothing$. We put $S=\widetilde{\tau}^{-1}(s_{0})$
and $B=\tau^{-1}(s_{0})$. Then $S=B\cap\widetilde{V}\subset B$, and the
divisor~$B$ is naturally isomorphic to the scroll
$$
\operatorname{Proj}\bigl(\Cal{O}_{\Bbb{P}^{1}}(2)\oplus
\Cal{O}_{\Bbb{P}^{1}}\oplus\Cal{O}_{\Bbb{P}^{1}}\bigr).
$$
Moreover, $g(B)\cong B$ and $g(B)\cap V=g(S)\cup Y_{3}$. However,
the surface $Y_{3}$ is determined by the equation $x_{1}=0$ on the scroll
$g(B)$, while $g(B)$ is a general divisor in the linear system $|M-L|$.
Therefore we have $S\sim 2T+F$ on the scroll~$B$, where $T$ is the
 tautological line bundle on~$B$ and $F$ is a fibre of the natural projection
of $B$ to~$\Bbb{P}^{1}$. It follows that $K_{S}^{2}=1$, $s_{0}\cdot\Delta=7$
and~$a=7$. Hence $\widetilde{V}$ is non-rational by Theorem~2.16.
\enddemo

\proclaim{Proposition 5.11}
Let $X$ be a $3$-fold $T_{25}$ from Theorem~$1.6$.
Then $X$ is rational.
\endproclaim

\demo{Proof}
We can repeat the construction of the conic bundle in Proposition~5.10
to get a conic bundle
$\widetilde{\tau}=\widetilde{V}\to\Bbb{F}_{3}$, where $\widetilde{V}$
is birationally equivalent to~$X$. However, we need not the condition
$\operatorname{Pic}(\,\widetilde{V})\cong\Bbb{Z}^{3}$ and smoothness
of~$\widetilde{V}$. Let $\Delta\subset\Bbb{F}_{3}$
be the degeneration divisor of~$\widetilde{\tau}$. Then elementary
calculations imply that $\Delta\cdot s_{0}=1$, where $s_{0}$ is a
sufficiently general section on~$\Bbb{F}_{3}$ which is disjoint from the
exceptional section of~$\Bbb{F}_{3}$. Therefore the 3-fold $\widetilde{V}$
is rational by Theorems~2.24, 2.22 and~2.23.
\enddemo

\proclaim{Proposition 5.12}
Let $X$ be a sufficiently general $3$-fold $T_{9}$ from Theorem~$1.6$.
Then $X$ is non-rational.
\endproclaim

\demo{Proof}
We can repeat the construction of the conic bundle in
the proof of Proposition~5.10 to get a conic bundle
$\widetilde{\tau}=\widetilde{V}\to\Bbb{F}_{0}\cong
\Bbb{P}^{1}\times\Bbb{P}^{1}$ such that $\widetilde{V}$ is
birationally equivalent to~$X$, \ $\widetilde{V}$ is smooth, and
$\operatorname{Pic}(\,\widetilde{V})\cong\Bbb{Z}^{3}$. Let
$\Delta\subset\Bbb{F}_{0}$ be the degeneration divisor of
$\widetilde{\tau}$. Then elementary calculations (see the proof of
Proposition~5.10) imply that the divisor
$\Delta\subset\Bbb{P}^{1}\times\Bbb{P}^{1}$ has bidegree $(5,4)$.
Therefore $\widetilde{V}$ is non-rational by Theorem~2.16.
\enddemo

\proclaim{Proposition 5.13}
Let $X$ be a $3$-fold $T_{11}$ from Theorem~$1.6$. Then
$X$ is rational.
\endproclaim

\demo{Proof}
We can repeat the construction of the conic bundle in
the proof of Proposition~5.10 to get a conic bundle
$\widetilde{\tau}=\widetilde{V}\to\Bbb{F}_{0}\cong
\Bbb{P}^{1}\times\Bbb{P}^{1}$ such that $\widetilde{V}$ is
birationally isomorphic to~$X$. However we do not need to prove that
$\widetilde{V}$ is smooth and
$\operatorname{Pic}(\,\widetilde{V})\cong\Bbb{Z}^{3}$. Let
$\Delta\subset\Bbb{F}_{0}$ be the degeneration divisor of~$\widetilde{\tau}$.
Then elementary calculations imply that the divisor
$\Delta\subset\Bbb{P}^{1}\times\Bbb{P}^{1}$ has bidegree
$(5,2)$. Hence we can consider the composite
$\theta\:\widetilde{V}\to\Bbb{P}^{1}$ of the conic bundle
$\widetilde{\tau}$ and one of the projections of
$\Bbb{P}^{1}\times\Bbb{P}^{1}$ onto~$\Bbb{P}^{1}$ such that a sufficiently
general fibre of~$\theta$ is a surface~$S$ with $K_{S}^{2}=6$. Then the
rationality of~$\widetilde{V}$ follows from Theorems~2.24, 2.22 and~2.23.
\enddemo

Thus Proposition~1.10 is proved. The approach to proving the non-rationality
of $H_{i}$ and~$T_{j}$ together with the standard degeneration technique
(see~\cite{54}, \cite{37}, \cite{94}) can be used as a pattern to prove
non-rationality of many 3-folds fibred into del~Pezzo surfaces of degree
$2$ and~$3$ (see~\cite{52}, \cite{30}, \cite{56},~\cite{43}).

\enddocument